\documentclass{article}[12pt]

\usepackage{amssymb}
\usepackage{amsmath}
\usepackage{enumerate}
\usepackage{latexsym}
\usepackage{mathrsfs}
\usepackage{amsthm}
\usepackage{verbatim}
\usepackage{graphicx}
\usepackage{epstopdf}
\usepackage{epsfig}

\usepackage{color}
\usepackage[dvipsnames]{xcolor}

\usepackage{subfig}
\usepackage{float}
\usepackage{titlesec}
\usepackage{bm}
\usepackage[colorlinks,linkcolor=blue,anchorcolor=green,citecolor=red]{hyperref}
\usepackage[T1]{fontenc}

\usepackage{charter}
\usepackage{pythonhighlight}
\usepackage{fancyhdr}
\usepackage{amscd}
\usepackage{cleveref}

\usepackage{caption}
\usepackage{algorithm} 
\usepackage{algorithmic} 

\newtheorem{theorem}{Theorem}[section]
\newtheorem{corollary}[theorem]{Corollary}
\newtheorem{lemma}[theorem]{Lemma}
\newtheorem{proposition}[theorem]{Proposition}

\newtheorem{remark}[theorem]{Remark}
\newtheorem{definition}[theorem]{Definition}

\usepackage{apptools}

\makeatletter
\@addtoreset{equation}{section}
\makeatother

\makeatletter 
\@addtoreset{equation}{section}
\makeatother  
\renewcommand\theequation{\oldstylenums{\thesection}%
	.\oldstylenums{\arabic{equation}}}

\usepackage{geometry}
\usepackage{listings}
\usepackage{setspace}

\newcommand{\ssubset}{\subset\joinrel\subset}

\begin{document}
	
	\title{A nonlinear homogenization-based perspective on the soft modes and effective energies of some conformal metamaterials}
	\author{Xuenan Li\thanks{Department of Applied Physics and Applied Mathematics, Columbia University, xl3383@columbia.edu}\qquad Robert V. Kohn \thanks{Courant Institute of Mathematical Sciences, New York University, kohn@cims.nyu.edu}}
	\date{}
	
	\maketitle
	
	\begin{abstract}
		There is a growing mechanics literature concerning the macroscopic properties of mechanism-based mechanical metamaterials. This amounts mathematically to a homogenization problem involving nonlinear elasticity. A key goal is to identify the "soft modes" of the metamaterial. We achieve this goal using methods from homogenization for some specific 2D examples -- including discrete models of the Rotating Squares metamaterial and the Kagome metamaterial -- whose soft modes are compressive conformal maps. The innovation behind this achievement is a new technique for bounding the effective energy from below, which takes advantage of the metamaterial's structure and symmetry.
	\end{abstract}
	
	\tableofcontents

	\section{Introduction}\label{sec:intro}
	
	This paper applies methods from homogenization to identify the soft modes
	and estimate the effective energy densities of certain mechanism-based
	mechanical metamaterials.
	\medskip
	
	\noindent {\sc Motivation.}
	While homogenization will be familiar to most readers, the study of
	mechanism-based mechanical metamaterials is rather new. Therefore we
	start with a brief discussion of this application area, emphasizing
	how the study of soft modes leads to the homogenization problems
	that are the focus of this paper. (For more detail on this topic,
	see Section 1.1 of \cite{li2025effective}.)
	
	We are interested in the continuum limits of certain spatially
	periodic lattices of springs, viewed as nonlinear elastic materials.
	We shall focus especially on two specific systems
	(introduced in Section \ref{subsec:rotating-squares-and-kagome}) -- the Kagome metamaterial
	and the Rotating Squares metamaterial. Both have attracted considerable
	attention in the mechanics literature as examples of ``auxetic materials.''
	In particular, they have periodic mechanisms (parametrized families of
	elastic-energy-free deformations) whose macroscopic effect is isotropic
	compression.
	
	While the macroscopic effect of a periodic mechanism is necessarily
	an affine map, these structures also have non-affine ``soft modes'' (that is,
	low-elastic-energy deformations) obtained by modulating a mechanism (see for example
	Figure \ref{fig:modulation}).
	
	It is natural to ask whether there are other soft modes, not
	obtained by modulating a mechanism. This question is especially
	relevant in systems with many mechanisms, not all of which are known
	explicitly (for example the Kagome metamaterial). To address it, one
	must begin with a rigorous definition of a ``soft mode'' -- one that
	works even for systems with many mechanisms, and which does not require
	a list or classification of the mechanisms. We proposed in
	\cite{li2025effective} that the soft modes of a mechanical metamaterial
	should be the elastic deformations whose (suitably-defined) effective energy
	density vanishes. Put differently: they are the elastic deformations whose
	spatially-averaged elastic energy vanishes in the continuum limit. This
	definition is consistent with the emerging understanding of soft modes
	in the mechanics literature, where a number of studies have considered
	specific examples \cite{czajkowski2022conformal,czajkowski2024duality,zheng2022continuum,zheng2022modeling}.
	
	For our proposal to be useful, it is crucial to learn how the
	effective energy can be estimated -- and how the soft modes can
	be identified -- in specific examples. This paper achieves that goal
	for the Kagome metamaterial, the Rotating Squares metamaterial, and
	some other examples. The innovation behind this achievement is a
	novel technique for bounding the effective energy density from below,
	which takes advantage of the metamaterial's structure and symmetry.
	We are particularly proud that our technique works for the Kagome
	example, even though it has infinitely many mechanisms (for which
	there is no known classification).
	\medskip
	
	\noindent {\sc Mathematical context.}
	Turning now to a more mathematical description of our goals: the
	structures that interest us are lattices of springs, and we are interested
	in the elastic energy associated with (potentially large) displacements
	of their nodes. Therefore the study of their macroscopic behavior
	combines the challenges associated with (i) periodic homogenization for
	nonlinear elasticity and (ii) discrete-to-continuous limits. The basics of
	homogenization for nonlinear elastic solids have been understood for about
	$40$ years \cite{braides1985homogenization,muller1987homogenization}, and discrete-to-continuous
	limits of various nonlinearly elastic spring systems have been understood
	for about $20$ years \cite{alicandro2004general}. Therefore it is
	not a surprise that our metamaterials' macroscopic behavior is hyperelastic,
	with an effective energy density that's characterized by a minimization over
	supercells. There are, however, some subtleties associated with the use of
	these tools for the study of metamaterials; in particular (a) it is important
	to include, in the elastic energy, an appropriate penalization for change of
	orientation; and (b) it is natural that the theory place no restrictions
	on the geometry of the lattice or the locations of its nodes. With such
	considerations in mind, our recent paper \cite{li2025effective} develops a new approach
	to the existence of an effective energy that models the macroscopic behavior
	of a mechanical metamaterial. To make the present paper self-contained, we shall
	summarize the framework of \cite{li2025effective} in Section \ref{sec:setup}.
	
	We are not the first to consider lower bounds for the effective energies
	of nonlinearly elastic systems. There is, in particular, a certain analogy
	between our work and the use of elastic energy minimization to study
	martensitic transformation (see for example \cite{bhattacharya2003microstructure,rindler2018calculus}).
	In that area, one considers an elastic energy
	density $W(Du)$ with finitely many ``nonlinear wells'' (one for each martensitic phase)\footnote{In more detail: one considers a frame-indifferent energy density $W \geq 0$ such that the set where $W=0$ has the form $SO(3)E_1 \cup \ldots \cup SO(3)E_k$; here $k$ is the number of martensite variants and each $E_j$ is a symmetric $3 \times 3$ matrix, which we call the ``energy-free strain'' of the $j$th variant.}. The analogue of our effective energy density is
	the quasiconvexification $QW$ of $W$, and the analogue of a soft mode
	is a deformation $u$ such that $QW(Du)=0$. (Roughly speaking, this means
	the deformation $u$ can be accommodated by an asymptotically energy-free
	mixture of martensite variants.) There are, however, some major differences
	between our problem and martensitic transformation: (a) the energy-free
	strains achieved by mechanisms form a continuum, while the energy-free
	strain of a martensite variant is a single symmetric matrix; and (b) in
	some systems (for example the Kagome metamaterial) there can be infinitely
	many mechanisms, for which we have no explicit formulas or classification; in
	martensitic transformation, by contrast, there are finitely many variants
	whose energy-free strains are known from experiments. In view of these
	differences, it is perhaps not surprising that the methods by which we
	identify the soft modes of mechanical metamaterials are quite different from
	those that have been used to understand energy-free mixtures of martensite variants.
	\medskip
	
	\noindent {\sc Main ideas and results.}
	We now offer an informal discussion of our main results (whose rigorous
	statements can be found in Section \ref{subsec:main-results}). Our arguments use the
	specific structure and symmetry of the metamaterials under
	consideration. To display the core idea as simply and accessibly as
	possible, our analysis starts in Section \ref{sec:geo-argument} by showing (for the Kagome and
	Rotating Squares examples) that the macroscopic effect of any
	periodic mechanism must be an isotropic compression (Proposition \ref{lemma-geometric}).
	While this is not a surprise for the Rotating Squares example (which is easily seen
	to have a single mechanism), it is already a striking result for
	the Kagome example since it applies to mechanisms with any periodicity,
	and makes no use of any formula for the mechanism. The proof has two main
	components:
	
	\begin{itemize}
		\item Since our structures are periodic lattices of springs, their images
		under periodic mechanisms are also periodic lattices of springs. Our argument
		identifies the fundamental parallelogram of the image, and shows that it
		has the same angles as that of the reference lattice (this is the geometric
		interpretation of Proposition \ref{prop:periodic-mechanism-macroscopic}). The isotropy of the macroscopic
		deformation follows using a linear-algebra-based argument
		(Proposition \ref{prop:our-geometric-claim}).
		
		\item Our structures include straight lines of springs. Since a mechanism is
		by definition energy-free, such a line cannot experience
		macroscopic extension. Combining this with isotropy, we deduce (in the proof of
		Proposition \ref{lemma-geometric}) that the macroscopic deformation must have the form
		$c R$ where $0 \leq c \leq 1$ is scalar and $R \in SO(2)$.
	\end{itemize}
	
	The study of soft modes is, of course, different from the
	study of mechanisms, since it requires considering low-energy
	deformations rather than zero-energy ones. Therefore our main
	result -- a lower bound on the effective energy (Theorem \ref{thm:lower-bound},
	stated at the beginning of Section \ref{subsec:main-results}) -- requires additional ideas besides
	those of  Section \ref{sec:geo-argument}. Basically, the arguments in Section \ref{sec:geo-argument} turn out
	to be rather robust: if a deformation has the form $u(x) = \lambda x + \psi(x)$
	where $\psi$ is $k$-periodic (for any $k$) and the elastic energy of $u$ is
	merely \emph{small} (rather than zero), then the arguments in Section \ref{sec:geo-argument}
	(combined with appropriate use of H\"{o}lder's inequality) show that $\lambda$ is
	\emph{close} to an isotropic compression. This argument relies heavily
	upon the variational characterization of the effective energy density, which
	permits us to consider only deformations of the form $u(x) = \lambda x + \psi(x)$
	where $\psi$ is $k$-periodic for some $k$ -- though we do not claim, and indeed
	it is not true, that all low-energy deformations have this form.
	
	While our metamaterials are lattices of springs, in defining their energies
	we include (in addition to the energies of the springs) a term that penalizes
	change of orientation. The presence of this term assures, for example,
	that a deformation which folds the lattice upon itself like an accordion does
	not have low energy. In proving our lower bound on the effective energy the presence
	of the penalization term is crucial, since it assures that a low-energy deformation
	is orientation-preserving except on a set of small measure. (We are, to be sure, not the first
	to use such a penalization term -- see for example \cite{alicandro2011integral}.)
	
	Throughout sections \ref{sec:setup}--\ref{sec:lower-bd} we discuss the Kagome metamaterial and the
	Rotating Squares metamaterial in parallel. These two examples have rather
	different microscopic structure, and yet our arguments work for both
	and they have essentially the same effective behavior. Their effective
	energy densities vanish exactly at isotropic compressions, from which it follows
	(see Section \ref{subsec:main-results}) that the soft modes are deformations $u$ such that
	$Du(x) = c(x) R(x)$ where $0 \leq c(x) \leq 1$ is scalar-valued and
	$R(x) \in SO(2)$. When $c>0$ such a map is \emph{conformal}, so it is
	natural to call these systems \emph{conformal metamaterials} \cite{czajkowski2022conformal}.
	
	To demonstrate that our methods are not limited to these examples, we apply them
	to some additional examples in Section \ref{sec:other-conformal}. However, we do not attempt
	to define a broad class of structures for which our methods work. In
	practice, all the systems discussed in this paper are conformal metamaterials,
	and it is not clear whether a similar method can be applied to systems that
	are not of this type.
	
	Concerning some practical consequences of our results, we offer the following observations: 
	
	\begin{enumerate}
		\item[(1)] Given a 2D domain $\Omega$ occupied a particular metamaterial, it is natural
		to ask which deformations of $\Omega$ are achievable by its soft modes. For the Kagome and Rotating Squares
		systems (as well as the other examples discussed in Section \ref{sec:other-conformal}), the answer is provided
		by Theorem \ref{thm:compressive-conformal}: a deformation $u : \Omega \rightarrow \mathbb{R}^2$ is a soft mode 
		if and only if it is a compressive conformal map. While this was previously known for the Rotating Squares
		system \cite{czajkowski2022conformal,zheng2022continuum}, prior to our work it was at best a conjecture for the others. 
		
		\item[(2)] It is also natural to ask what will be seen when a Dirichlet-type boundary condition is 
		imposed at $\partial \Omega$ that's not consistent with a soft mode. The deformation that minimizes our effective energy 
		$\int_{\Omega} \overline{W}^\eta (Du) \, dx $ is a natural candidate. What kind of variational problem is this? While
		we do not have a formula for $\overline{W}^\eta$, we expect it to resemble the right hand side of 
		our lower bound \eqref{eqn:eff-lower-bd}.
	\end{enumerate}

	Does knowing the soft modes of a mechanical metamaterial permit one to predict its
	response to loading? The answer may be case-dependent. In recent work on the Rotating Squares
	system \cite{czajkowski2022conformal} and on a related but more general
	family of structures \cite{zheng2022continuum,zheng2022modeling},
	the models used to match experiments involve more than knowledge of the soft
	modes -- they also take into account the leading-order elastic energy due to
	modulation of the mechanism. For the Kagome system, it seems that the energetic
	cost of modulation may be negligible \cite{tobasco2025work}. This leads us to speculate that for this system the response to loading might
	be predicted by maximizing the work done by the load within the class of soft modes.
	
	A word is in order about our use of the term ``mechanism.'' We wrote above that a mechanism is a parametrized family of 
	elastic-energy-free deformations. However, in our mathematical analysis we sometimes want to consider 
	\emph{individual} energy-free deformations. These, too, will be called mechanisms (see, for example, the 
	definition of a periodic mechanism at the beginning of Section\ref{sec:geo-argument}). We hope and expect that this 
	double usage will not lead to confusion, since the context should make our meaning clear. 
	
	\medskip
	
	\noindent {\sc Related literature.} To put our work in context, let us briefly note some
	related research threads.
	
	\begin{itemize}
		
		\item As we have already mentioned, the Rotating Squares and Kagome metamaterials have attracted
		considerable attention as examples of 2D auxetic materials  (in other words, materials which when
		compressed in a given direction, will compress rather than expand in the orthogonal direction).
		The literature also includes many other examples of auxetic materials, most of which are mechanism-based;
		see for example \cite{borcea2015geometric,borcea2018periodic,gatt2015hierarchical,grima2000auxetic,grima2005auxetic,grima2006auxetic,milton2015new}.
		A different but related challenge is the design of systems whose mechanisms achieve a specified
		family of macroscopic strains; work of this type can be found, for example, in \cite{milton2013complete,milton2013adaptable}.
		
		\item While the design of systems with interesting mechanisms has a long history, the study of such systems'
		soft modes is much more recent. In \cite{bertoldi2017flexible} one finds the viewpoint that soft modes are
		achieved by coordinated buckling of thin necks. While this is certainly true, it suggests that soft
		modes are best understood using tools from bifurcation theory. Our variational viewpoint -- that soft
		modes are deformations whose effective energy vanishes -- is quite different. While we are perhaps the
		first to formulate and study it mathematically, a similar viewpoint is implicit in \cite{czajkowski2022conformal,czajkowski2024duality,zheng2022continuum,zheng2022modeling}.
		
		\item It is natural to ask how the soft modes of mechanism-based mechanical metamaterials can be
		used for the design of interesting devices or structures. An extended answer lies far beyond the
		scope of this paper (and also beyond the competence of its authors). However we want to highlight
		the recent work of Konakovic-Lukovic et al, which designs inflatable structures using
		a variant of the Kagome metamaterial \cite{konakovic2018rapid}.
		
		\item While this paper considers only lattices of springs, periodic origami-type structures can
		also be viewed as mechanism-based mechanical metamaterials (see e.g. \cite{rocklin2024shortcuts}).
		They, too, have soft modes, whose analysis has begun to attract
		attention \cite{mcinerney2025coarse,xu2024derivation}.
		
		\item This paper studies the \emph{nonlinear elastic properties} of mechanism-based
		mechanical metamaterials. We note, however, that it is also interesting to consider other
		aspects of these systems. One topic that has received much attention is the band structure
		of a metamaterial's linear elastic waves; since the
		image of a lattice under a periodic mechanism is an entirely new lattice, mechanism-based metamaterials
		are systems with controllable spectral properties (see e.g.
		\cite{bertoldi2017flexible,fruchart2020dualities,shan2014harnessing}). In an entirely different
		direction, the Kagome system is a favorite example in the
		emerging area of topological mechanics (see e.g. \cite{mao2018maxwell,nassar2020microtwist,rocklin2017transformable}).
		
	\end{itemize}
	
	\section{Problem setup and main results}\label{sec:setup}
	We start with a brief introduction to our two key examples, the Rotating Squares metamaterial and the Kagome metamaterial (Section \ref{subsec:rotating-squares-and-kagome}). Also, to make this paper self-contained, we summarize the framework developed in \cite{li2025effective} as it applies to these examples
	(Sections \ref{subsec:lattice-setup}--\ref{subsec:eff-energy}).\footnote{For simplicity, the discussion of
		lattice metamaterials in this section is less general than the one in \cite{li2025effective}; in
		particular, it does not consider long-range springs. It is, however, general enough to treat the examples
		discussed in this paper.} Then we summarize this paper's main results in Section \ref{subsec:main-results}.
	
	\subsection{The Rotating Squares and Kagome metamaterials} \label{subsec:rotating-squares-and-kagome}
	
\begin{figure}[!htb]
		\centering
		\subfloat[]{\includegraphics[width=0.75\linewidth]{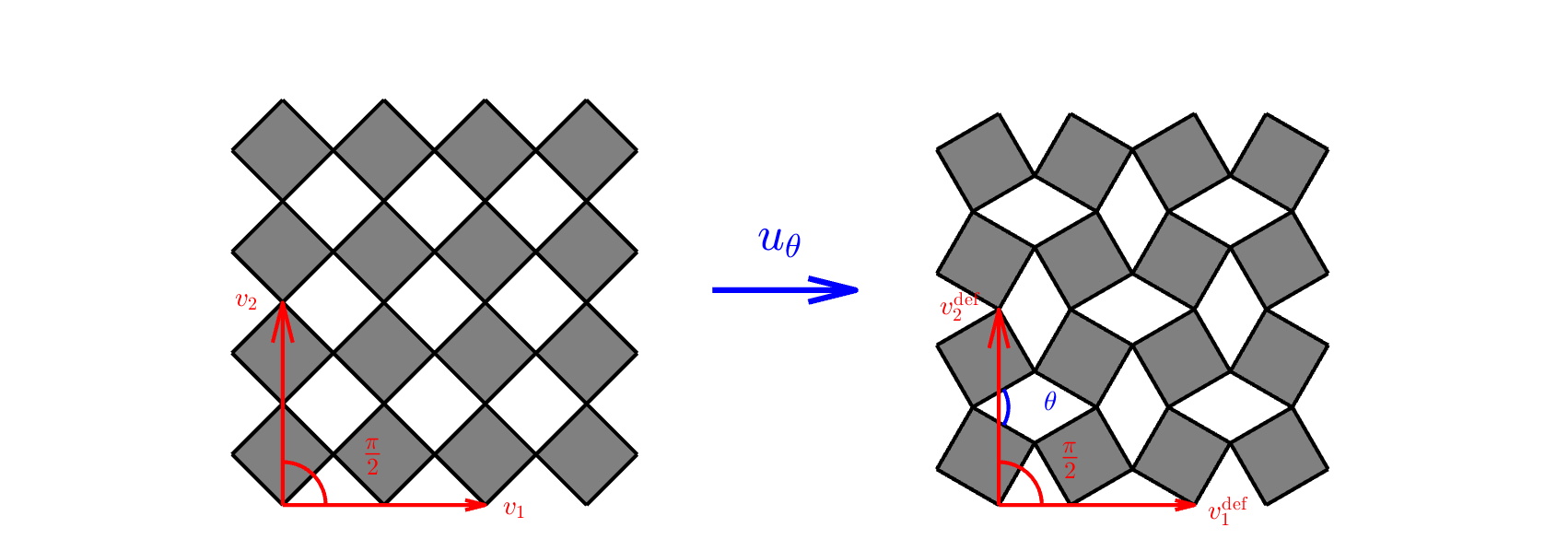}}
		\subfloat[]{\raisebox{0.25cm}{\includegraphics[width=0.2\linewidth]{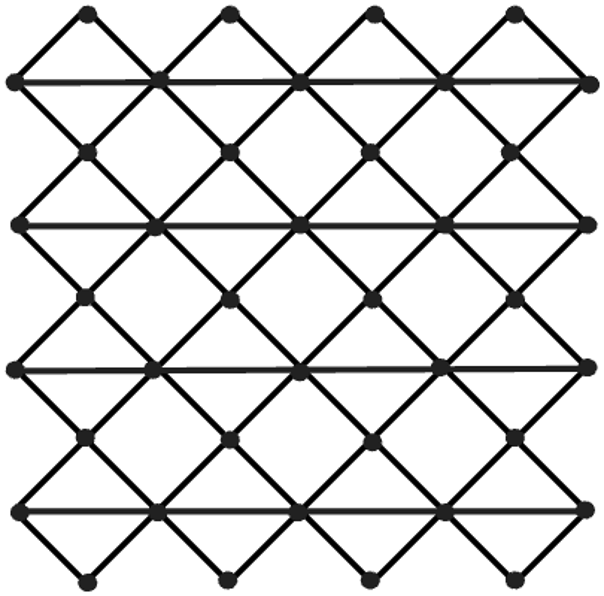}}}
		\caption{The Rotating Squares metamaterial: (a) its reference state as a cut-out and its image
			under the mechanism; (b) its reference state as a spring network, with springs connecting the diagonals of certain squares.}
		\label{fig:rs-conformal-mechanism}
	\end{figure}

	The Rotating Squares example was, it seems, first introduced by Grima and Evans \cite{grima2000auxetic}. It is among the oldest and most-studied mechanism-based mechanical metamaterials (see e.g. \cite{czajkowski2022conformal,deng2020characterization,duell2024variational,kochmann2013homogenized,zheng2022continuum,zheng2022modeling}). Viewed as a \emph{cut-out}, it is obtained by patterning a planar elastic sheet like a checkerboard and removing the white squares -- leaving the black squares, connected at their corners (which would in practice be thin necks). Viewed as a 2D elastic system, its energy-free deformations must move each square by a rigid motion; while the squares should remain connected at their corners, rotation is viewed as being elastic-energy-free (a natural idealization). This structure has a one-parameter family of mechanisms, which deform the holes from squares to parallograms; Figures \ref{fig:rs-conformal-mechanism}a and \ref{fig:rs-conformal-mechanism}b  show the reference state and its image under the mechanism, for a particular choice of the parameter. This structure also has soft modes obtained by modulating the mechanism, closely analogous to soft mode of the Kagome system shown in Figure \ref{fig:modulation}b.
	
	The analysis of this cut-out model has many subtleties, since the black squares interact via thin necks \cite{Engl2025variational}. Fortunately, there is a simpler, discrete model which captures the essential physics. This \emph{lattice-of-springs} model has nodes at the corners of the squares, with springs connecting them along the edges and diagonals of the black squares.
	Compression or extension of the springs costs elastic energy, but rotation at the nodes is free. This model's reference state -- the analogue of left side of Figure \ref{fig:rs-conformal-mechanism}a -- is shown
	in Figure \ref{fig:rs-conformal-mechanism}b.
	
	In either setting, we are interested in 2D deformations that preserve the orientation of each black square. (If change of orientation were permitted, the effective behavior would be very different -- indeed, quite degenerate -- since the structure could fold like an accordion along the lines of springs that run through the lattice.) It is therefore crucial that our variational problem include, besides the energies of the springs, a penalization for change of orientation. Note that the penalization should only be applied on the black squares (in the cut-out model) or the triangles that represent them (in the spring model), since the white squares are viewed as holes. We penalize rather than prohibit change of orientation, since the theory we are using to study the continuum limit permits penalization but not prohibition.
	
	From the cut-out perspective, it would be reasonable to impose the condition that the images of the black squares 
	not overlap. We alert the reader that our lattice-of-springs-based model does not include such a condition. (The
	difficulty of including non-interpenetration conditions in variational problems from elasticity is well-known. We
	do not know the existence of an effective energy in such a setting, though we wonder whether the method of 
	\cite{ciarlet1987injectivity} could be helpful.) 
	
	We note that while the spring model is simpler than the cut-out model, the two have the same mechanisms. Indeed, in the spring model an elastic deformation is determined by its values at
	the nodes. We take the view that the deformation of a triangle is the affine interpolant of
	the deformations of its corners. Now, in the spring model a deformation is energy-free when (i) its spring energies vanish -- i.e. it leaves the lengths of the triangles' sides unchanged, and (ii) it preserves the orientation of each triangle. It follows by elementary geometry that each triangle moves by a rigid motion.
	
\begin{figure}[!htb]
		\centering
		\subfloat[]{\includegraphics[width=0.3\textwidth]{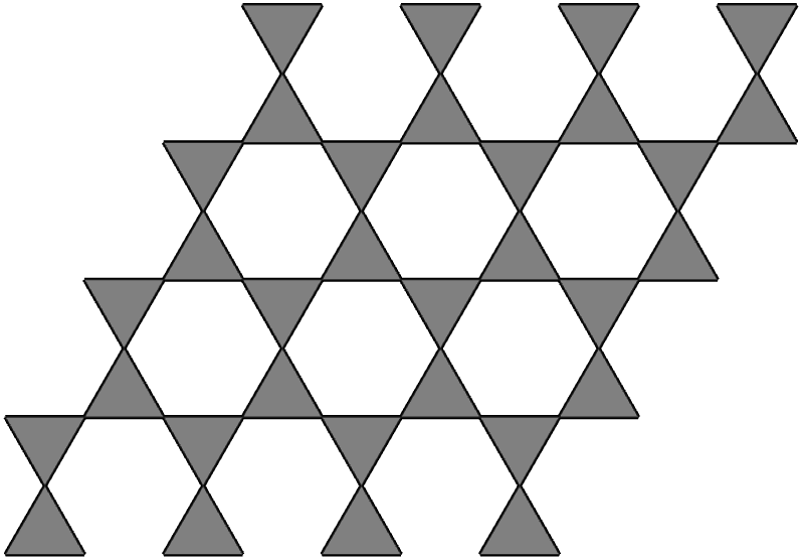}}\hfil
		\subfloat[]{\includegraphics[width=0.3\textwidth]{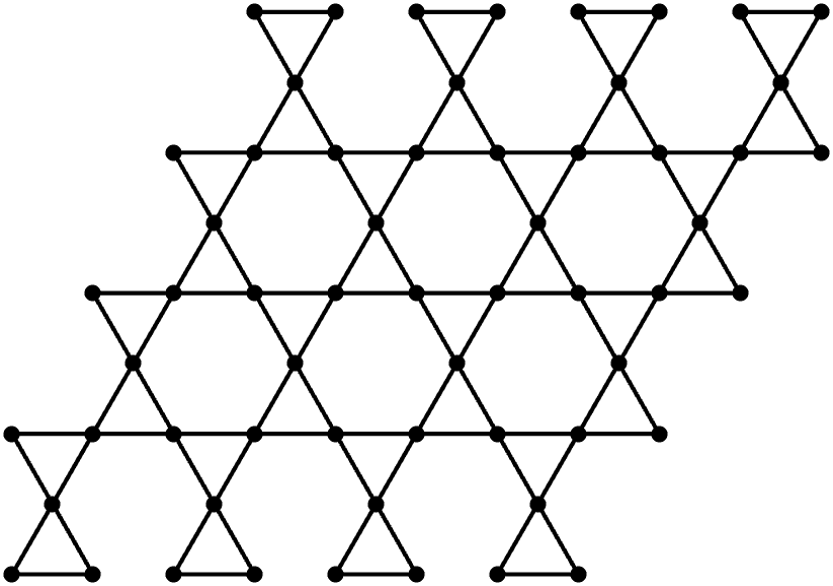}}
		\caption{The reference state of the Kagome metamaterial, viewed (a) as a cut-out, and (b) as a lattice of springs in which each node is connected to its nearest neighbors by springs.}\label{fig:kagome-cutout-and-spring-models}
	\end{figure}

	We turn now to the Kagome example, which has also been studied extensively (see
	e.g. \cite{kapko2009collapse,li2023some,shan2014harnessing,sun2012surface}). It, too has a cut-out model and a spring model. The former is obtained by patterning 2D elastic sheet as shown in Figure \ref{fig:kagome-cutout-and-spring-models}a
	and removing the hexagons -- leaving the triangles, connected at their corners where rotation is free. The spring model has nodes at the triangles' corners and springs at their sides, as shown in Figure \ref{fig:kagome-cutout-and-spring-models}b. As in the Rotating Squares setting, our
	variational problem will include, besides the spring energies, a term penalizing
	change of orientation in the triangles. The cut-out and spring models have the same mechanisms, by the same argument given above for the Rotating Squares example.
	
\begin{figure}[!htb]
		\centering
		\subfloat[]{\includegraphics[width=0.3\textwidth]{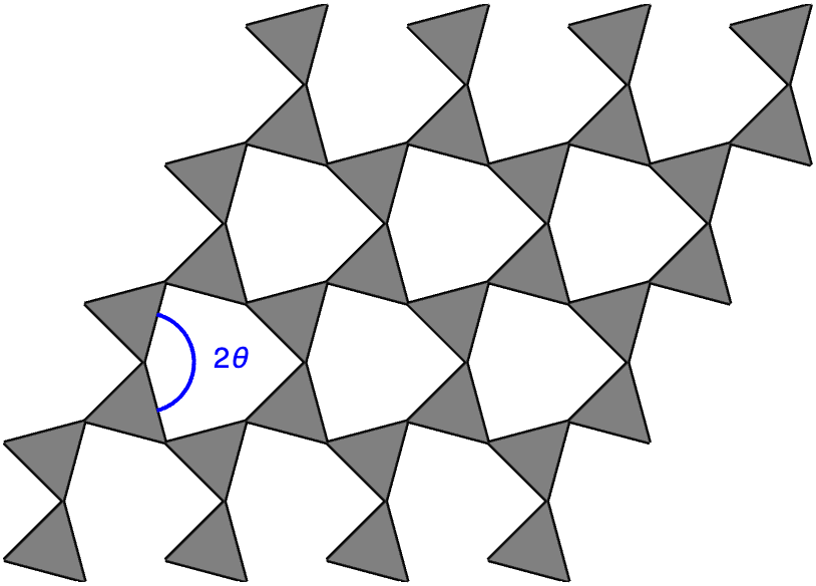}}\hfil
		\subfloat[]{\includegraphics[width=0.33\textwidth]{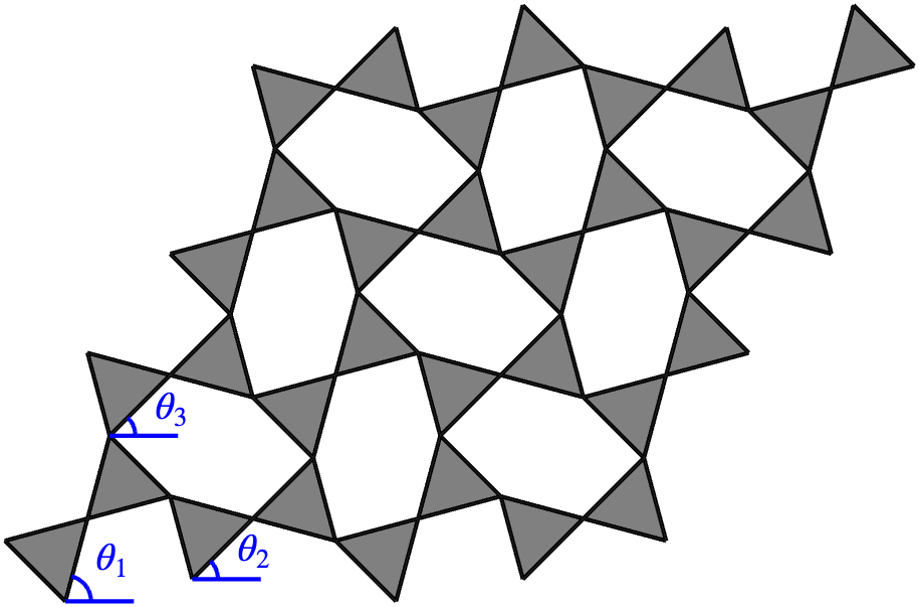}}\\
		\subfloat[]{\includegraphics[width=0.5\textwidth]{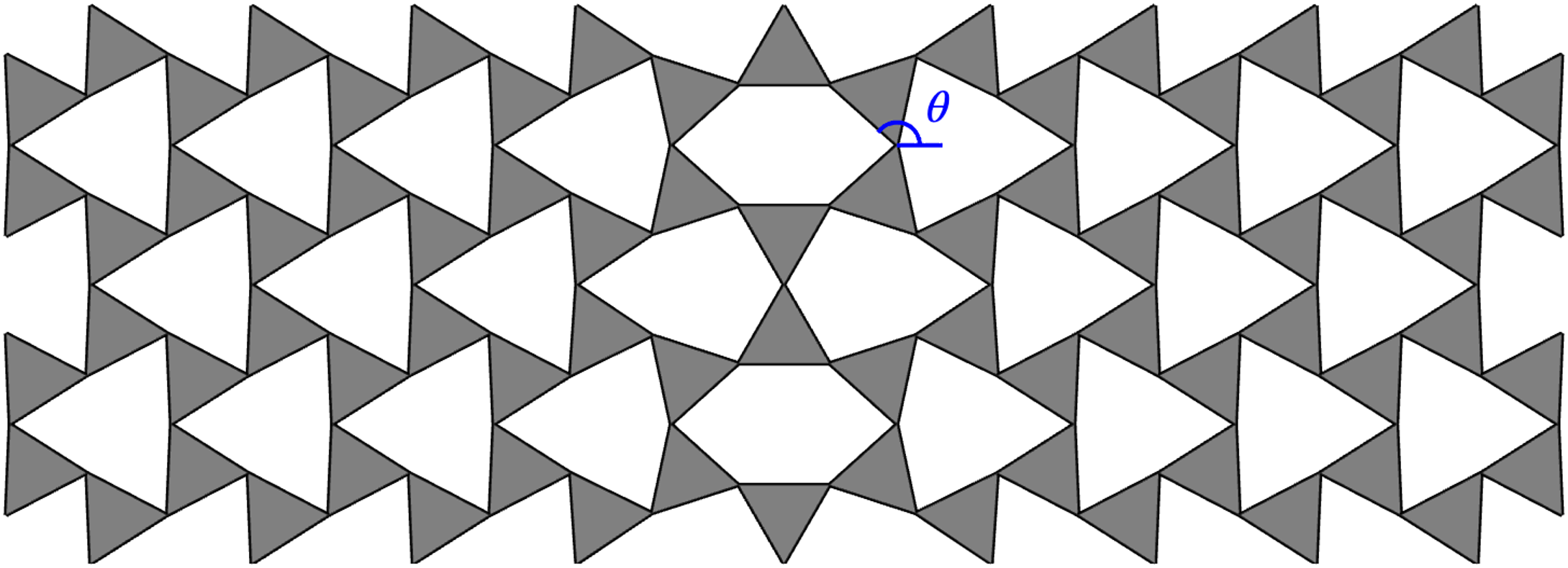}}
		\caption{Some mechanisms of the Kagome metamaterial: (a) the deformed state of the one-periodic mechanism $u_\theta$, controlled by a single angle $\theta$; (b) the deformed state of a two-periodic
			mechanism $u_{\theta_1, \theta_2, \theta_3}$, controlled by three angles $\theta_1, \theta_2, \theta_3$ (for more detail see \cite{li2023some}); (c) the deformed state of a non-periodic mechanism $u_\theta$, controlled by a single angle $\theta$ (for more detail see Section \ref{subsec:non-periodic}).} \label{fig:kagome-mechanisms}
	\end{figure}

\begin{figure}[!htb]
		\centering
		\subfloat[]{\includegraphics[width=0.48\textwidth]{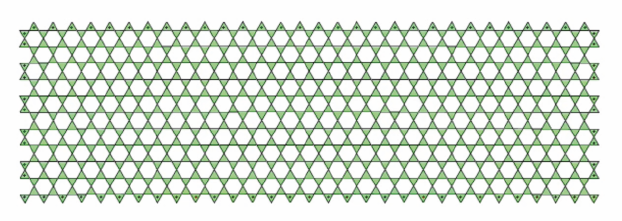}}\\
		\subfloat[]{\includegraphics[width=0.42\textwidth]{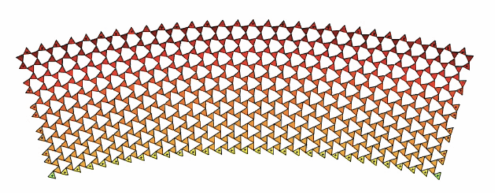}}
		\subfloat[]{\includegraphics[width=0.42\textwidth]{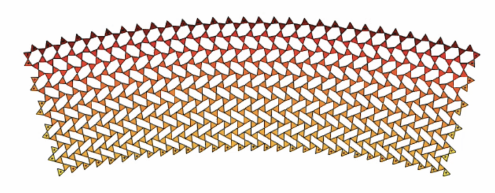}}
		\caption{Soft modes in the Kagome lattice: (a) a rectangle filled with the Kagome metamaterial
			in its reference state; (b) the image of a soft mode, achieved microscopically by modulating
			the mechanism shown in Figure \ref{fig:kagome-mechanisms}a; (c) the image of the same soft mode, achieved
			microscopically by modulating a different mechanism. The colors in (b) and (c) indicate
			the rotation angles of the deformed structure's triangles.}
		\label{fig:modulation}
	\end{figure}

	Unlike the Rotating Squares system, the Kagome metamaterial has infinitely many mechanisms (see e.g. \cite{kapko2009collapse,li2023some}). There is a one-parameter family of mechanisms with the same periodicity as the reference lattice; its image is shown in Figure \ref{fig:kagome-mechanisms}a for a particular value of the parameter. But there is also a 3-parameter family of mechanisms whose periodicity is twice that of the reference lattice (Figure \ref{fig:kagome-mechanisms}b), and there are mechanisms with larger periodicity as well. While the macroscopic behavior of each explicitly-known periodic mechanisms is isotropic compression, prior to the present work it was an open question whether \emph{every} periodic mechanism had this property. Our Proposition \ref{lemma-geometric} provides an affirmative answer. We note in passing that the Kagome system also has mechanisms that are not periodic -- including one in which two symmetry-related versions of the one-periodic mechanism meet at an energy-free wall (Figure \ref{fig:kagome-mechanisms}c).
	
	As noted earlier, mechanism-based mechanical metamaterials have soft modes obtained by modulating the mechanism. The associated macroscopic deformation must (as we explain in Section \ref{subsec:main-results}) be a compressive conformal map. Figure \ref{fig:modulation} shows soft modes of the Kagome lattice obtained by modulating two different mechanisms. We note, however, that since the Kagome system has many mechanisms, the microscopic character of a soft mode need not be so simple -- for example, it can have domains that use
	different mechanisms.
	
	\subsection{Introduction to  some lattice systems of springs}\label{subsec:lattice-setup}
	We begin by introducing some notation for a class of 2D lattice systems of springs that includes the Kagome and the Rotating Squares metamaterial. We start with a \textit{unit cell} $U$ (a parallelogram containing the origin) and two vectors
	$v_1, v_2 \in \mathbb{R}^2$ such that the translated cells $U + \alpha$ tile\footnote{The translated copies of the unit cell may have overlapping boundaries, but their interiors remain distinct and non-intersecting.}
	the entire plane when $\alpha = \alpha_1 v_1 + \alpha_2 v_2$ with $\alpha_1,\alpha_2 \in \mathbb{Z}$. We call $v_1, v_2$ the \textit{lattice vectors}. To identify the nodes of the lattice, we fix a basic set of nodes in the unit cell, $V = \{p_1,\dots, p_{|V|}\} \subset U$;
	the full set $\mathcal{V}$ of nodes consists of all translates of elements of $V$:
	\begin{equation} \label{eqn:nodes-of-lattice-unscaled}
		\mathcal{V} = \bigcup_{\alpha_1,\alpha_2 \in \mathbb{Z}} (V + \alpha_1 v_1 + \alpha_2 v_2).
	\end{equation}
	We assume that no two elements of $V$ are lattice translates of one another, so each node of the lattice is \emph{uniquely} expressible as $p + \alpha$ for some $p \in V$ and
	$\alpha  = \alpha_1 v_1 + \alpha_2 v_2$. 
	
\begin{figure}[!htb]
		\begin{minipage}{.45\linewidth}
			\centering
			\subfloat[]{\label{fig:kagome-unit}\includegraphics[width=0.85\textwidth]{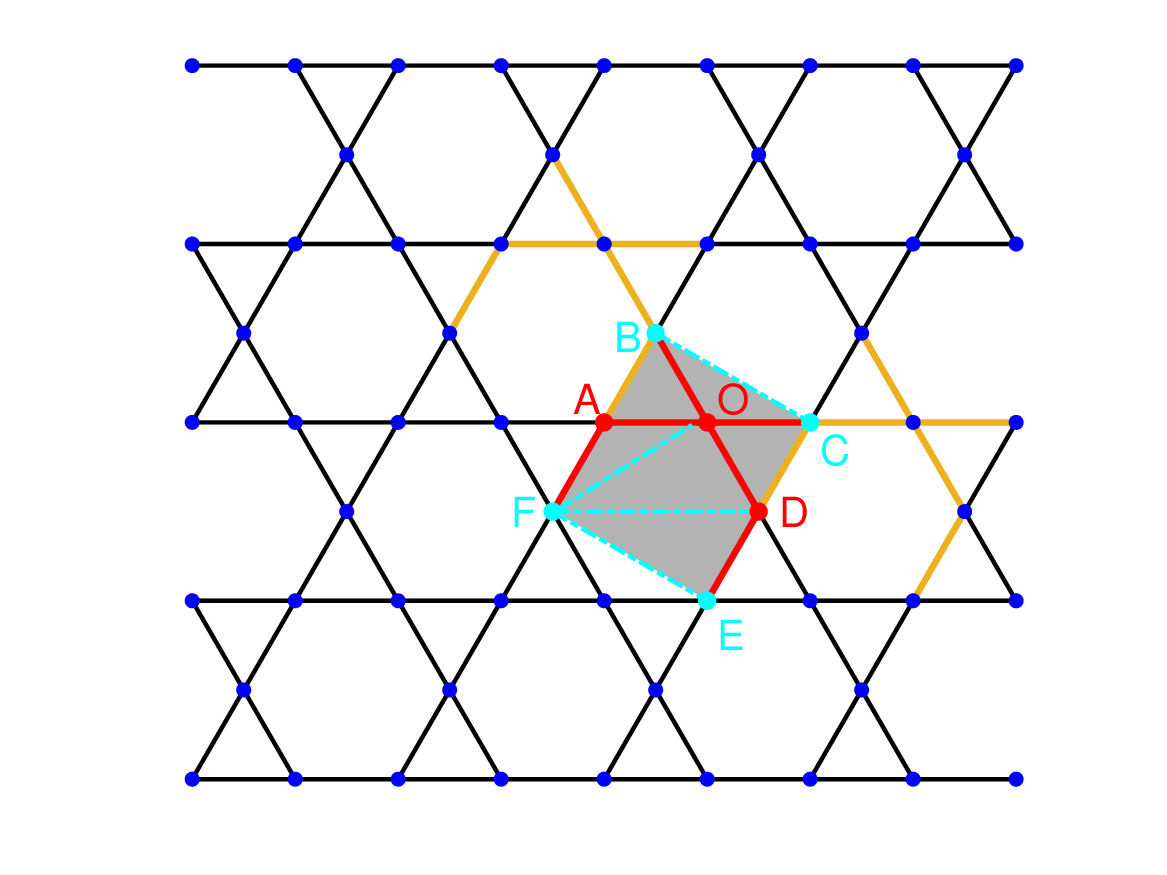}}
		\end{minipage}
		\begin{minipage}{.45\linewidth}
			\centering
			\subfloat[]{\label{fig:rs-unit}\includegraphics[width=0.85\textwidth]{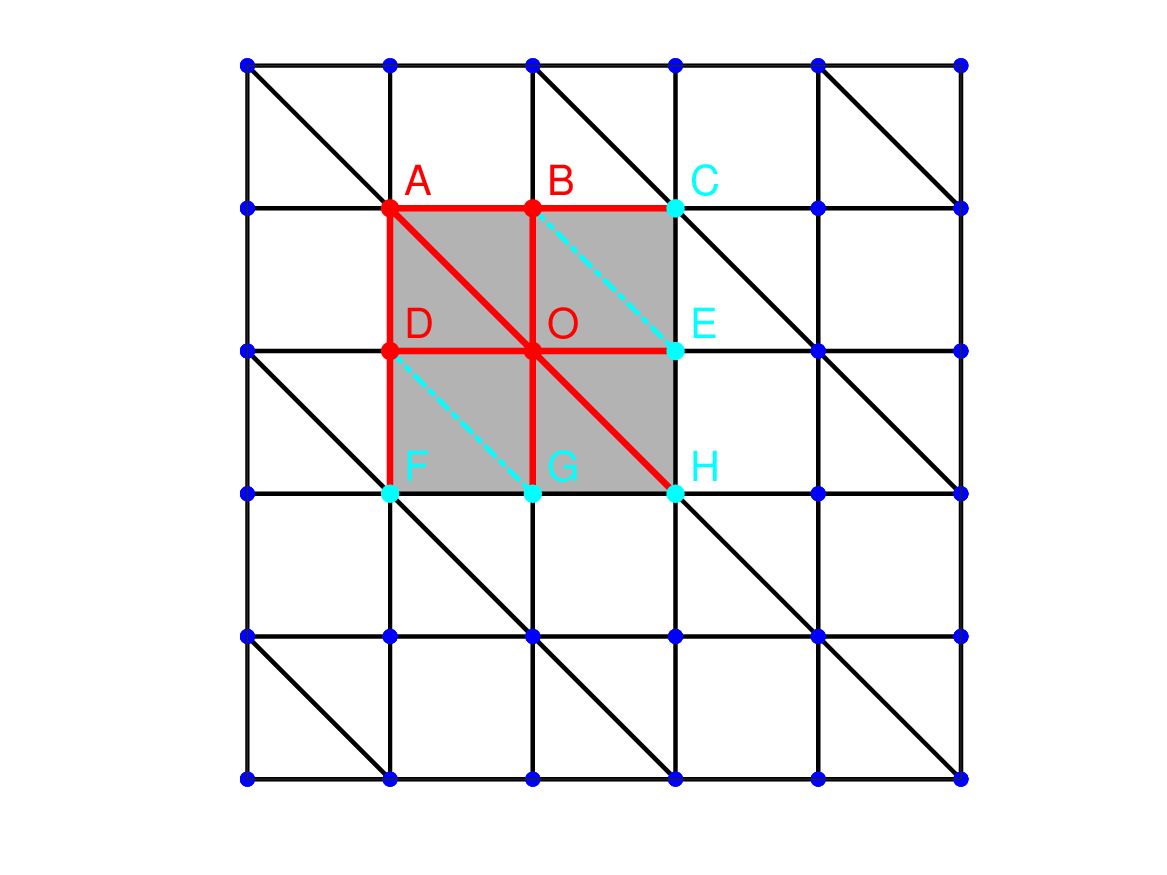}}
		\end{minipage}
		\caption{(a) The Kagome lattice: The shaded rectangle represents the unit cell $U$ for the Kagome lattice, which contains three vertices $A,O,D$ marked in red. These vertices can be translated to obtain the entire lattice. The solid red edges are those included in the energy $E_\text{spr}(u,U)$ in equation \eqref{eqn:kagome-spr-energy}. A translated copy of these edges is marked in yellow to illustrate that all edges in the Kagome lattice can be viewed as translated copies of the red solid edges. The dotted lines indicate the triangular mesh used to interpolate the admissible deformations. (b) The Rotating Squares lattice. The nodes associated with the unit cell (our set $V$) are marked in red; nodes not in $V$ but used in the energy $E(u,U)$ are marked in cyan; springs counted in $E_\text{spr}(u,U)$ are marked by red solid lines; artificial edges used only for the triangularization of $U$ are marked by cyan dotted lines; the shaded area is $U$.}
	\end{figure}

	As an example, consider the 2D Kagome lattice shown in Figure \ref{fig:kagome-unit}. A convenient choice of its
	unit cell $U$ is the rectangle with vertices $B,C,E,F$, and a convenient choice of the basic set of nodes is $V = \{A,O,D\}$.
	If we choose the distance between two nearest nodes to be $1$, then the lattice vectors are $v_1 = (2,0)^T$ and $v_2 = (1,\sqrt{3})^T$.
	
	We want to endow such a lattice with an elastic energy. To do so, it is important to be clear about what we mean by an elastic
	deformation. We take the view that a deformation is an $\mathbb{R}^2$-valued function defined \emph{only at nodes}.
	(Our situation is thus different from the theory of ``reticulated structures,'' discussed e.g.
	in \cite{cioranescu2012homogenization}, where the deformations are defined on sets with nonzero volume.)
	
	To formulate the elastic energy in terms of deformed positions of the nodes, we model the lattice system as a two-dimensional elastic material by viewing the connections between nodes as Hookean springs. It is convenient to define the energy of a unit cell, $E(u,U)$, as the sum of the spring energies associated with the springs contained within the unit cell $U$\footnote{The lattice materials considered in this paper are sufficiently simple that no springs within a unit cell cross its boundary. For a more general treatment of cases with boundary-crossing springs, we refer the reader to our paper \cite{li2025effective}.}. Since we consider spring energies, it is not hard to observe that $E(u,U)$ is nonnegative 
	\begin{equation*}
		E(u,U) \geq 0
	\end{equation*}
	and translation-invariant 
	\begin{equation*}
		E(u,U) = E(u+c,U)
	\end{equation*}
	when $c\in \mathbb{R}^2$ is a translation vector, i.e. it takes the same value at every node. 
	
	As an example, consider our spring model of the Kagome metamaterial, with Hookean springs connecting each pair of nearest-neighbor nodes. If the unit cell is chosen as shown in Figure \ref{fig:kagome-unit}, then it is convenient to let $E(u,U)$ be the energy of the six springs $AO,BO,CO,DO,AF,DE$, since each spring in the lattice is (uniquely) a translate of one of these. With this choice (and taking all the springs to be the same, and making a choice of the spring constant) the spring energy of a translated unit cell $U + \alpha$ is
	\begin{equation}\label{eqn:kagome-spr-energy}
		\begin{aligned}
			E_\text{spr}(u,U+\alpha) &= \Bigg(\Big|u(A+\alpha)-u(O+\alpha)\Big|-|A-O|\Bigg)^2 + \Bigg(\Big|u(B+\alpha)-u(O+\alpha)\Big|-|B-O| \Bigg)^2\\
			& + \Bigg(\Big|u(C+\alpha)-u(O+\alpha)\Big|-|C-O|\Bigg)^2 + \Bigg(\Big|u(D+\alpha)-u(O+\alpha)\Big|-|D-O|\Bigg)^2\\
			& + \Bigg(\Big|u(A+\alpha)-u(F+\alpha)\Big|-|A-F|\Bigg)^2 + \Bigg(\Big|u(D+\alpha)-u(E+\alpha)\Big|-|D-E|\Bigg)^2 .
		\end{aligned}
	\end{equation}
	However, as discussed in Section \ref{sec:intro}, the elastic energy with only spring energies permits change of local orientations. In the Kagome example, a non-physical folding deformation that flips the triangle $\Delta AOB$ towards $\Delta COD$ has zero spring energy. To avoid considering such folding deformations, we add a penalty energy with penalty constant $\eta > 0$, defining the new unit cell energy as $E^\eta(u,U+\alpha)$
	\begin{equation}\label{eqn:kagome-intro-energy}
		E^\eta(u,U+\alpha) := E_\text{spr}(u,U+\alpha) + E_\text{pen}^\eta(u,U+\alpha)
	\end{equation}
	where the penalty energy is defined as
	\begin{equation}\label{eqn:kagome-pen-energy}
		E_\text{pen}^\eta(u,U+\alpha) = |\Delta AOB| \; f^\eta(\det Du(x+\alpha)|_{\Delta AOB}) + |\Delta COD| \; f^\eta(\det Du(x+\alpha)|_{\Delta COD})
	\end{equation}
	and the function $f^\eta(t)$ is chosen as
	\begin{equation}\label{eqn:pen-f}
		f^\eta(t) = \begin{cases}
			1/\eta, & t \leq 0,\\
			0, &t> 0.
		\end{cases}
	\end{equation}
	By adding the penalty energy $E_\text{pen}^\eta(u,U+\alpha)$ and summing over all $\alpha$, we penalize change of orientation for all the equilateral triangles in the Kagome metamaterial. Here we have chosen $f^\eta$ to be a step function -- a simple choice that satisfies two key criteria: (1) there is no penalty energy when the orientation is preserved, i.e. $f^\eta(\det D u) = 0$ when $\det D u > 0$; and (2) when the penalty energy is active on a triangle, the penalty energy can be large by choosing $\eta$ small. In general, one can choose any function $f^\eta(t)$ that satisfies these two criteria.
	
	For the Rotating Squares metamaterial, we shall work with the reference lattice shown in Figure \ref{fig:rs-unit} (which is rotated by $\pi/4$ relative to the one in Figure \ref{fig:rs-conformal-mechanism}). It amounts to a square lattice augmented by diagonal edges on alternating squares (a diagonal edge is added to make the square rigid). A convenient choice of the unit cell $U$ is a big square containing four small squares (the region $ACHF$ in Figure \ref{fig:rs-unit}) and the basic set of nodes is $V = \{O,A,B,D\}$. If we choose the distance between two nearest nodes to be 1, then the lattice vectors are $v_1 = (2,0)^T$ and $v_2 = (0,2)^T$. The spring energy in the Rotating Squares metamaterial is the aggregate energy of the $10$ springs $AB$, $AO$, $AD$, $BC$, $BO$, $DO$, $EO$, $DF$, $OG$, and $OH$. It will be convenient in Section \ref{subsec:averaged-energy-kU} for the diagonal springs AO and OH to have spring constants that are twice those of the other springs, so for Rotating Squares metamaterial we take the spring energy of the unit cell to be:
	\begin{equation}\label{eqn:rs-spring-energy}
		\begin{aligned}
			&E_\text{spr}(u,U) = \Bigg(\Big|u(A)-u(B)\Big|-|A-B|\Bigg)^2 + 2\Bigg(\Big|u(A)-u(O)\Big|-|A-O| \Bigg)^2 \\
			+ \, &\Bigg(\Big|u(B)-u(O)\Big|-|B-O| \Bigg)^2 + \Bigg(\Big|u(D)-u(O)\Big|-|D-O| \Bigg)^2 + \Bigg(\Big|u(A)-u(D)\Big|-|A-D| \Bigg)^2\\
			+ \, &\Bigg(\Big|u(B)-u(C)\Big|-|B-C| \Bigg)^2 + \Bigg(\Big|u(D)-u(F)\Big|-|D-F| \Bigg)^2 + \Bigg(\Big|u(O)-u(G)\Big|-|O-G| \Bigg)^2\\
			+ \,  & 2\Bigg(\Big|u(O)-u(H)\Big|-|O-H| \Bigg)^2 + \Bigg(\Big|u(O)-u(E)\Big|-|O-E| \Bigg)^2 \; .
		\end{aligned}
	\end{equation}
	This determines the spring energy on each translated cell $U+\alpha$ by periodicity. To avoid change of orientation in the rigid squares (the squares with diagonal springs), we add penalty terms on four triangles $\Delta AOB$, $\Delta AOD$, $\Delta EOH$ and $\Delta GOH$ with penalty energy 
	\begin{equation}\label{eqn:rs-penalty-energy}
		\begin{aligned}
			E_\text{pen}^\eta(u,U)&= |\Delta AOB| f^\eta(\det D u(x)|_{\Delta AOB}) + |\Delta AOD| f^\eta(\det D u(x)|_{\Delta AOD}) \\
			&+ |\Delta EOH| f^\eta(\det D u(x)|_{\Delta EOH}) + |\Delta GOH| f^\eta(\det D u(x)|_{\Delta GOH}) \;,
		\end{aligned}
	\end{equation}
	and we use the same $f^\eta$ as in \eqref{eqn:pen-f}. The unit cell energy $E^\eta(u,U)$ is defined as the sum of the spring energy \eqref{eqn:rs-spring-energy} and the penalty energy \eqref{eqn:rs-penalty-energy}, as in \eqref{eqn:kagome-intro-energy}.
	
	We alert the reader that while our energies penalize change of orientation, they do not penalize or prohibit interpenetration. Thus, for example, in the mechanism of the Rotating Squares metamaterial shown in Figure \ref{fig:rs-conformal-mechanism}, the parameter $\theta$ need not be positive. In the cut-out version of this structure, the holes are reduced to slits when $\theta = 0$, and the macroscopic deformation is isotropic compression by a factor of $1/\sqrt{2}$ (consistent with reduction of area by a factor of $1/2$). As $\theta$ becomes negative the elastic squares begin to interpenetrate and the macroscopic deformation is isotropic compression by a factor $c < 1/\sqrt{2}$. When $\theta = -\pi/2$ the elastic squares all lie on top of one another (i.e. they are mapped to a single elastic square) and the macroscopic deformation is $0$. While it would be natural in the cut-out setting to avoid interpenetration, it seems more natural in our lattice-of-springs model setting to permit it. Moreover, from a strictly mathematical viewpoint, a theorem should be stated with minimal hypotheses; therefore it is favorable that our main result (a lower bound on the effective energy in Theorem \ref{thm:lower-bound}) holds even in a model that permits interpenetration.
	
	While an elastic deformation is characterized by its values at the nodes of the lattice, we want to also view it as a piecewise linear function defined on a suitable mesh. This is useful (a) because deformations of the scaled lattice can then be viewed as functions in a finite-dimensional subspace of $H^1$, and (b) because our theorem on the existence of an effective energy requires that $E^\eta(u,U)$ satisfy certain conditions, whose statement involves the piecewise linear version of $u$. Therefore for a given 2D lattice metamaterial, we fix -- in addition to the structure introduced so far -- a triangulation of the unit cell $U$. The theory in \cite{li2025effective} requires that 
	\begin{equation} \label{eqn:unscaled-upper-bound}
		E^\eta(u,U) \leq C_1 \Big( |D u|^2_{L^2(U)} + |U| \Big )
	\end{equation}
	and
	\begin{equation} \label{eqn:unscaled-lower-bound}
		E^\eta(u,U) \geq \max \Big\{ C_2 \Big( |D u|^2_{L^2(U))} - D_2 |U| \Big) , 0 \Big\}
	\end{equation}
	for some positive constants $C_1$, $C_2$, and $D_2$ that only depend on $\eta$ and the geometry of the lattice; here the $L^2$ norms on the right hand side refers to the piecewise linear version of $u$. The lattice systems considered in this paper are simple enough that (i) we can always use a triangulation of $U$ whose vertices are the lattice nodes in $\overline{U}$, and (ii) the triangulation can be chosen so that \eqref{eqn:unscaled-upper-bound} and \eqref{eqn:unscaled-lower-bound} hold. (These bounds are proved in \cite{li2025effective} for the Kagome example using the unit cell shown in Figure \ref{fig:kagome-unit} and the triangular mesh consisting of $\Delta AOB, \Delta BOC$, $\Delta COD$, $\Delta AOF, \Delta DOF, \Delta DEF$. They are also proved there for the Rotating Squares example with the triangularization shown in Figure \ref{fig:rs-unit}. The method used to prove these bounds is quite flexible, extending easily to the examples discussed in Section \ref{sec:other-conformal}.) 
	
	We offer two comments on the upper and lower bounds \eqref{eqn:unscaled-upper-bound} -- \eqref{eqn:unscaled-lower-bound}.  First: since we are interested in structures with mechanisms, the presence of a negative term $D_2 |U|$ in the lower bound is crucial. Indeed, for the lower bound to hold, $D_2 |U|$ must clearly be larger than the maximum of $|D u|^2_{L^2(U)}$ as $u$ ranges over mechanisms (that is, over deformations such that $E^\eta(u,U) = 0$). Our second comment is that it is natural to use the $L^2$ norm of $D u$ (rather than some other $L^p$ norm) since the elastic energy comes from Hookean springs (whose energy grows quadratically when $D u$ gets large). However, our theory permits the use of springs whose response is non-Hookean at small strains; indeed, the conditions \eqref{eqn:unscaled-upper-bound} and \eqref{eqn:unscaled-lower-bound} place no constraint on the springs' character at small or even moderate strains, due to the terms involving $|U|$ on the right hand side.

	\begin{remark} \label{rmk:other-spring-constants}
		We have made specific choices of the spring constants in \eqref{eqn:kagome-spr-energy} and \eqref{eqn:rs-spring-energy} only for simplicity. Our results would remain
		valid with any (strictly positive) choice of the $6$ spring constants used for the Kagome unit cell, or the $10$ spring constants used
		for the Rotating Squares unit cell. This is elementary, since the unit cell energies obtained this way would be bounded below and 
		above by constant multiples of the ones we analyze.
	\end{remark}
	
	\subsection{The scaled energy and our admissible deformations} \label{subsec:scaled-energy}
	
	Now we introduce the scaled lattice and its scaled energy. The nodes of the scaled lattice are $\mathcal{V}^\epsilon := \epsilon \mathcal{V}$ and the translated unit cells of this lattice are $\epsilon U + \alpha$ with $\alpha = \alpha_1 v_1 + \alpha_2 v_2$ for $\alpha_1, \alpha_2 \in \epsilon \mathbb{Z}$. The scaled unit cell energy $E^{\epsilon,\eta}(u^\epsilon, \epsilon U + \alpha)$ is defined via elasticity
	scaling
	\begin{equation}\label{eqn:elasticity-scaling}
		E^{\epsilon,\eta}(u^\epsilon, \epsilon U + \alpha) = \epsilon^2 E^\eta(u,U), \qquad u^\epsilon(x) = \epsilon u\left(\frac{x-\alpha}{\epsilon}\right),
	\end{equation}
	which matches the behavior of Hookean springs. As an example: for the Kagome metamaterial with $U$ and $E^\eta(u,U)$ given by Figure \ref{fig:kagome-unit} and Figure \ref{eqn:kagome-intro-energy}, if all the springs have length $1$ in the unscaled
	setting then
	\begin{equation*}
		\begin{aligned}
			E^{\epsilon,\eta}(u^\epsilon,\epsilon U+\alpha) &= \Bigg(\Big|u^\epsilon(\epsilon A+\alpha)-u^\epsilon(\epsilon O+\alpha)\Big|-\epsilon\Bigg)^2 + \Bigg(\Big|u^\epsilon(\epsilon B+\alpha)-u^\epsilon(\epsilon O+\alpha)\Big|-\epsilon\Bigg)^2\\
			&+ \Bigg(\Big|u^\epsilon(\epsilon C+\alpha)-u^\epsilon(\epsilon O+\alpha)\Big|-\epsilon\Bigg)^2 + \Bigg(\Big|u^\epsilon(\epsilon D+v_\alpha)-u^\epsilon(\epsilon O+\alpha)\Big|-\epsilon\Bigg)^2\\
			&+ \Bigg(\Big|u^\epsilon(\epsilon A+\alpha)-u^\epsilon(\epsilon F+\alpha)\Big|-\epsilon\Bigg)^2 + \Bigg(\Big|u^\epsilon(\epsilon D+\alpha)-u^\epsilon(\epsilon E+\alpha)\Big|-\epsilon\Bigg)^2\\
			&+ \epsilon^2 |\Delta AOB| \Bigg(f^\eta(D u^\epsilon(x+\alpha)|_{\epsilon \Delta AOB}) + f^\eta(D u^\epsilon(x+\alpha)|_{\epsilon \Delta COD})\Bigg),
		\end{aligned}
	\end{equation*}
	where $\epsilon \Delta AOB$ and $\epsilon \Delta COD$ are $\epsilon$-scale triangles. Our unscaled upper and lower bounds \eqref{eqn:unscaled-upper-bound} and \eqref{eqn:unscaled-lower-bound} have scaled versions. Their right hand sides involve the piecewise linearization of $u^\epsilon$ (determined by our unscaled piecewise linearization scheme and elasticity scaling).
	
	Since we will primarily work with the scaled formulation, we state here the full list of conditions required when the theory of \cite{li2025effective} is specialized to systems (like those considered here) where the springs associated with the basic cell $U$ connect nodes in $\overline{U}$ and the triangularization of $U$ uses only nodes in $\overline{U}$. They are that 
	
	\begin{enumerate}[(1)]
		\item the energy on the $\epsilon$-scale unit cell is periodic, i.e. we have
		\begin{equation}\label{eqn:unit-cell-eps-periodicity}
			E^{\epsilon,\eta} (u^\epsilon(x+\alpha),\epsilon U+\alpha) = E^{\epsilon,\eta} (u^\epsilon,\epsilon U)
		\end{equation}
		for any $\alpha = \alpha_1 v_1 + \alpha_2 v_2$ with $\alpha_1,\alpha_2 \in \epsilon \mathbb{Z}$ ;
		
		\item the energy on the $\epsilon$-scale unit cell is translation-invariant, in the sense that for any vector
		$c \in \mathbb{R}^N$, we have
		\begin{equation}\label{eqn:unit-cell-translation-invariant}
			E^{\epsilon,\eta} (u^\epsilon,\epsilon U+\alpha) = E^{\epsilon,\eta} (u^\epsilon+c,\epsilon U+\alpha) \, ;
		\end{equation}
		
		\item an upper bound: there exists $C_1 > 0$ (independent of $\alpha$ and $\epsilon$) such that
		\begin{equation}\label{eqn:unit-cell-upper}
			E^{\epsilon,\eta} (u^\epsilon,\epsilon U+\alpha) \leq C_1 \Big(|D u^\epsilon|^2_{L^2(\epsilon U+\alpha)} +
			|\epsilon U+\alpha|\Big);
		\end{equation}
		
		\item a lower bound: there exist $C_2 > 0$ and $D_2 \geq 0$ (independent of $\alpha$ and $\epsilon$) such that
		\begin{equation}\label{eqn:unit-cell-lower}
			E^{\epsilon,\eta} (u^\epsilon,\epsilon U + \alpha) \geq \max \Big\{C_2 \Big(|D u^\epsilon|^2_{L^2(\epsilon U+\alpha)} -
			D_2|\epsilon U + \alpha|\Big), 0 \Big\} \, .
		\end{equation}
	\end{enumerate}
	We emphasize that these conditions are valid for the scaled versions of $E^\eta(u,U)$ that we introduced above for the Kagome and Rotating Squares metamaterials. 
	
	The energy of a domain $\Omega$ filled by the scaled structure is, roughly speaking, the sum of the scaled energies of all translates of $\epsilon U$ that lie inside $\Omega$. Since the lattice systems we consider have the property that the springs associated with the unit cell are inside the unit cell, when calculating the energy of $\Omega$ we don't need to be concerned about springs that cross $\partial \Omega$. 
	
	We now make this discussion more precise. When considering the limiting energy of a fixed domain $\Omega$, it is natural to focus on deformations that are defined at lattice nodes in $\Omega$, in other
	words $u^\epsilon$ in
	\begin{equation}\label{eqn:admissible-eps}
		\mathcal{A}_\epsilon(\Omega) = \{u^\epsilon(x)\; |\; u^\epsilon(x) \text{ has values on } \mathcal{V}^\epsilon \cap \Omega\} \, .
	\end{equation}
	The energy of $\Omega$ for a deformation $u^\epsilon \in \mathcal{A}_\epsilon(\Omega)$ takes the form
	\begin{equation}\label{eqn:intro-energy-of-Omega}
		E^{\epsilon,\eta}(u^\epsilon, \Omega) := \sum_{\alpha \in R_\epsilon(\Omega)} E^\epsilon (u^\epsilon, \epsilon U + \alpha) \;,
	\end{equation}
	where $R_\epsilon(\Omega)$ is defined as
	\begin{equation}\label{eqn:R_eps}
		R_\epsilon(\Omega) := \{\alpha = \alpha_1 v_1 + \alpha_2 v_2  \: : \: \alpha_1,\alpha_2 \in \epsilon \mathbb{Z} \quad \mbox{and} \quad
		\epsilon U + \alpha \ssubset \Omega \} \,.
	\end{equation}
	We use the usual convention that $A \ssubset B$ means $\overline{A} \subset B$.
	
	We note that $\Omega$ need not be an open set for $E^{\epsilon,\eta}(u^\epsilon,\Omega)$ to be well-defined, and no regularity
	is needed for $\partial \Omega$. However, our theorem on the existence of an effective energy (Theorem \ref{thm:effective-energy}) requires that $\Omega$ be a Lipschitz domain.
	
	We wrote above that $u^\epsilon$ should be defined at lattice nodes in $\Omega$ (see \eqref{eqn:admissible-eps}); however, as noted earlier we also want to view $u^\epsilon$ as a piecewise linear function. Therefore our full definition of an admissible deformation is as follows: 
	
	\begin{definition} \label{defn:admissible-deformation}
		An admissible deformation is a pair $(u^\epsilon, \tilde{u}^\epsilon)$ such that
		\begin{enumerate}
			\item[(a)] $u^\epsilon$ belongs to $\mathcal{A}_\epsilon (\Omega)$, i.e. it is a deformation defined at all nodes of the scaled
			lattice that lie in $\Omega$;
			\item[(b)] $\tilde{u}^\epsilon \in H^1(\Omega)$ is the restriction to $\Omega$ of a piecewise linear function
			obtained by applying our piecewise linearization scheme to $u^\epsilon$, which is defined at nodes of the scaled lattice.
		\end{enumerate}
	\end{definition}
	For the Kagome and Rotating Squares examples (as well as for the other examples considered in this paper) we have $\tilde{u}^\epsilon = u^\epsilon$ at all nodes of the scaled lattice that lie in $\Omega$; moreover, $\tilde{u}^\epsilon$ agrees with our piecewise linearization of $u^\epsilon$ at all vertices of the triangulation where the piecewise linearization of $u^\epsilon$ is fully determined. In practice we shall drop the tilde, writing $u^\epsilon$ instead of $\tilde{u}^\epsilon$ since the two functions agree wherever they are both well-defined. The presence of $|D u^\epsilon |_{L^2(\epsilon U + \alpha)}$ on the right hand side of the scaled lower bound \eqref{eqn:unit-cell-lower} makes it natural to work in the function space $H^1(\Omega)$. Therefore, $u^\epsilon$ will henceforth denote a function in $H^1(\Omega)$ that represents a piecewise linearization of some deformation in $\mathcal{A}_\epsilon(\Omega)$.
	
	\subsection{Existence of an effective energy} \label{subsec:eff-energy}
	Since our results on the effective energy use the notion of $\Gamma$-convergence, we start by defining what this means in the present context.
	Here and throughout the paper, the notation $u^\epsilon \rightharpoonup u$ means that $\{ u^\epsilon \}$ remains uniformly
	bounded in $H^1(\Omega)$ and converges \emph{weakly} to $u$.
	
	\begin{definition}[$\Gamma$-convergence] \label{defn:gamma-convergence}
		We say that the family of discrete functionals $\{E^\epsilon(u^\epsilon,\Omega)\}$ $\Gamma$-converges to a functional $E_{\text{eff}}(u,\Omega)$ (with respect to the weak topology of $H^1(\Omega)$) if
		\begin{enumerate}[(i)]
			\item for every admissible sequence $\{u^\epsilon\}_{\epsilon > 0}$ with $u^\epsilon \rightharpoonup u$ in
			$H^1(\Omega)$, we have
			\begin{equation*}
				\liminf_{\epsilon \rightarrow 0} E^\epsilon(u^\epsilon, \Omega)   \geq E_{\text{eff}}(u,\Omega) \, , \quad \mbox{and}
			\end{equation*}
			\item for every $u \in H^1(\Omega)$, there is an admissible sequence $\{u^\epsilon\}_{\epsilon > 0}$ such that
			$u^\epsilon \rightharpoonup u$ in $H^1(\Omega)$ and
			\begin{equation*}
				\lim_{\epsilon \rightarrow 0} E^\epsilon(u^\epsilon, \Omega)  = E_{\text{eff}}(u,\Omega) \, .
			\end{equation*}
		\end{enumerate}
	\end{definition}
	
	Since our unit cell energies $E^{\epsilon,\eta}(u^\epsilon,\epsilon U)$ for the Kagome and Rotating Squares metamaterials satisfy \eqref{eqn:unit-cell-eps-periodicity}-\eqref{eqn:unit-cell-lower}, Theorem 2.11 and Lemma 2.16 of \cite{li2025effective} give the following result concerning the existence of an effective energy for each of these systems:
	\begin{theorem}\label{thm:effective-energy}
		Let $E^{\epsilon,\eta} (u^\epsilon, \Omega)$ be our discrete energy for either the Kagome or the Rotating Squares metamaterial, for any $\eta > 0$. Then for any bounded, Lipschitz domain $\Omega$, $E^{\epsilon,\eta}(u,\Omega)$ $\Gamma$-converges in $H^1(\Omega)$ as $\epsilon \rightarrow 0$ (with respect to the weak topology of $H^1(\Omega)$) to an effective energy of the form
		\begin{equation}\label{eqn:effective-energy}
			E_{\text{eff}}^\eta(u, \Omega) = \int_\Omega \overline{W}^\eta(D u)\: dx \, .
		\end{equation}
		Moreover, the effective energy density $\overline{W}^\eta (\lambda)$ is independent of the domain $\Omega$, and it has the following variational characterization:
		\begin{equation}\label{eqn:effective-energy-density}
			\overline{W}^\eta(\lambda) = \inf_{k \in \mathbb{N}} \inf_{\psi \in \mathcal{A}^{\#}(kU)} \frac{1}{k^2 |U|}
			\sum_{\alpha_1, \alpha_2=0}^{k-1} E^\eta(\lambda x + \psi,U+\alpha_1 v_1 + \alpha_2 v_2) \, .
		\end{equation} 
		where
		\begin{equation} \label{eqn:defn-kU}
			kU = \bigcup_{\alpha_1, \alpha_2=0}^{k-1} \Big(U+ \alpha_1 v_1 + \alpha_2 v_2 \Big)
		\end{equation}
		and $\mathcal{A}^{\#}(kU)$ is the set of deformations defined at all nodes of our lattice that are $k$-periodic (that is, deformations $\psi$ such that $\psi(x) = \psi(x+k \alpha_1 v_1 + k\alpha_2 v_2)$ for any $\alpha_1,\alpha_2 \in \mathbb{Z}$).
	\end{theorem} 
	
	Our spring systems have the property that their energies are \emph{frame-indifferent}, in other words that $E^\eta (Ru, U) = E^\eta(u,U)$ for any orientation-preserving rotation $R$. It follows from \eqref{eqn:effective-energy-density} that the effective energy density $\overline{W}^\eta$ is also frame-indifferent (that is, $\overline{W}^\eta(R\lambda) = \overline{W}^\eta (\lambda) $ for any orientation-preserving rotation $R$). 
	
	In Theorem \ref{thm:effective-energy} we have not imposed any boundary condition. A similar result holds, however, in the presence of a (suitably-imposed) Dirichlet-type boundary condition; see Theorem 2.13 of \cite{li2025effective}. It is worth noting that the effective energy density remains the same, regardless of whether or not a boundary condition is imposed.
	
	\subsection{Our main results}\label{subsec:main-results}
	
	From the definition of $\Gamma$-convergence, we see that the macroscopic deformations with effective energy $0$ represent \emph{soft modes}, in the sense that they can be
	achieved asymptotically (as $\epsilon \rightarrow 0$) by deformations ($u^\epsilon$) whose area-averaged energy tends to zero. We use the term ``area-averaged energy'' because one easily checks using \eqref{eqn:elasticity-scaling} that the scaled energy $E^{\epsilon,\eta}(u^\epsilon,\Omega)$ is essentially the sum of (i) $|\Omega|$ times the average energy of all the springs in $\Omega$, and (ii) $1/\eta$ times the areas of the triangles in the penalization term whose orientations are reversed by $u^\epsilon$. 
	
	Our main accomplishment in this paper is a unified proof that for both the Kagome and Rotating Squares metamaterials, the soft modes are compressive conformal maps (see Theorem \ref{thm:compressive-conformal}). Our approach is to consider the zero set of the effective energy density $\overline{W}^\eta(\lambda)$ for $\eta$ small enough. In fact, using the geometry of the two metamaterials, we are able to provide a lower bound on $\overline{W}^\eta(\lambda)$ for small $\eta$, which vanishes only at isotropic compressions:
	\begin{theorem}\label{thm:lower-bound}
		For the Kagome and Rotating Squares metamaterials, there exists a threshold $\eta_0 > 0$ for the penalty constant such that when $0 < \eta < \eta_0$, the effective energy density $\overline{W}^{\eta}(\lambda)$ is lower bounded by
		\begin{equation}\label{eqn:eff-lower-bd}
			\begin{aligned}
				\overline{W}^\eta(\lambda) & \geq 
				C \Big[(\lambda_1 - \lambda_2)^2 + (\lambda_1 - 1)_+^2 + (\lambda_2 - 1)_+^2\Big], & \det(\lambda) \geq 0,\\
				\overline{W}^\eta(\lambda) & \geq  C \Big[(\lambda_1 + \lambda_2)^2 + (\lambda_1 - 1)_+^2 + (\lambda_2 - 1)_+^2\Big], & \det(\lambda) < 0,
			\end{aligned}
		\end{equation}
		where the constant $C$ is independent of $\lambda$ and $\eta$. Here $\lambda_1$ and $\lambda_2$ are the principal stretches associated with the $2 \times 2$ matrix $\lambda$, in other words the eigenvalues of $(\lambda^T \lambda)^{1/2}$, and we use the notation $(x-1)_+^2 = (\max\{x-1,0\})^2$.
	\end{theorem}
	
	The proof of Theorem \ref{thm:lower-bound} is presented in Section \ref{sec:lower-bd}. As a consequence of Theorem \ref{thm:lower-bound}, we have the following result about the zero set of the effective energy density:
	\begin{proposition} \label{prop:isotropic-compression}
		For the Kagome and Rotating Squares metamaterials, when $0 < \eta < \eta_0$, the effective energy density functional $\overline{W}^{\eta}(\lambda)$ vanishes if and only if $\lambda$ is an isotropic compression, i.e. $\lambda = cR$ with $0 \leq c \leq 1$ and $R \in SO(2)$.
	\end{proposition}
	\begin{proof}
		One direction is straightforward since by applying Theorem \ref{thm:lower-bound}, the effective energy density $\overline{W}^{\eta}(\lambda)$ vanishes only when $\det \lambda \geq 0$ and the principal strains satisfy $\lambda_1 =  \lambda_2 \leq 1$. Since $\lambda_1 = \lambda_2$, the singular decomposition of $\lambda$ becomes $\lambda = \lambda_1 UV^T$ with $U,V\in O(2)$. Since $\det \lambda \geq 0$, we have $\det UV^T \geq 0$ and $UV^T \in SO(2)$. Therefore, $\lambda = cR$ for some $R \in SO(2)$ and $c=\lambda_1 \in [0,1]$.
		
		For the other direction, we observe that for any $\lambda =cR$ with $0 \leq c \leq 1$ and $R \in SO(2)$, we can choose a one-periodic mechanism of the metamaterial that achieves the compression ratio $c$, and then apply a global rotation $R$. This rotated one-periodic mechanism has the form $u(x) = cR x + \psi(x)$ for some 1-periodic $\psi$. Therefore we have 
		\begin{equation*}
			0 \leq \overline{W}^\eta(cR) \leq \frac{1}{|U|} E^\eta(cRx + \psi, U) = 0,
		\end{equation*}
		which completes our proof.
	\end{proof}
	
	Having characterized the zero set of $\overline{W}^\eta(\lambda)$ explicitly, we now turn to the zero-energy set of the effective energy$E_\text{eff}^\eta(u,\Omega)$, which we claim to be the class of compressive conformal maps. Here we are using the following terminology: 
	\begin{definition}[Compressive conformal maps]
		A deformation $u: \Omega \subset \mathbb{R}^2 \rightarrow \mathbb{R}^2$ is a compressive conformal map if $u(x_1,x_2) = (u_1(x_1,x_2), u_2(x_1,x_2))$ satisfies the following conditions
		\begin{itemize}
			\item the complex function $f(z) = u_1(x_1,x_2) + iu_2(x_1,x_2)$ with $z = x_1+ix_2$ is analytic;
			
			\item the complex function $f(z)$ has its derivative bounded by 1, i.e. $|f'(z)| \leq 1$ for any $z = x_1+ix_2$ with $(x_1,x_2) \in \Omega$.
		\end{itemize}
	\end{definition}
	
	\begin{remark}
		A compressive conformal map by our definition is slightly different from the usual complex-variable-based notion of a conformal map, which requires $f'(z) \neq 0$. A compressive conformal map under our definition can have zero gradient, i.e. $|f'(z)| = 0$ is allowed.
	\end{remark}
	
	Now we state and prove our main theorem:
	\begin{theorem}\label{thm:compressive-conformal}
		For the Kagome and Rotating Squares metamaterials, when $0 < \eta < \eta_0$, the effective energy $E^\eta_\text{eff}(u,\Omega)$ vanishes if and only if $u$ is a compressive conformal map.
	\end{theorem}
	
	\begin{proof}
		The key idea of the proof is to show that $E_\text{eff}^\eta (u,\Omega) = 0$ is equivalent to $D u$ satisfying the Cauchy-Riemann equation as a distribution. For a deformation $u \in H^1(\Omega)$, by applying Proposition \ref{prop:isotropic-compression}, the effective energy $E_\text{eff}^\eta (u,\Omega)$ vanishes if and only if there exists $c(x), R(x) \in L^2$ such that
		\begin{equation}\label{eqn:grad-u}
			D u = c(x) R(x) \quad \text{a.e. in }\Omega, \qquad 0 \leq c(x)\leq 1, \quad R(x) \in SO(2).
		\end{equation}
		We claim that \eqref{eqn:grad-u} is equivalent to $u$ being compressive conformal. To see this, we express the components of $D u = c(x) R(x)$ explicitly in terms of an angle function $\alpha (x)$ defined almost everywhere, i.e.
		\begin{equation}
			D u = \begin{pmatrix}
				\partial_{x_1} u_1 & \partial_{x_2} u_1\\
				\partial_{x_1} u_2 & \partial_{x_2} u_2
			\end{pmatrix} = \begin{pmatrix}
				c(x) \cos \alpha(x) & -c(x) \sin \alpha(x)\\
				c(x) \sin \alpha(x) & c(x) \cos \alpha(x)
			\end{pmatrix} = c(x) R(x),
		\end{equation}
		which implies the Cauchy-Riemann equations for the complex function $f(z) = u_1 + iu_2$, i.e.
		\begin{equation} \label{eqn:CR-eq}
			\partial_{x_1} u_1 - \partial_{x_2} u_2 = \partial_{x_1} u_2 +\partial_{x_2} u_1 = 0 \qquad \text{a.e. in } \Omega.
		\end{equation}
		It follows that $\Delta u_1 = \Delta u_2 = 0$ in the distributional sense. To explain why, let us show that $\Delta u_1 =0$. For any $\varphi \in C_0^\infty(\Omega)$, we have
		\begin{equation*}
			\int_\Omega u_1 \Delta \varphi \; dx = -\int_\Omega \partial_{x_1} u_1 \partial_{x_1} \varphi + \partial_{x_2} u_1 \partial_{x_2} \varphi \; dx = -\int_\Omega \partial_{x_2} u_2 \partial_{x_1} \varphi - \partial_{x_1} u_2 \partial_{x_2} \varphi \; dx= 0.
		\end{equation*}
		Therefore, by Weyl's Lemma (stated below, just after the end of this proof), we obtain the smoothness of $u_1$ (and similarly that of $u_2$). As a consequence of this smoothness, it follows that $f(z) = u_1+iu_2$ is analytic since \eqref{eqn:CR-eq} is the Cauchy-Riemann equation for $f(z)$. Finally, to show the compressive part, we compute that $|f'(z)| = c(x) \leq 1$. Thus, we have shown that \eqref{eqn:grad-u} is equivalent to $u$ being compressive conformal.
	\end{proof}
	
	The preceding argument used the following well-known result (see e.g. Corollary 2.2.1 in \cite{jost2012partial}):
	\begin{lemma}[Weyl's lemma]
		Let $f:\Omega \rightarrow \mathbb{R}$ be measurable and locally integrable in $\Omega$. Suppose that $\int_\Omega u(x) \Delta \varphi(x) dx = 0$ holds for all $\varphi \in C_0^\infty(\Omega)$. Then $u$ is harmonic and smooth.
	\end{lemma}
	
	The importance of Theorem \ref{thm:compressive-conformal} is that when we view the (Kagome or the Rotating Squares) metamaterial as a continuum object, it can sustain any compressive conformal map as an elastic deformation with negligible energy (in the sense that the area-averaged energy tends to $0$). Furthermore, the theorem allows us to analyze sequences of deformations $u^\epsilon$ with vanishing energy, which must then converge weakly to a compressive conformal map. We state this result as follows.
	\begin{corollary}\label{coro:compressive-conformal}
		Consider the Kagome or the Rotating Squares metamaterial and fix $\eta$ in the range $0 < \eta < \eta_0$. If a sequence of admissible deformations $u^\epsilon \in H^1(\Omega)$ satisfies the following conditions:
		\begin{itemize}
			\item there exists a uniform $L^\infty$ bound on $u^\epsilon$, i.e. there exists a constant $D_0$ such that
			\begin{equation}\label{eqn:uniform-bound-u}
				|u^\epsilon|_{L^\infty(\Omega)} \leq D_0;
			\end{equation}
			
			\item there exists a larger domain $\Omega_c$ with $\overline{\Omega} \ssubset \Omega_c$ such that $u^\epsilon$ can be extended to $\Omega_c$ as a piecewise linear function on the corresponding mesh with a bounded energy, i.e. 
			\begin{equation}\label{eqn:energy-bound}
				E^{\epsilon, \eta}(u^\epsilon,\Omega_c) \leq E_0
			\end{equation}
			for some $E_0 > 0$;
			
			\item on the original set $\Omega$, the energy $E^{\epsilon,\eta}(u^\epsilon, \Omega)$ is small and vanishes in the limit as $\epsilon \rightarrow 0$, i.e.
			\begin{equation}\label{eqn:compressive-conformal-condition}
				\liminf_{\epsilon \rightarrow 0} E^{\epsilon,\eta}(u^\epsilon, \Omega) = 0,
			\end{equation}
		\end{itemize}
		then deformations $u^\epsilon$ lie in a weakly compact subset of $H^1(\Omega)$, and all their weak limit points are compressive conformal maps $u$.
	\end{corollary}
	
	\begin{proof}
		To show weak compactness, it suffices to show that $u^\epsilon$ remains bounded in $H^1(\Omega)$ as $\epsilon \rightarrow 0$. Since $\overline{\Omega} \ssubset \Omega_c$, the set $\Omega$ is contained within the union of $\epsilon$-scale unit cells that lie entirely in $\Omega_c$ for sufficiently small $\epsilon$. Therefore, we can control the $H^1$ norm of $u^\epsilon$ on $\Omega$ by the energy on $\Omega_c$, by making use of the lower bound \eqref{eqn:unit-cell-lower}:
		\begin{equation*}
			C_2 \Big(|D u^\epsilon|^2_{L^2(\Omega)} - D_2 |\Omega|\Big) \leq \sum_{\alpha \in R_\epsilon(\Omega_c)} E^{\epsilon, \eta}(u^\epsilon, \epsilon U + \alpha) = E^{\epsilon, \eta}(u^\epsilon, \Omega_c) \leq E_0.
		\end{equation*}
		Combining the above bound with \eqref{eqn:uniform-bound-u}, we obtain the uniform $H^1$ bound for $u^\epsilon$.
		
		The remaining assertions of the corollary follow directly from the $\Gamma$-convergence of $E^{\epsilon,\eta}$ to $E^\eta_\text{eff}$. Indeed, consider any weakly-convergence subsequence (still denoted $u^\epsilon$, for simplicity) and let $u$ be its limit. To see why $u$ must be a compressive conformal map, we use part (i) of 
		Definition \ref{defn:gamma-convergence}  to get 
		\begin{align*}
			0 \leq E_\text{eff}^\eta(u, \Omega) &\leq \liminf_{\epsilon \rightarrow 0} E^{\epsilon, \eta} (u^\epsilon, \Omega).
		\end{align*}
		Thus, the effective energy $E^\eta_\text{eff}(u,\Omega)$ vanishes by assumption \eqref{eqn:compressive-conformal-condition}, and  we conclude using Theorem \ref{thm:compressive-conformal} that $u$ is a compressive conformal map.
	\end{proof}
	
	\begin{remark}
		The first condition \eqref{eqn:uniform-bound-u} serves to rule out sequences with unbounded translations -- for example, $u^\epsilon(x) = x + 1/\epsilon$, which has scaled energy $0$ (so that \eqref{eqn:energy-bound} and \eqref{eqn:compressive-conformal-condition} are satisfied) though it doesn't stay bounded in $H^1(\Omega)$. One could of course get compactness by imposing a weaker condition -- for instance, by controlling the average of $u^\epsilon$. Here we have chosen a very simple condition, which is sufficient for our purposes since the soft modes encountered in practice are bounded in $L^\infty$ norm.
		
		The second condition \eqref{eqn:energy-bound} is also necessary to ensure compactness, since the energy $E^{\epsilon,\eta}(u^\epsilon,\Omega)$ only captures the behavior of $u^\epsilon$ on the $\epsilon$-scale unit cells that lie entirely in $\Omega$; thus, it neglects information on an $\epsilon$-scale boundary layer. Let us use a 1D lattice example to help visualize this phenomenon.  
		Consider the interval $\Omega = (0,1)$ and let the lattice size be $\epsilon = 1/N$ for $N \in \mathbb{N}$. We take the (scaled) lattice nodes to be $x_i = i\epsilon$ for $i=0,1,\dots,N$, and we take the (scaled) unit cell to be of size $\epsilon$. To align with our discrete energy framework, the admissible deformation $u^\epsilon$ is defined only at the lattice nodes $x_i$. We use the standard Hookean spring energy and, for simplicity, omit the penalty term. Let the resulting energy be denoted by $E^\epsilon(u^\epsilon,\Omega)$. The key point is that $E^\epsilon(u^\epsilon,\Omega)$ only sees springs that lie inside $\Omega = (0,1)$; in particular, the nodal values of $u^\epsilon$ at the ends $x_0 = 0$ and $x_N = 1$ are not considered. If we choose $u^\epsilon(x)$ to be the piecewise linearization of
		\begin{equation*}
			v^\epsilon(x_i)= \begin{cases}
				x_i + 1, & x_i =0 \text{ or }1,\\
				x_i, & x_i \in (0,1),
			\end{cases}
		\end{equation*}
		then the discrete energy $E^\epsilon(u^\epsilon,\Omega) = 0$, while the $H^1$ norm of $u^\epsilon$ is not uniformly bounded since $|(u^\epsilon)'|_{L^2(\Omega)}$ $=$ $(2/\epsilon)^{1/2}$ blows up as $\epsilon \rightarrow 0$. The point of \eqref{eqn:energy-bound} is to eliminate such examples.
	\end{remark}
	
	
	\section{A geometric argument: periodic mechanisms are isotropic compressions}\label{sec:geo-argument}
	
	This section's main goal is Proposition \ref{lemma-geometric}. It characterizes the macroscopic deformations achievable
	by periodic mechanisms of the Kagome or Rotating Squares metamaterials, using an argument that applies equally
	to these two systems. While its assertion was already known for the Rotating Squares example, it was \emph{not}
	previously known for the Kagome metamaterial (which has infinitely many mechanisms, for which there is no known
	classification).
	
	The proof of Proposition \ref{lemma-geometric} is relatively simple and geometric. It captures -- in the simplest possible
	setting -- how we will use the symmetry and structure of our metamaterials to lower-bound their effective
	energies.
	
	Let us explain why studying periodic mechanisms is related to estimating the effective energy.
	Our variational characterization \eqref{eqn:effective-energy-density} of the effective energy density
	$\overline{W}^\eta$ involves only deformations of the form $u(x) = \lambda x + \psi(x)$ where $\psi$ is $k$-periodic
	for some $k$. The characterization involves energy minimization, and deformations with energy $0$ are automatically
	minimizers (since our energies are nonnegative). So it is quite natural, as a first step, to consider the existence
	of zero-energy deformations of the form $u(x) = \lambda x + \psi(x)$ where $\psi$ is $k$-periodic. These are
	precisely the periodic mechanisms (see Definition \ref{def:periodic-mechanism}).
	
\begin{figure}[!htb]
		\centering
		\includegraphics[width=0.7\linewidth]{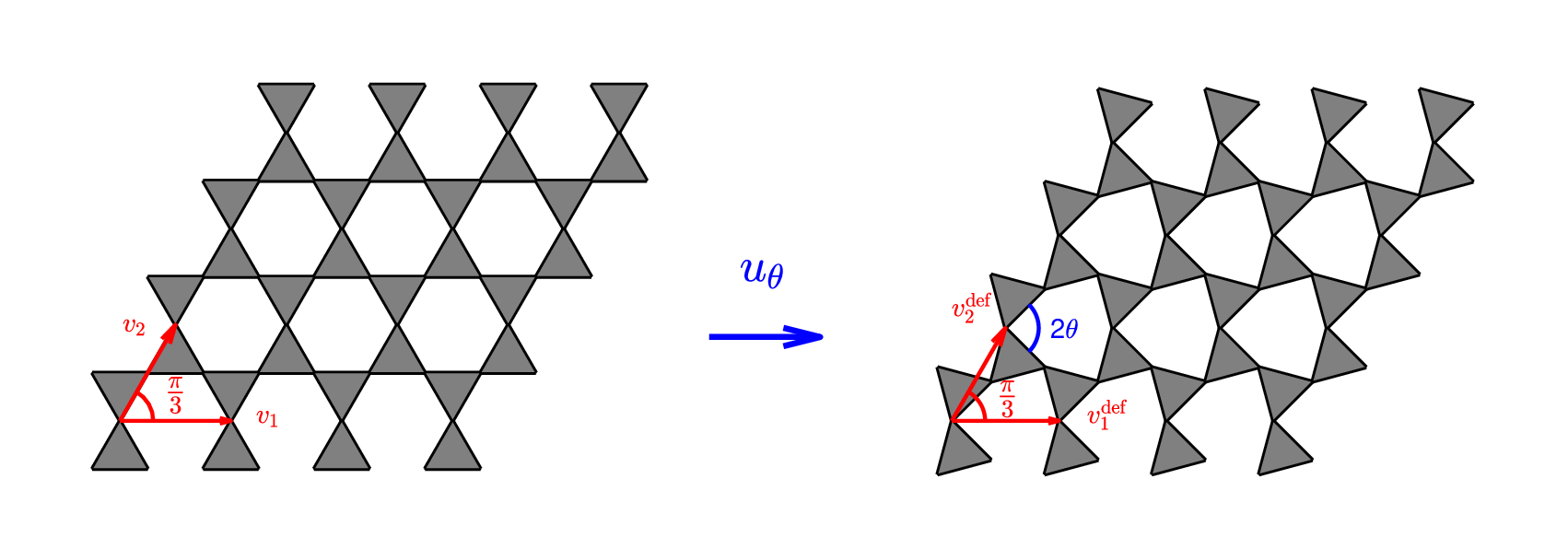}
		\caption{The ``conformality'' of the one-periodic mechanism of the Kagome metamaterial: 
			the angle between the lattice vectors $v_1$ and $v_2$ in the reference state is the same as the 
			angle between $v_1^{\text{def}}$ and $v_2^{\text{def}}$ in the deformed state (in both cases it is $\pi/3$).}
		\label{fig:kagome-conformal-mechanism}
	\end{figure}

	A key feature of our analysis is that it doesn't require any formula for the mechanism under consideration. To
	provide some intuition about what is at stake, let us briefly discuss the one-periodic mechanism of the Kagome
	lattice in a formula-free way. Figure \ref{fig:kagome-conformal-mechanism} shows the reference and deformed
	configurations for a particular value of its parameter $\theta$. Both are periodic lattices of springs, with
	lattice vectors $v_1,v_2$ for the former and $v_1^\text{def},v_2^\text{def}$ for the latter. The macroscopic
	action of the mechanism is the linear transformation $\lambda_\theta$ such that $\lambda_\theta v_i = v_i^\text{def}$
	for $i=1,2$ (see Lemma \ref{lemma:atilde-vs-lambda}). To prove that $\lambda_\theta$ is an isotropic compression, we 
	will have to show that the parallelogram
	with sides $v_1^\text{def}, v_2^\text{def}$ is a scaled copy of the one with sides $v_1, v_2$ (with a scale
	factor in $[0,1]$). The arguments we use for this will apply to $k$-periodic mechanisms for any $k$, and will be
	formula-free (in the sense that they don't require a formula for the mechanism).
	
	We turn now to a more careful statement of this section's goal. Let us start with the definition of a periodic mechanism.
	While we sometimes use the term ``mechanism'' for a parametrized family of energy-free deformations, we use it here
	for a \emph{particular} energy-free deformation (which might or might not come from a parametrized family):
	
	\begin{definition} \label{def:periodic-mechanism}
		Let $E^\eta(u,U)$ be the unit cell energy of the Kagome or Rotating Squares metamaterial, as discussed in
		Section \ref{subsec:lattice-setup}. We say that $u(x) = \lambda x + \psi(x)$ is a $k$-periodic mechanism achieving
		macroscopic deformation $\lambda$ if
		\begin{enumerate}
			\item[(a)] $\lambda$ is a $2 \times 2$ matrix,
			\item[(b)] $\psi$ is an $\mathbb{R}^2$-valued function defined at all nodes of the lattice
			which is $k$-periodic in the sense that $\psi(x+kv_1) = \psi(x+kv_2) = \psi(x)$ (where $v_1,v_2$ are
			the lattice vectors of the reference configuration), and
			\item[(c)] $E^\eta (u,U+\alpha_1 v_1 + \alpha_2 v_2) = 0$ \ for every $(\alpha_1,\alpha_2) \in \mathbb{Z}^2$.
		\end{enumerate}
	\end{definition}
	Since $E^\eta$ consists of spring energies and a penalization term, (c) amounts the condition that $u$ doesn't
	stretch or compress any of the springs and it preserves the orientations of the triangles that enter the
	penalization term. Since $\psi$ is $k$-periodic and $E^\eta \geq 0$, (c) is equivalent to the assertion that
	the \emph{average energy on $kU$} vanishes, in other words that $\overline{E^\eta}(\lambda, \psi, kU) = 0$ where
	\begin{equation}\label{eqn:avg-energy}
		\overline{E^\eta}(\lambda, \psi, kU) := \frac{1}{k^2 |U|}
		\sum_{\alpha_1, \alpha_2=0}^{k-1} E^\eta(\lambda x + \psi,U+\alpha_1 v_1 + \alpha_2 v_2) \,.
	\end{equation}
	(Note that in the definition \eqref{eqn:effective-energy-density} of the effective energy density
	$\overline{W}^\eta (\lambda)$,
	it is precisely $\overline{E^\eta}(\lambda,\psi,kU)$ that is being minimized with respect to $k$ and $\psi$.)
	
	We now state the main result of this section:
	\begin{proposition}\label{lemma-geometric}
		For either the Kagome or Rotating Squares metamaterial, if $u(x) = \lambda x + \psi(x)$ is a
		$k$-periodic mechanism then its macroscopic deformation $\lambda$ must be an isotropic compression, i.e.
		$\lambda = c R$ for some $0 \leq c \leq 1$ and $R \in SO(2)$.
	\end{proposition}
	
	Our plan for the section is as follows. We start in Section \ref{subsec:k-periodic-def} by introducing some notation
	that lets us treat both the Kagome and Rotating Squares lattices at the same time. The utility of that notation
	is demonstrated in Lemma \ref{lemma:atilde-vs-lambda}, which identifies the
	analogues in the $k$-periodic setting of the vectors $v_i^\text{def}$ shown in
	Figure \ref{fig:kagome-conformal-mechanism}. Then we present the proof of Proposition \ref{lemma-geometric} in
	Section \ref{subsec:proof-of-geometric-lemma}.
	
	\subsection{Deformations of the form $u(x) = \lambda x + \psi$ with $k$-periodic $\psi$} \label{subsec:k-periodic-def}
	In Section \ref{subsec:lattice-setup} we introduced our energies $E^\eta(u,U)$ for the Kagome and Rotating-Squares examples. In doing
	so, we identified certain basic sets of nodes and springs (whose translates by lattice vectors were disjoint and
	gave all the lattice's nodes and springs). Those choices were not arbitrary -- they were constrained by the
	requirement that $E^\eta(u,U)$ satisfy the basic upper and lower bounds
	\eqref{eqn:unscaled-upper-bound}--\eqref{eqn:unscaled-lower-bound} required by our theorem on the existence of an
	effective energy. For this section and Section \ref{sec:lower-bd}, however, it is convenient to use different choices of
	the basic sets, which we now present. (Throughout this subsection we consider deformations of the form
	$u(x) = \lambda x + \psi(x)$ where $\psi$ is $k$-periodic, however \emph{we do not assume that $u$ is a mechanism}. This is
	important, since the framework we develop here will also be used in Section \ref{sec:lower-bd}.)
	
\begin{figure}[!htb]
		\begin{minipage}{.45\linewidth}
			\centering
			\subfloat[]{\label{fig:kagome-unit-1-wo}\includegraphics[scale=.35]{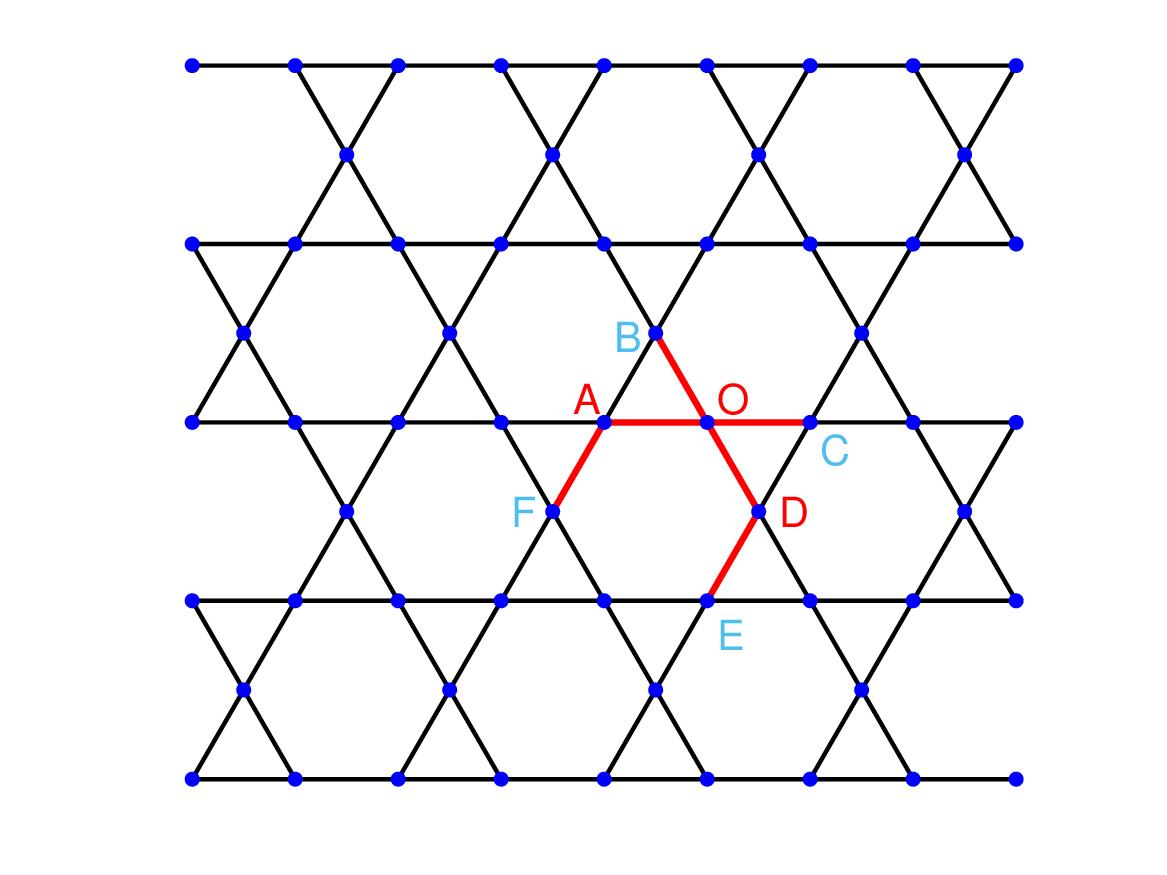}}
		\end{minipage}
		\begin{minipage}{.45\linewidth}
			\centering
			\subfloat[]{\label{fig:kagome-unit-2}\includegraphics[scale=.35]{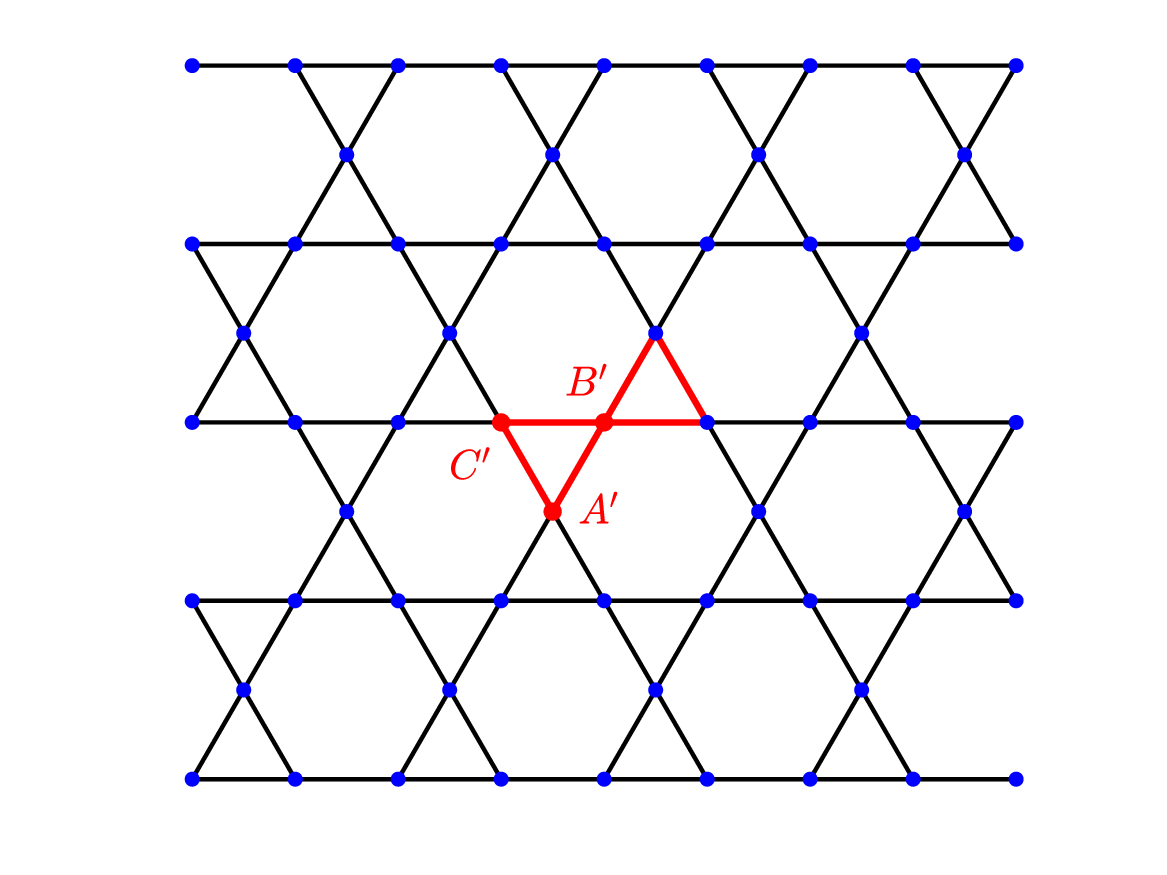}}
		\end{minipage}
		\caption{(a) The basic sets of nodes and springs from Figure \ref{fig:kagome-unit}; (b) our new basic set of nodes (labeled
			$A$, $B$, and $C$) and springs (shown in red).}
	\end{figure}

	We start with the Kagome metamaterial. Figure \ref{fig:kagome-unit-1-wo} shows the basic sets of nodes 
	and springs discussed in Section \ref{subsec:lattice-setup}. Our new choice is shown in Figure \ref{fig:kagome-unit-2}: the basic set of nodes is
	labeled $A', B', C'$, and the basic set of springs is shown in red. These choices have, once again, the property
	that their translates by lattice vectors are disjoint and give all the lattice's nodes and
	springs.
	
\begin{figure}[!htb]
		\begin{minipage}{.45\linewidth}
			\centering
			\subfloat[]{\label{fig:kagome-k-period}\includegraphics[scale=.35]{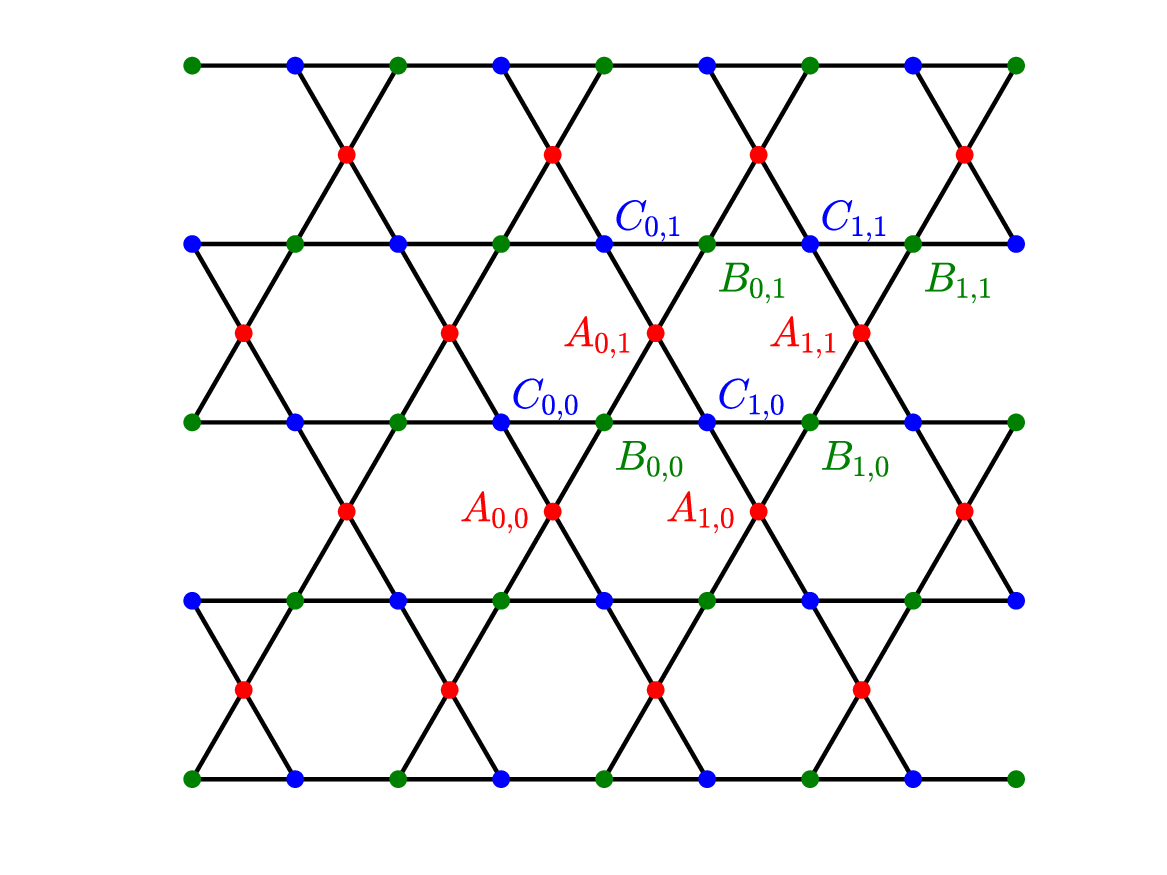}}
		\end{minipage}
		\begin{minipage}{.45\linewidth}
			\centering
			\subfloat[]{\label{fig:kagome-k-unit}\includegraphics[scale=.3]{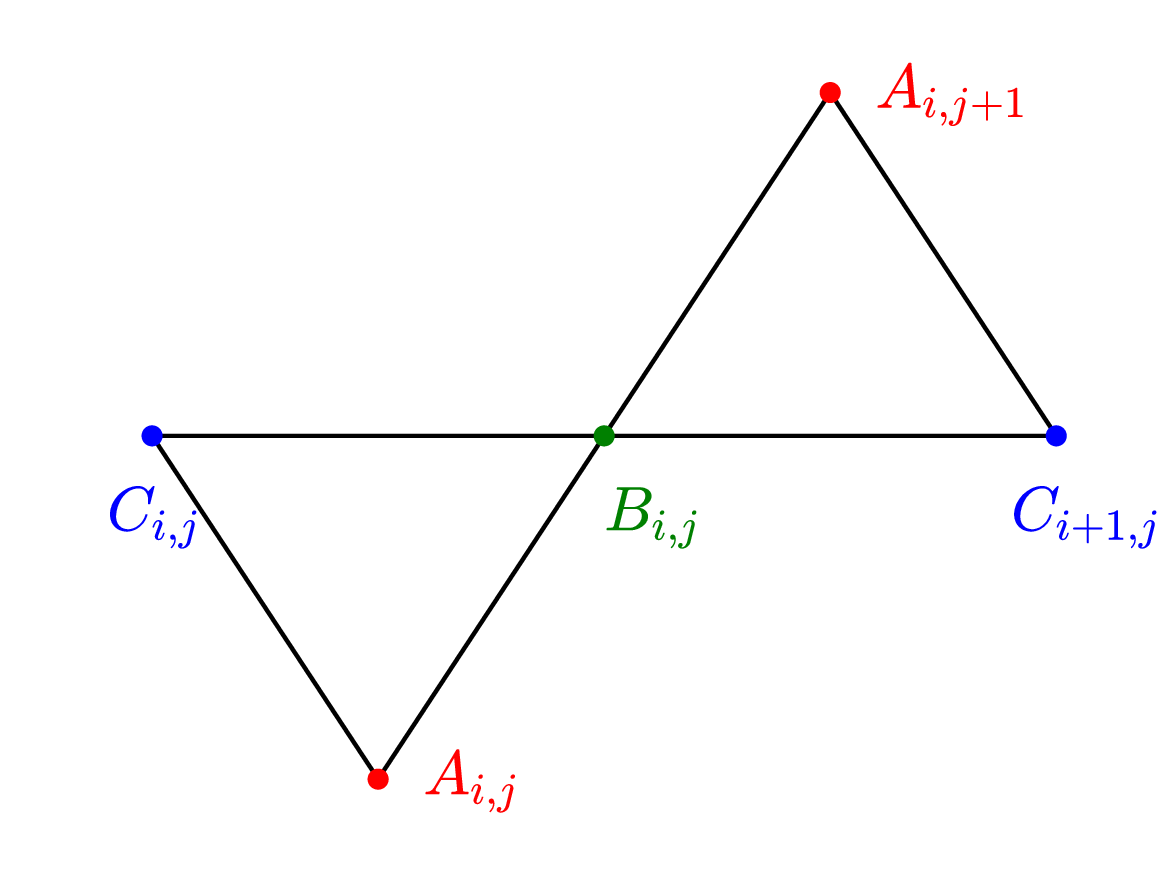}}
		\end{minipage}
		\caption{(a) The Kagome lattice the labeling convention \eqref{eqn:notation-for-translates} for its nodes. Labels
			are given here for a basic set of nodes under equivalence by $2$-periodicity.  (b) The $(i,j)$th pair of
			equilateral triangles.}
	\end{figure}

	For a deformation of the form $u(x) = \lambda x + \psi(x)$ with $\psi$ $k$-periodic, the quantities that interest us
	(for example the energies of its springs) will be $k$-periodic. We will need to average them (for example to calculate
	the average energy \eqref{eqn:avg-energy}). This requires choosing one representative from each of the equivalence
	classes associated with $k$-periodicity. Our new choice of basic set makes this very easy. To average a function
	defined on nodes, for example, it suffices to consider its values at the $3k^2$ nodes
	\begin{equation} \label{eqn:notation-for-translates}
		A_{i,j} = A' + iv_1 + jv_2, \quad B_{i,j} = B' + iv_1 + jv_2, \quad C_{i,j} = C' + iv_1 + jv_2,
	\end{equation}
	where $i,j=0,1,\dots,k-1$ (see Figure \ref{fig:kagome-k-period} for an example with $k=2$).
	This labeling convention for the nodes makes it natural to use a similar convention for their
	images under the deformation $u$; we therefore define
	\begin{equation*}
		\widetilde{A}_{i,j} = u(A_{i,j}), \quad \widetilde{B}_{i,j} = u(B_{i,j}), \quad \widetilde{C}_{i,j} = u(C_{i,j}) .
	\end{equation*}
	
	A similar convention will be used for the equilateral triangles in the Kagome lattice. Our new basic set of springs
	determine two triangles (see Figure \ref{fig:kagome-unit-2}), so in the $k$-periodic setting the average of a function
	defined on triangles is the normalized sum of its values on their translates
	\begin{equation} \label{eqn:translates-of-triangles}
		T_{i,j,1} := \Delta {A}_{ij} {B}_{ij} {C}_{ij} \qquad T_{i,j,2} := \Delta {A}_{i,j+1} {B}_{i,j} {C}_{i+1,j},
	\end{equation}
	for $i,j = 0,\dots,k-1$ (see Figure \ref{fig:kagome-k-unit}). The corresponding triangles in the deformed lattice
	are
	\begin{equation*}
		\widetilde{T}_{i,j,1} := \Delta \widetilde{A}_{i,j} \widetilde{B}_{i,j} \widetilde{C}_{i,j}, \qquad 
		\widetilde{T}_{i,j,2} := \Delta \widetilde{A}_{i,j+1} \widetilde{B}_{i,j} \widetilde{C}_{i+1,j},
	\end{equation*}
	for $i,j = 0,\dots,k-1$. (Note that if $u$ is not a mechanism then these may not be equilateral.)
	
\begin{figure}[!htb]
		\begin{minipage}{.48\linewidth}
			\centering
			\subfloat[]{\label{fig:ref-kagome}\includegraphics[width=0.9\linewidth]{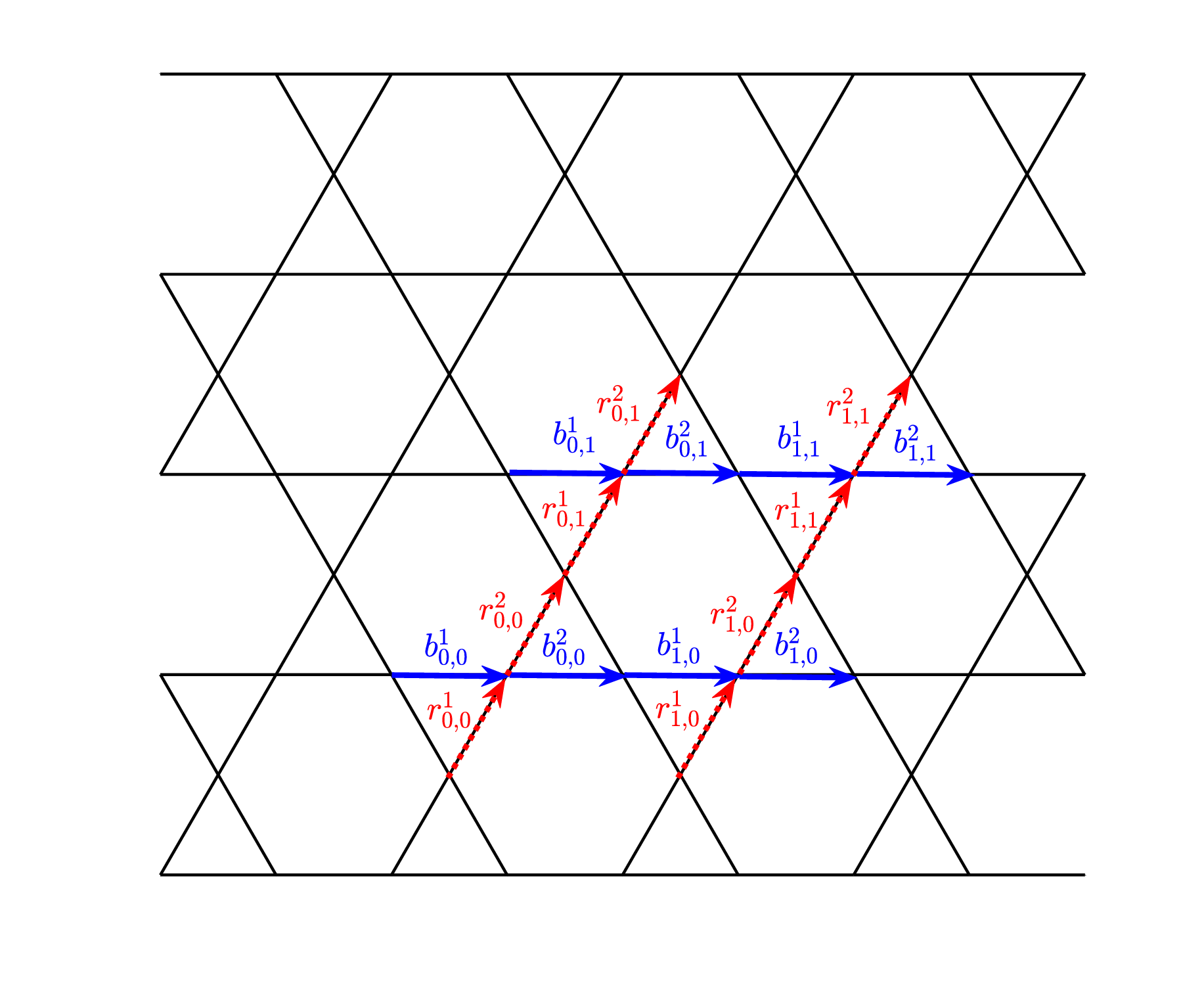}}
		\end{minipage}
		\begin{minipage}{.48\linewidth}
			\centering
			\subfloat[]{\label{fig:def-kagome}\includegraphics[width=0.9\linewidth]{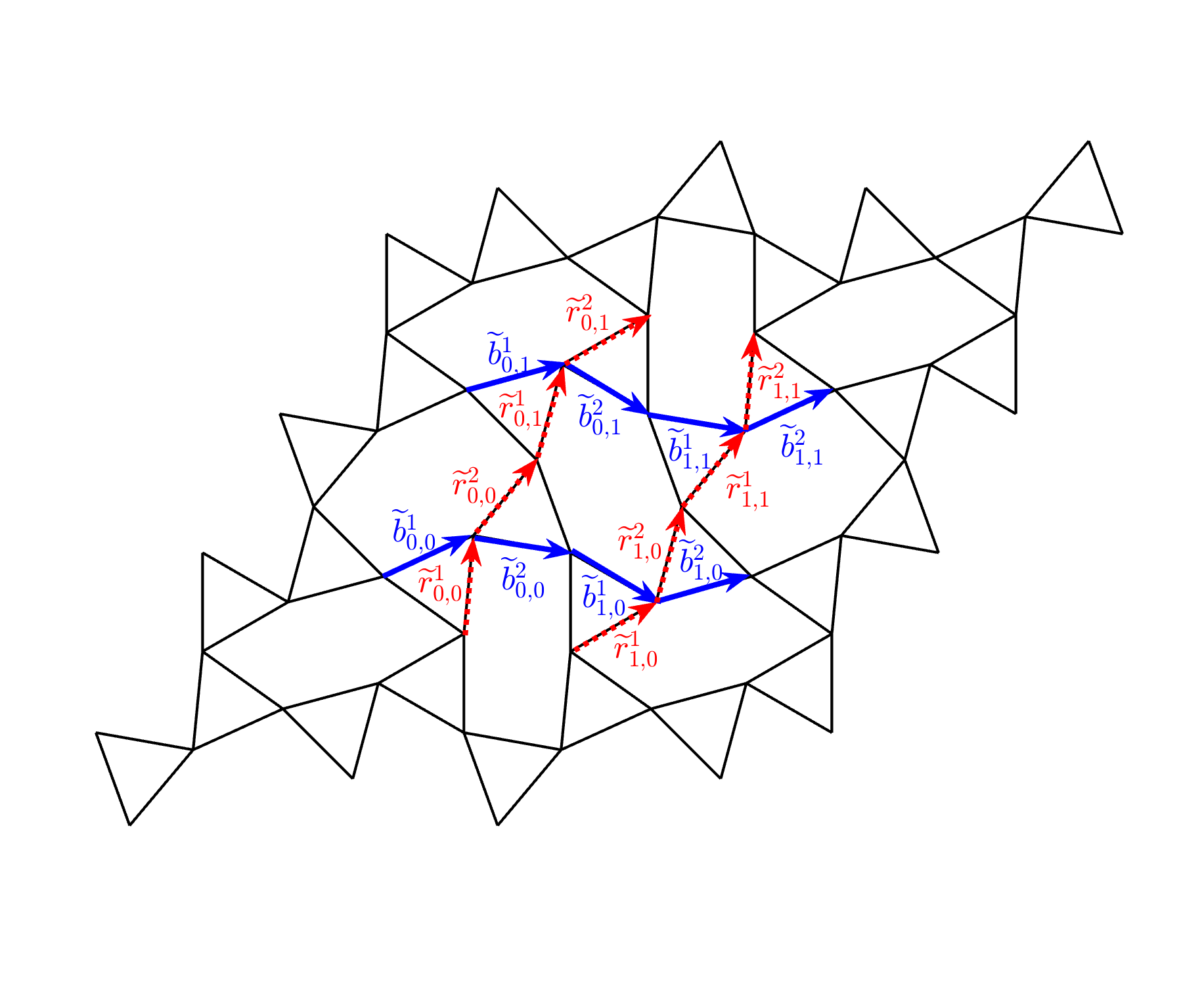}}
		\end{minipage}
		\caption{(a) The Kagome lattice in the reference state; (b) the image of the Kagome lattice by a 
			two-periodic mechanism. The marked vectors are $b_{i,j}^t, r_{i,j}^t$ in the reference state and 
			$\widetilde{b}_{i,j}^t, \widetilde{r}_{i,j}^t$ in the deformed state with $t=1,2$.}
		\label{fig:kagome-2-periodic}
	\end{figure}

	We now introduce names for some of the sides of the triangles. Notice that each of the
	triangles $T_{i,j,1}$ and $T_{i,j,2}$ has a horizontal side and a side in the $60$ degree direction. We
	let $b_{i,j}^1, b_{i,j}^2$ be the horizontal sides of triangles $T_{i,j,1}$ and $T_{i,j,2}$ respectively,
	we let $r_{i,j}^1, r_{i,j}^2$ be the sides in the $60$ degree direction:
	\begin{equation}\label{eqn:br-vec-ref}
		\begin{aligned}
			& b_{i,j}^1 :=\overrightarrow{C_{i,j} B_{i,j}}, \qquad b_{i,j}^2 := \overrightarrow{B_{i,j} C_{i+1,j}},\\
			& r_{i,j}^1 :=\overrightarrow{A_{i,j} B_{i,j}}, \qquad r_{i,j}^2 := \overrightarrow{B_{i,j} A_{i,j+1}}.
		\end{aligned}
	\end{equation}
	As usual in geometry, we view these as vectors and identify them with points in the plane
	(for example, $b_{i,j}^1 = B_{i,j} - C_{i,j}$). An illustration for the case $k = 2$ is shown
	in Figure \ref{fig:ref-kagome}. The corresponding sides of the triangles in the deformed lattice will be called
	\begin{equation}\label{eqn:br-vec-def}
		\begin{aligned}
			& \widetilde{b}_{i,j}^1 :=\overrightarrow{\widetilde{C}_{i,j} \widetilde{B}_{i,j}}, \qquad \widetilde{b}_{i,j}^2 := \overrightarrow{\widetilde{B}_{i,j} \widetilde{C}_{i+1,j}},\\
			& \widetilde{r}_{i,j}^1 :=\overrightarrow{\widetilde{A}_{i,j} \widetilde{B}_{i,j}}, \qquad 
			\widetilde{r}_{i,j}^2 := \overrightarrow{\widetilde{B}_{i,j} \widetilde{A}_{i,j+1}}.
		\end{aligned}
	\end{equation}
	(see Figure \ref{fig:def-kagome} as an example with $k=2$).
	
	\begin{figure}[!htb]
		\begin{minipage}{.48\linewidth}
			\centering
			\subfloat[]{\label{fig:ref-rs}\includegraphics[width=0.85\linewidth]{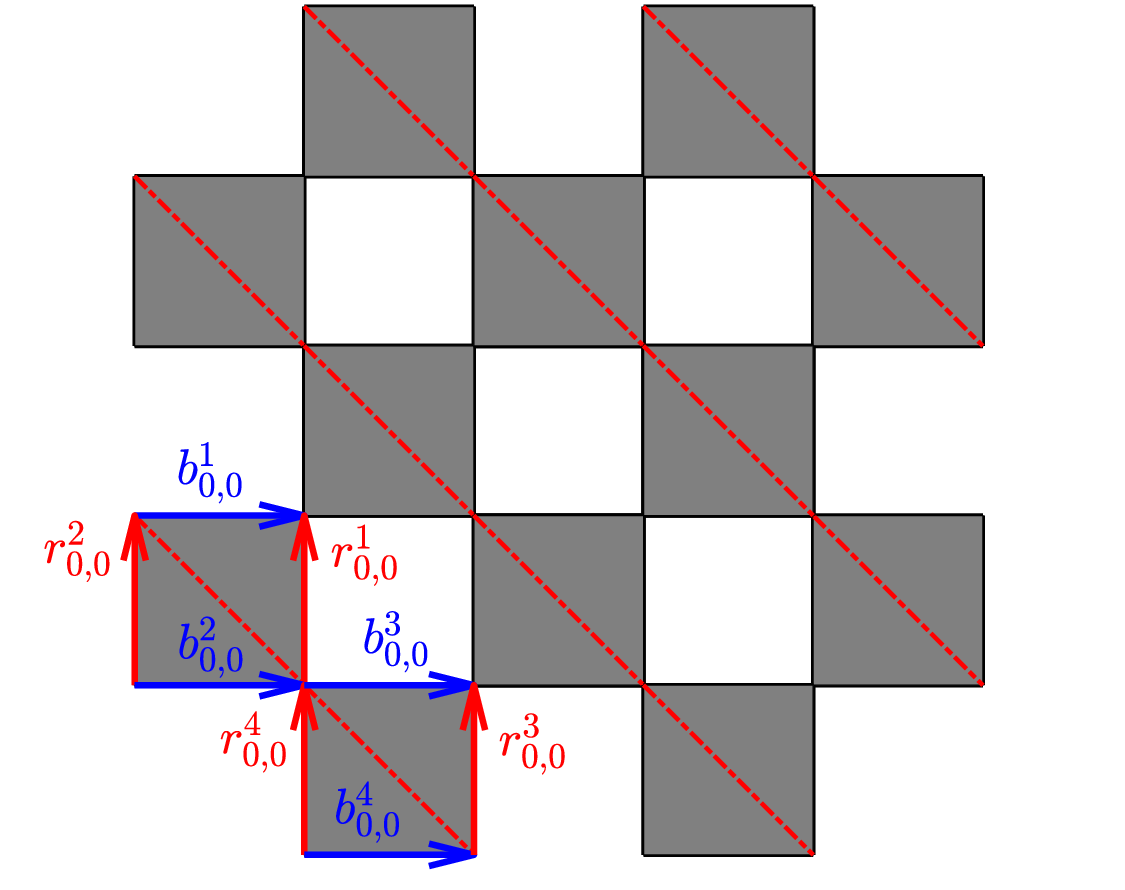}}
		\end{minipage}
		\begin{minipage}{.48\linewidth}
			\centering
			\subfloat[]{\label{fig:def-rs}\includegraphics[width=0.8\linewidth]{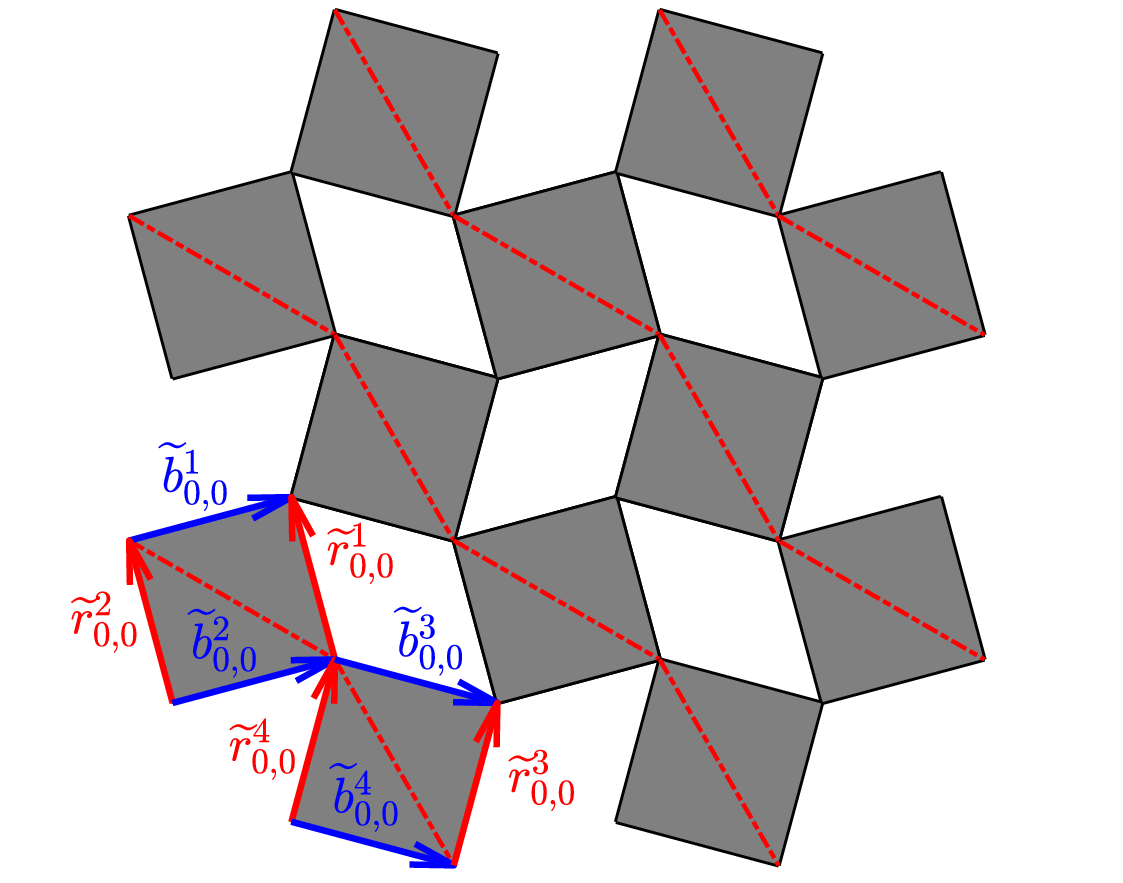}}
		\end{minipage}
		\caption{(a) The Rotating Squares lattice in the reference state; (b) the image of the Rotating Squares 
			lattice by the one-periodic mechanism. The marked vectors are $b_{0,0}^t, r_{0,0}^t$ in the reference state 
			and $\widetilde{b}_{0,0}^t, \widetilde{r}_{0,0}^t$ in the deformed state with $t=1,2,3,4$. Our basic set for 
			the springs of the Rotating Squares lattice has $10$ springs: $8$ corresponding to $b_{0,0}^t$ and $r_{0,0}^t$ 
			and $2$ diagonal springs along the diagonals of the squares bounded by $b_{0,0}^t$ and $r_{i,j}^t$.}
		\label{fig:rs-periodic}
	\end{figure}
	
	Turning now to the Rotating Squares lattice: we shall use the same conventions, but of course we must identify our basic
	sets of springs and triangles, and our definitions of $b_{i,j}^t$ and $r_{i,j}^t$. This is accomplished
	by Figure \ref{fig:rs-periodic}. (While the figure shows only $b_{0,0}^t$ and $r_{0,0}^t$, these fully determine
	$b_{i,j}^t$ and $r_{i,j}^t$ since the latter are the translates of the former by $iv_1 + jv_2$, where
	$v_1 = (2,0)^T$ and $v_2 = (0,2)^T$ are the periods of the Rotating Squares lattice.) There is one small difference
	between the Kagome and Rotating squares settings: for Kagome the basic set of triangles has just two elements (and therefore
	$b_{i,j}^t, r_{i,j}^t$ are defined for $t=1,2$) while for Rotating Squares the basic set of triangles has four elements
	(and therefore $b_{i,j}^t, r_{i,j}^t$ are defined for $t=1,2,3,4$). To permit a discussion that applies to both, we
	we let $\mathcal{T}$ be the set of indices of the triangles in our basic set; thus $\mathcal{T} = \{1,2\}$ for Kagome
	and $\mathcal{T} = \{1,2,3,4\}$ for Rotating Squares. With this convention, our notation is that for both the Kagome
	and Rotating Squares metamaterials
	\begin{itemize}
		
		\item in the reference triangle $T_{i,j,t}$, the horizontal vector is $b_{i,j}^t$
		and the vector in the other lattice direction is $r_{i,j}^t$;
		
		\item in the deformed triangle $\widetilde{T}_{i,j,t}$, the images of $b_{i,j}^t$
		and $r_{i,j}^t$ are $\widetilde{b}_{i,j}^t$ and $\widetilde{r}_{i,j}^t$.
	\end{itemize}
	
	We now use this framework to give a geometric characterization of the matrix $\lambda$, for deformations of the form
	$u(x) = \lambda x + \psi(x)$ where $\psi$ is $k$-periodic. To this end, let
	\begin{equation}\label{eqn:avg-br-ref}
		a_1 := \frac{1}{k^2} \sum_{i,j=0}^{k-1} \sum_{t \in \mathcal{T}} b_{i,j}^t, \qquad
		a_2 := \frac{1}{k^2} \sum_{i,j=0}^{k-1} \sum_{t \in \mathcal{T}} r_{i,j}^t .
	\end{equation}
	(We shall call these the ``averaged'' $b$ and $r$ vectors, though we have chosen not to normalize
	by $|\mathcal{T}|$.) Similarly, let
	\begin{equation}\label{eqn:avg-br-def}
		\widetilde{a}_1:= \frac{1}{k^2}\sum_{i,j=0}^{k-1} \sum_{t \in \mathcal{T}} \widetilde{b}_{i,j}^t, \qquad
		\widetilde{a}_2:= \frac{1}{k^2} \sum_{i,j=0}^{k-1} \sum_{t \in \mathcal{T}} \widetilde{r}_{i,j}^t.
	\end{equation}	
	be the averages of the corresponding sides of the deformed triangles.
	
	The values of $a_1$ and $a_2$ are easy to identify, since each $b_{i,j}^t$ is the vector $(1,0)^T$ and each
	$r_{i,j}^t$ is a unit-length vector in the other lattice direction; it follows that
	\begin{equation}\label{eqn:av-relation}
		\begin{aligned}
			&a_1 = |\mathcal{T}| \, (1,0)^T = v_1, \quad && a_2 = |\mathcal{T}| \, (1/2,\sqrt{3}/2)^T = v_2, \qquad && \text{for Kagome},\\
			&a_1 = |\mathcal{T}| \, (1,0)^T = 2v_1, \quad && a_2 = |\mathcal{T}| \, (0,1)^T = 2v_2, \qquad && \text{for Rotating Squares}.
		\end{aligned}
	\end{equation}
	To set the stage for our later arguments we observe that the relation
	\begin{equation}\label{eqn:a2-vs-a1-elementary}
		a_2 = R_\alpha a_1 \quad \mbox{where $\alpha = \pi/3$ for Kagome and $\alpha = \pi/2$ for Rotating Squares}
	\end{equation}
	(which follows from \eqref{eqn:av-relation}) also has another, rather different proof: it follows immediately
	(by summation) from the fact that $r_{i,j}^t = R_\alpha b_{i,j}^t$ for each $i,j,t$.
	
	The characterization of $\widetilde{a}_i$ is more interesting:
	
	\begin{lemma} \label{lemma:atilde-vs-lambda}
		For either the Kagome or Rotating Squares metamaterial, if $u(x) = \lambda x + \psi(x)$ with $k$-periodic $\psi$ then
		\begin{equation}\label{eqn:lambda-avec}
			\widetilde{a}_1 = \lambda a_1, \quad \widetilde{a}_2 = \lambda a_2.
		\end{equation}
	\end{lemma}
	\begin{proof}
		We begin by discussing the Kagome system. The key point is that
		\begin{align}
			(i) & \quad \mbox{for each $j=0,\ldots,k-1$, the segments $\{ b_{i,j}^t \}_{0\leq i \leq k-1, t\in \mathcal{T}}$ form
				a continuous path} \nonumber \\[-2pt]
			& \quad \mbox{whose two ends are related by $k$-periodicity; and similarly} \label{eqn:key-pt-kagome-i} \\[2pt]
			(ii) & \quad \mbox{for each $i=0,\ldots,k-1$, the segments $\{ r_{i,j}^t \}_{0\leq j \leq k-1, t\in \mathcal{T}}$ form
				a continuous path} \nonumber \nonumber\\[-2pt]
			& \quad \mbox{whose two ends are related by $k$-periodicity.} \label{eqn:key-pt-kagome-ii}
		\end{align}
		To explain why this implies \eqref{eqn:lambda-avec}, let us explain in detail why point (i) implies that
		$\widetilde{a}_1 = \lambda a_1$, focusing on the case $k=2$ (depicted by Figure \ref{fig:kagome-2-periodic}) for the
		purpose of visualization. Taking $j=0$ in (i)
		(that is, focusing on the lower row of blue edges in Figure \ref{fig:ref-kagome}, which maps to the lower blue 
		zigzag line in Figure \ref{fig:def-kagome}) and recalling \eqref{eqn:br-vec-ref}--\eqref{eqn:br-vec-def} (which amount, when we view
		$b_{i,j}^t$ and $\widetilde{b}_{i,j}^t$ as vectors, to $b_{i,j}^1 = B_{i,j} - C_{i,j}$, etc.)
		we see that
		$$
		\sum_{i=0}^{1} \sum_{t \in \mathcal{T}} b_{i,0}^t \quad \mbox{and} \quad
		\sum_{i=0}^{1} \sum_{t \in \mathcal{T}} \widetilde{b}_{i,0}^t
		$$
		are telescoping sums that reduce to the differences of the associated paths' endpoints
		$$
		C_{2,0} - C_{0,0} \quad \mbox{and} \quad u(C_{2,0}) - u(C_{0,0}) .
		$$
		Now, $C_{2,0}$ and $C_{0,0}$ are related by $2$-periodicity (indeed, $C_{2,0} =  C_{0,0} + 2v_1$). It follows that
		$\psi(C_{2,0}) = \psi(C_{0,0})$, so the preceding expressions reduce to
		$$
		2 v_1 \quad \mbox{and} \quad \lambda (2v_1).
		$$
		In short: we have shown that when $j=0$ and $k=2$,
		\begin{equation} \label{eqn:kagome-horiz-averages}
			\sum_{i=0}^{k-1} \sum_{t \in \mathcal{T}} \widetilde{b}_{i,j}^t =
			\lambda  \sum_{i=0}^{k-1} \sum_{t \in \mathcal{T}} b_{i,j}^t  .
		\end{equation}
		This relation is in fact true for any $k$ and for each $j=0 \ldots k-1 $, by exactly the
		same argument, and summing over $j$ gives the desired result $\widetilde{a}_1 = \lambda a_1$
		for the Kagome metamaterial. The proof that $\widetilde{a}_2 = \lambda a_2$ is entirely parallel.
		
		Turning now to the Rotating Squares metamaterial, we must proceed slightly differently since the analogues of
		\eqref{eqn:key-pt-kagome-i}--\eqref{eqn:key-pt-kagome-ii} are not true. The fix is easy, once we recognize that
		$\{b_{i,j}^t\}_{0\leq i,j \leq k-1, t \in \mathcal{T}}$ are a basic set for the lattice's horizontal springs
		under the equivalence relation associated with $k$-periodicity (and the analogous statement holds for
		$\{r_{i,j}^t\}$ with ``horizontal springs'' replaced by
		``vertical springs''). The averages
		\eqref{eqn:avg-br-ref}--\eqref{eqn:avg-br-def} defining $a_1$ and $\widetilde{a}_1$ can
		be calculated using \emph{any} basic set for horizontal segments under the equivalence relation
		associated with $k$-periodicity (and similarly, $a_2$ and $\widetilde{a}_2$ can be calculating using any
		basic set for the vertical springs). The argument we used to prove $\widetilde{a}_1 = \lambda a_1$
		for Kagome works perfectly for Rotating Squares, if as the basic set of horizontal springs we use a
		union of horizontal chains of $2k$ springs. It is easy to see that there is such a choice. (When $k=1$, for
		example, rather than the basic set of horizontal springs $\{ b^1_{0,0}, b^2_{0,0}, b^3_{0,0}, b^4_{0,0} \}$
		shown in Figure \ref{fig:ref-rs}, we would use $\{b^1_{0,-1}, b^2_{0,0}, b^3_{0,0}, b^4_{0,0} \}$.) Similarly, the
		argument we used to prove $\widetilde{a}_2 = \lambda a_2$ for Kagome works also for Rotating Squares, when we
		use a basic set of vertical springs that's a union of chains of $2k$ springs.
	\end{proof}
	
	To provide a geometrical understanding of the preceding result, let us connect its proof to 
	Figure \ref{fig:kagome-conformal-mechanism}. Equation \eqref{eqn:kagome-horiz-averages} reduces when $k=2$ and $j=0$ to the
	statement that $v_1^{\text{def}} = \lambda v_1$ for Kagome. Similarly, the analogous statement for $r_{i,j}^t$
	reduces when $k=2$ and $j=0$ to $v_2^\text{def} = \lambda v_2$. Moreover, \eqref{eqn:kagome-horiz-averages} holds
	for every $j$, not just $j=0$ -- so it does not matter which horizontal line of springs is used to
	define $v_1^{\text{def}}$ and which $60$-degree line of springs is used to define $v_2^{\text{def}}$. These observations
	are, of course, related to the fact that when $u = \lambda x + \psi(x)$ with $\psi$ $k$-periodic, the macroscopic
	deformation gradient is $\lambda$.
	
	\subsection{The proof of Proposition \ref{lemma-geometric}} \label{subsec:proof-of-geometric-lemma}
	
	The following proposition is, roughly speaking, the analogue of our elementary observation \eqref{eqn:a2-vs-a1-elementary}
	for the deformed lattice associated with a $k$-periodic mechanism.
	
	\begin{proposition}\label{prop:periodic-mechanism-macroscopic}
		In the Kagome and Rotating Squares settings, if $u(x) = \lambda x + \psi(x)$ be a $k$-periodic mechanism 
		(in the sense of Definition \ref{def:periodic-mechanism}) then the matrix $\lambda$ must satisfy
		\begin{equation}\label{eqn:lemma-geometric}
			R_\alpha \lambda e_1 = \lambda R_\alpha e_1, \qquad R_\alpha =
			\begin{pmatrix}
				\cos \alpha & -\sin \alpha\\
				\sin \alpha & \cos \alpha
			\end{pmatrix}, \qquad
			e_1 =
			\begin{pmatrix}
				1\\0
			\end{pmatrix},
		\end{equation}
		where $\alpha=\frac{\pi}{3}$ for the Kagome metamaterial and $\alpha=\frac{\pi}{2}$ for the Rotating Squares metamaterial.
	\end{proposition}
	
	\begin{proof}
		For a $k$-periodic mechanism $u(x)$, all triangles $T_{i,j,t}$ retain the same shape after deformation and none 
		of the triangles is flipped. Therefore, the deformed vectors $\widetilde{b}_{i,j}^t$ and $\widetilde{r}_{i,j}^t$ satisfy
		\begin{equation}\label{eqn:def-relationship}
			\widetilde{r}_{i,j}^t = R_{\alpha} \widetilde{b}_{i,j}^t, \qquad i,j=0,\dots,k-1, \quad t \in \mathcal{T}.
		\end{equation}
		It follows that the averages $\widetilde{a}_1,\widetilde{a}_2$ (defined by \eqref{eqn:avg-br-def}) satisfy
		\begin{equation*}
			\widetilde{a}_2 = R_\alpha \widetilde{a}_1.
		\end{equation*}
		Since $\widetilde{a}_1 = \lambda a_1$ by Lemma \ref{lemma:atilde-vs-lambda}, we conclude that
		\begin{equation*} 
			\widetilde{a}_2 = R_\alpha \lambda a_1 .
		\end{equation*}
		But we also have, using Lemma \ref{lemma:atilde-vs-lambda} first then \eqref{eqn:avg-br-ref}, that
		\begin{equation*} 
			\widetilde{a}_2 = \lambda a_2 = \lambda R_{\alpha} a_1 .
		\end{equation*}
		Combining these relations, we conclude that
		$$
		(\lambda R_\alpha - R_\alpha \lambda) a_1 = 0 .
		$$
		Since $a_1 = |\mathcal{T}|e_1$ (see \eqref{eqn:av-relation}), this is equivalent to the desired assertion
		\eqref{eqn:lemma-geometric}.
	\end{proof}
	
	It is obvious that if $\lambda = c R$ for some $R \in SO(2)$ then $(\lambda R_\alpha - R_\alpha \lambda) e_1 = 0$. We
	need to show that there are no \emph{other} choices of $\lambda$ with this property. At the heart of our proof lies
	the following result:
	
	\begin{proposition}\label{prop:our-geometric-claim}
		For any matrix $\lambda \in \mathbb{R}^{2\times 2}$, we have
		\begin{equation}\label{eqn:direct-computation}
			\Big|\big(\lambda R_{\alpha}-R_{\alpha} \lambda\big) e_1\Big| =
			\begin{cases}
				|\sin\alpha| |\lambda_2 - \lambda_1|, & \det\lambda\geq 0,\\
				|\sin\alpha| (\lambda_2 + \lambda_1), & \det\lambda < 0.
			\end{cases}
		\end{equation}
		where $\lambda_1, \lambda_2$ are the principal stretches of $\lambda$ (the eigenvalues of $(\lambda^T \lambda)^{1/2}$)
		and $e_1 = (1,0)^T$.
	\end{proposition}
	\begin{proof}
		The proof of \eqref{eqn:direct-computation} comes from a standard calculation using a singular value decomposition
		of $\lambda$. To deal with the sign of the determinant, we use a particular choice of singular value decomposition
		such that
		\begin{equation}\label{eqn:lambda-svd}
			\lambda =
			\begin{cases}
				USV^T, \qquad \text{with }U,V \in SO(2), & \det \lambda \geq 0,\\
				USV^T, \qquad \text{with }U \in O(2)\setminus SO(2) \text{ and } V \in SO(2), & \det\lambda < 0.
			\end{cases}
		\end{equation}
		where $ S = \text{diag}(\lambda_1, \lambda_2)$.
		Such a singular value decomposition always exists: when $\det \lambda \geq 0$, the two orthonormal matrices 
		$U,V$ must satisfy $\det U \det V > 0$. If their determinants are both negative, then we can choose 
		$\widetilde{U} = U F$ and $\widetilde{V} = VF$, where $F = \text{diag}(1,-1) \in O(2)$. The matrices 
		$\widetilde{U}, S, \widetilde{V}$ form a singular decomposition of $\lambda$, i.e. 
		$\lambda = \widetilde{U} S \widetilde{V}^T$ since $F S F^T = S$. Similarly, when 
		$\det \lambda < 0$, we have $\det U \det V < 0$. If $\det U > 0$ and $\det V <0$, then we obtain 
		\eqref{eqn:lambda-svd} by choosing the new matrices $\widetilde{U}, \widetilde{V}$ as above.
		
		Now we prove \eqref{eqn:direct-computation} by a direct calculation. When $\det \lambda \geq 0$ we use
		\eqref{eqn:lambda-svd} to write the left side of \eqref{eqn:direct-computation} as
		\begin{equation}
			\Big|\big(\lambda R_{\alpha}-R_{\alpha} \lambda\big) e_1\Big| = 
			\Big|U\big(S R_{\alpha}-R_{\alpha} S\big) V^T e_1\Big|= 
			\Big|\big(S R_{\alpha}-R_{\alpha} S\big) V^T e_1\Big|, \label{eqn:direct-computation-2}
		\end{equation}
		where $S = \text{diag}(\lambda_1, \lambda_2)$ and $U,V \in SO(2)$. Since $V^T e_1$ is a unit vector, we write it as $V^T e_1 = (\cos\theta, \sin\theta)^T$ for some angle $\theta$. Plugging this into \eqref{eqn:direct-computation-2}, we obtain
		\begin{equation*}
			\begin{aligned}
				&\big(S R_{\alpha}-R_{\alpha} S\big) V^T e_1 = 
				\begin{pmatrix}
					\lambda_1 \cos(\theta + \alpha)\\
					\lambda_2 \sin(\theta + \alpha)
				\end{pmatrix}
				-
				\begin{pmatrix}
					\cos\alpha & -\sin\alpha\\
					\sin\alpha & \cos\alpha
				\end{pmatrix}\begin{pmatrix}
					\lambda_1 \cos\theta\\
					\lambda_2 \sin\theta
				\end{pmatrix} \\
				=& \begin{pmatrix}
					\lambda_1 \cos\theta \cos\alpha - \lambda_1 \sin\theta \sin\alpha - \lambda_1 \cos\theta \cos\alpha + \lambda_2 \sin\theta \sin\alpha\\
					\lambda_2 \sin\theta \cos\alpha + \lambda_2 \cos\theta \sin\alpha - \lambda_1 \cos\theta \sin\alpha - \lambda_2 \sin \theta \cos\alpha
				\end{pmatrix}
				= 
				\begin{pmatrix}
					(\lambda_2 - \lambda_1) \sin\theta \sin\alpha\\
					(\lambda_2 - \lambda_1) \cos\theta \sin\alpha
				\end{pmatrix} .
			\end{aligned}
		\end{equation*}
		This completes the argument for \eqref{eqn:direct-computation} when $\det \lambda \geq 0$, since
		\begin{equation}\label{eqn:positive-lambda-det}
			\Big|\big(\lambda R_{\alpha}-R_{\alpha} \lambda\big) e_1\Big| = 
			\Big|\big(S R_{\alpha}-R_{\alpha} S\big) V^T e_1\Big| = |\lambda_1 - \lambda_2| |\sin \alpha|.
		\end{equation}
		
		When $\det \lambda < 0$, we use the singular value decomposition $\lambda = USV^T$ with 
		$U \in O(2)\setminus SO(2)$ and $V \in SO(2)$. By choosing $\widetilde{U} = UF$ with 
		$F = \text{diag}(1,-1)$, we obtain that $\widetilde{U} \in SO(2)$ and $\lambda = \widetilde{U} F S V^T$. 
		Then the left hand side of \eqref{eqn:direct-computation} becomes
		\begin{equation*}
			\begin{aligned}
				\Big|\big(\lambda R_{\alpha}-R_{\alpha} \lambda\big) e_1\Big| &= 
				\Big|\widetilde{U}\big(FS R_{\alpha}-R_{\alpha} FS\big) V^T e_1\Big| = 
				\Big|\big(FS R_{\alpha}-R_{\alpha} FS\big) V^T e_1\Big|.
			\end{aligned}
		\end{equation*}
		Since $FS = \text{diag}(\lambda_1, -\lambda_2)$, we can substitute $S$ with $FS$ in 
		\eqref{eqn:positive-lambda-det} and obtain that
		\begin{equation*}
			\Big|\big(S R_{\alpha}-R_{\alpha} S\big) V^T e_1\Big|= |\sin \alpha|(\lambda_1 + \lambda_2).
		\end{equation*}
		This completes the proof of \eqref{eqn:direct-computation}.
	\end{proof}
	
	\begin{remark}\label{rmk:diff-direction}
		In Proposition \ref{prop:periodic-mechanism-macroscopic} and Proposition \ref{prop:our-geometric-claim}, we can choose $e_1$ to be any unit vector $e$ with $|e|=1$. In fact, substituting $e$ into \eqref{eqn:direct-computation-2} yields $V^T e$, which is still a unit vector and can be written as $(\cos\theta, \sin\theta)^T$. The proof then proceeds unchanged.
	\end{remark}
	
	Finally, we combine Proposition \ref{prop:periodic-mechanism-macroscopic} and Proposition \ref{prop:our-geometric-claim} to 
	prove Proposition \ref{lemma-geometric}.
	
	\begin{proof}[Proof of Proposition \ref{lemma-geometric}]
		Our goal is to show, for both Kagome and Rotating Squares, that if $u(x) = \lambda x + \psi$ is a $k$-periodic
		mechanism then $\lambda$ is an isotropic compression, i.e. $\lambda = cR$ with $0 \leq c \leq 1$ and $R \in SO(2)$.
		
		We start by showing that $\lambda$ is isotropic. From Proposition \ref{prop:periodic-mechanism-macroscopic}, we know
		that $(\lambda R_{\alpha} - R_{\alpha} \lambda) e_1 = 0$ (with $\alpha = \frac{\pi}{3}$ for Kagome,
		$\alpha = \frac{\pi}{2}$ for Rotating Squares). It follows from Proposition \ref{prop:our-geometric-claim} that the principal
		stretches of $\lambda$ must be equal and $\det \lambda \geq 0$. Using a singular value decomposition of $\lambda$, we 
		obtain $\lambda = c UV^T$ with $U,V\in O(2)$. The constant $c$ is automatically non-negative, since it equals the repeated
		singular value: $c=\lambda_1 = \lambda_2 \geq 0$. Since $\det\lambda = c^2 \det(UV^T)\geq 0$, we obtain that
		$UV^T$ is orientation-preserving, i.e. $UV^T \in SO(2)$. This proves our assertion that
		$\lambda=cR$ for some $c \geq 0$ and $R \in SO(2)$.
		
		Our remaining task is to show that $ c \leq 1$. We know from Lemma \ref{lemma:atilde-vs-lambda} that $\lambda a_1 = \widetilde{a}_1$. 
		Since $\lambda = cR$, we have $|\widetilde{a}_1| = c |a_1|$. Now, for both the Kagome and Rotating Squares metamaterials, 
		the reference state (see the left sides of Figure \ref{fig:kagome-2-periodic} and Figure \ref{fig:rs-periodic}) has straight lines of 
		springs in the lattice direction $v_1$, which is a multiple of $a_1$ (see \eqref{eqn:av-relation}). Since a periodic 
		mechanism deforms those straight lines into zigzag lines (each segment being of length one), the
		lines must experience macroscopic compression, i.e. $|\widetilde{a}_1| \leq |a_1|$. It follows that $c \leq 1$, as
		claimed.
	\end{proof}
	
	\section{Lower bounds for the effective energy density $\overline{W}^\eta(\lambda)$}\label{sec:lower-bd}
	
	In this section we prove Theorem \ref{thm:lower-bound}, which gives our lower bound for the effective energy density 
	$\overline{W}^\eta(\lambda)$. As already noted in Section \ref{sec:geo-argument}, our task is basically to prove a lower bound
	on the averaged energy $\overline{E^\eta}(\lambda,\psi,kU)$ (defined by \eqref{eqn:avg-energy}) that's uniform in $k$
	and $\psi$.
	
	Our lower bound \eqref{eqn:eff-lower-bd} has two parts, the first being $(\lambda_1 \pm \lambda_2)^2$ (depending upon the sign of $\det \lambda$) and the second being $(\lambda_1 - 1)_+^2 + (\lambda_2 - 1)_+^2$. Therefore it is natural to separate our proof into two parts by showing the 
	following two propositions:
	
	\begin{proposition}\label{prop:isotropic-bd}
		For the Kagome and Rotating Squares metamaterials, there exist positive constants $\eta_0$ and $C_1$
		such that for any $k \in \mathbb{N}$, the averaged energy $\overline{E^\eta}(\lambda,\psi,kU)$ is lower bounded by
		\begin{equation}\label{eqn:istropic-bd}
			\overline{E^\eta}(\lambda,\psi,kU) \geq \begin{cases}
				C_1 (\lambda_1 - \lambda_2)^2, & \det(\lambda) \geq 0, \\
				C_1 (\lambda_1 + \lambda_2)^2, & \det(\lambda) < 0,
			\end{cases}
		\end{equation}
		for any $\lambda$, any $k$-periodic $\psi$, and any positive $\eta < \eta_0$.
	\end{proposition}
	
	\begin{proposition}\label{prop:compressive-bd}
		For the Kagome and Rotating Squares metamaterials, there exist positive constants $\eta_0 $ and $C_2 $ such that
		for any $k \in \mathbb{N}$, the averaged energy $\overline{E^\eta}(\lambda,\psi,kU)$ is lower bounded by
		\begin{equation}\label{eqn:compressive-bd}
			\overline{E^\eta}(\lambda,\psi,kU) \geq C_2 \Big((\lambda_1 - 1)_+^2 + (\lambda_2 - 1)_+^2\Big)
		\end{equation}
		for any $\lambda$, any $k$-periodic $\psi$, and any positive $\eta < \eta_0$.
	\end{proposition}
	
	We start, in Section \ref{subsec:averaged-energy-kU}, by representing the averaged energy in terms of the
	basic sets $b_{i,j}^t, r_{i,j}^t$ introduced in Section \ref{subsec:k-periodic-def}. Then we prove 
	Proposition \ref{prop:isotropic-bd} in Section \ref{subsec:isotropic-bd} and Proposition \ref{prop:compressive-bd} in Section \ref{subsec:compress-bd}.
	
	\subsection{The averaged energy on $kU$} \label{subsec:averaged-energy-kU}
	
	Throughout this section we are considering only deformations of the form $u(x) = \lambda x + \psi(x)$ where $\psi$ is
	$k$-periodic. We recall the definition of $\overline{E^\eta}(\lambda, \psi, kU)$ for the reader's convenience:
	\begin{equation} \label{eqn:avg-energy-repeated}
		\overline{E^\eta}(\lambda, \psi, kU) := \frac{1}{k^2 |U|}
		\sum_{\alpha_1, \alpha_2=0}^{k-1} E^\eta(\lambda x + \psi,U+\alpha_1 v_1 + \alpha_2 v_2) .
	\end{equation}
	We also recall that for both Kagome and Rotating Squares, the energy $E^\eta (u,U)$ of our unit cell has
	two parts: a spring energy (obtained by adding the energies of the springs associated with the unit cell) and
	a penalization energy (obtained by adding $|T| f^\eta ( \det ( D u |_T ) ) $ for certain triangles $T$).
	
	We now observe that when $u(x) = \lambda x + \psi(x)$ with $\psi(x)$ $k$-periodic,
	\begin{itemize}
		\item the energies of the individual springs are $k$-periodic,
		\item the penalization energies of the triangles are also $k$-periodic, and
		\item for both the springs and the penalized triangles, the sum on the right side of
		\eqref{eqn:avg-energy-repeated} samples each of the equivalence classes associated with $k$-periodicity exactly once.
	\end{itemize}
	Therefore the right side of \eqref{eqn:avg-energy-repeated} is really a sum over equivalence classes, and to evaluate it
	we can use any basic sets of springs and penalized triangles that we wish. (Here, as in Section \ref{subsec:k-periodic-def},
	we define a basic set to be a list that samples each equivalence class exactly once.)
	
	For Kagome, we have already introduced a basic set for the equilateral triangles, namely $T_{i,j,t}$ for $0 \leq i,j \leq k-1$
	and $t \in \mathcal{T}$. We have also introduced basic sets of horizontal and $60$-degree springs, namely the edges labeled
	$b_{i,j}^t$ and $r_{i,j}^t$ for $0 \leq i,j \leq k-1$ and $t \in \mathcal{T}$. A basic set for the springs in the
	remaining ($120$-degree) direction is obtained by taking, for each $i,j$ and $t$, the third side of the triangle
	of which $b_{i,j}^t$ and $r_{i,j}^t$ are edges. Using these choices, we get
	\begin{equation}\label{eqn:avg-energy-rewritten}
		\overline{E^\eta}(\lambda,\psi,kU) = \frac{1}{k^2 |U|}\sum_{i,j=0}^{k-1} \sum_{t \in \mathcal{T}} E^\eta_{i,j,t} (u)
	\end{equation}
	where $E^\eta_{i,j,t} (u)$ is the sum of the spring and penalty energies associated with triangle 
	$T_{i,j,t}$ for $i,j=0,\dots,k-1$ and $t \in \mathcal{T}$, i.e.
	\begin{equation}\label{eqn:energy-ijs}
		\begin{aligned}
			E^\eta_{i,j,t} (u) &= E^{\text{spr}}_{i,j,t} (u) + E^{\eta,\text{pen}}_{i,j,t} (u)
		\end{aligned}
	\end{equation}
	with
	\begin{equation}\label{eqn:energy-ijs-kagome}
		\begin{aligned}
			E^{\text{spr}}_{i,j,t} (u) &= 
			\bigg(|\widetilde{b}_{i,j}^t| - 1\bigg)^2 + \bigg(|\widetilde{r}_{i,j}^t| - 1\bigg)^2 + 
			\bigg(|\widetilde{b}_{i,j}^t - \widetilde{r}_{i,j}^1| - 1\bigg)^2,\\
			E^{\eta,\text{pen}}_{i,j,t} (u) &= \frac{\sqrt{3}}{4}f^\eta(\det D u|_{T_{i,j,t}}) .
		\end{aligned}
	\end{equation}
	(Here $\frac{\sqrt{3}}{4}$ is the area of the triangle $T_{i,j,t}$ in the Kagome setting.)
	
	The situation for Rotating Squares is only slightly different. There, too, the right side of
	\eqref{eqn:avg-energy-repeated} is really a sum over equivalence classes associated with $k$-periodicity. Moreover
	there, too, we have already identified basic sets for the triangles that get penalization terms
	($T_{i,j,t}$ for $0 \leq i,j \leq k-1$ and $t \in \mathcal{T}$) and for the horizontal and vertical springs
	($b_{i,j}^t$ and $r_{i,j}^t$ for $0 \leq i,j \leq k-1$ and $t \in \mathcal{T}$). Howeever, each diagonal spring in
	the Rotating Squares metamaterial belongs to \emph{two} distinct triangles. Therefore, choosing for each $i,j,t$ the
	third side of the triangle associated with $b_{i,j}^t$ and $r_{i,j}^t$ does not give a basic set for the diagonal springs,
	since it counts each diagonal spring twice. \emph{It was for this reason that we used $2$ for the spring 
		constants of the diagonal springs when we defined $E_{\text{spr}}(u,U)$ for the Rotating Squares example in 
		\eqref{eqn:rs-spring-energy}.} Taking this into account and arguing as
	we did for Kagome, one easily sees that \eqref{eqn:avg-energy-rewritten} holds for the Rotating Squares metamaterial
	with the convention \eqref{eqn:energy-ijs}, when
	\begin{equation}\label{eqn:energy-ijs-rs}
		\begin{aligned}
			E^{\text{spr}}_{i,j,t} (u) &= \bigg(|\widetilde{b}_{i,j}^t| - 1\bigg)^2 + \bigg(|\widetilde{r}_{i,j}^t| - 1\bigg)^2 + \bigg(|\widetilde{b}_{i,j}^t - \widetilde{r}_{i,j}^1| - \sqrt{2}\bigg)^2,\\
			E^{\eta,\text{pen}}_{i,j,t} (u) &= \frac{1}{2}f^\eta(\det \nabla u|_{T_{i,j,t}}) .
		\end{aligned}
	\end{equation}
	(Here $\frac{1}{2}$ is the area of the triangle $T_{i,j,t}$ in the Rotating Squares setting.)
	
	\subsection{The isotropic bound}\label{subsec:isotropic-bd}
	Our goal in this subsection is to prove Proposition \ref{prop:isotropic-bd}. 
	The underlying idea is relatively simple. We showed in 
	Section \ref{sec:geo-argument} that when $u(x) = \lambda x + \psi(x)$ 
	with $\psi$ $k$-periodic, if (a) the length of each spring
	is preserved, and (b) none of the penalized triangles has 
	its orientation reversed, then the macroscopic deformation 
	gradient $\lambda$ must be isotropic. This was done by 
	combining Proposition \ref{prop:periodic-mechanism-macroscopic} and 
	Proposition \ref{prop:our-geometric-claim}. To prove 
	Proposition \ref{prop:isotropic-bd} we must show that our arguments 
	were robust, in the sense that if (a) the spring energy is 
	small, and (b) the penalization energy is small, then the 
	failure of isotropy is also small. 
	
	We start by introducing some notation, to distinguish 
	the triangles whose orientations are preserved vs. reversed: 
	\begin{equation}\label{eqn:op-set}
		\begin{aligned}
			OP_t &= \{(i,j) \;|\; i,j = 0, \dots, k-1 \text{ and } \det(D u|_{T_{i,j,t}}) \geq 0\}, \qquad t \in \mathcal{T},\\
			OR_t &= \{(i,j) \;|\; i,j = 0, \dots, k-1 \text{ and } \det(D u|_{T_{i,j,t}}) < 0\}, \qquad t \in \mathcal{T}.
		\end{aligned}
	\end{equation}
	(Thus a triangle of type $t$ is in $OP_t$ if its orientation is preserved, and in $OR_t$  if its orientation is reversed.)
	
\begin{figure}[!thb]
		\centering
		\includegraphics[width=0.6\linewidth]{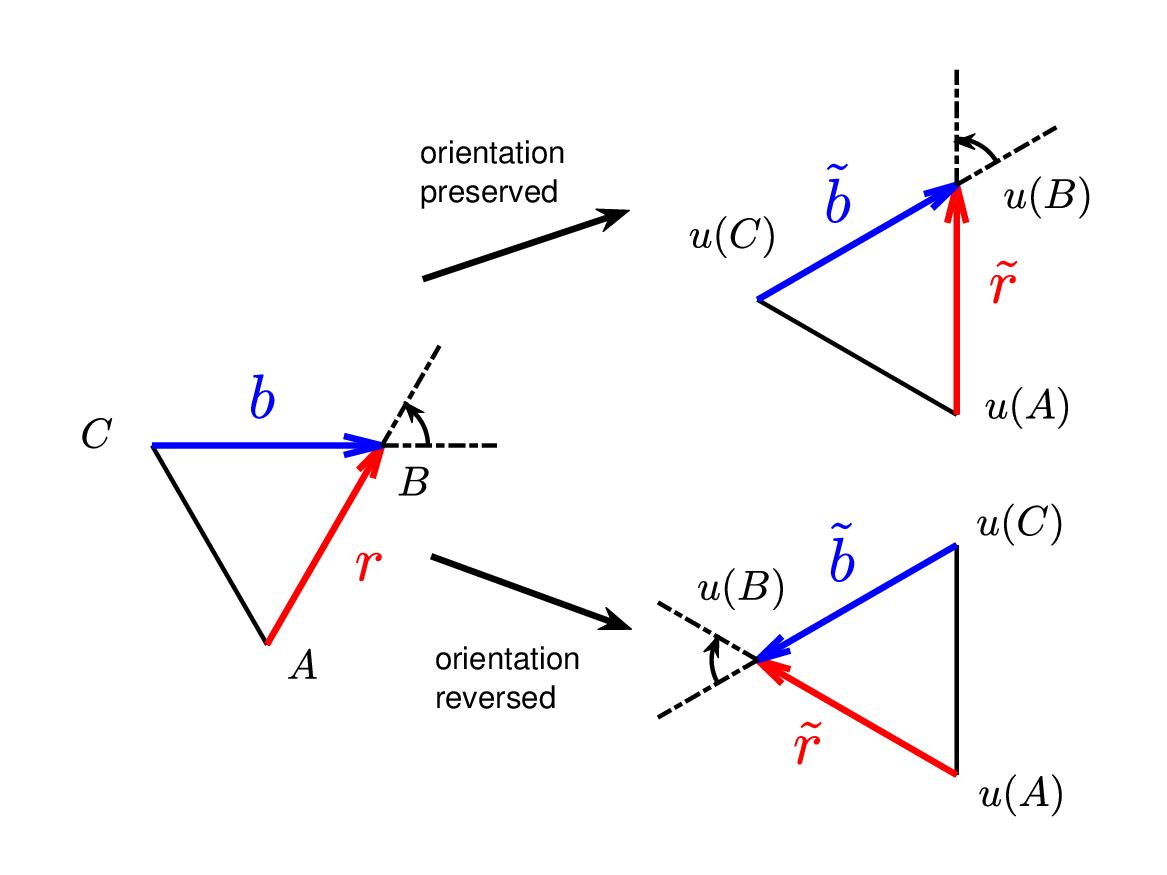}
		\caption{The difference between orientation preserving and orientation reversing deformations: for an equilateral triangle, if a deformation almost preserves the lengths of the sides, then the angle from $\widetilde{b}$ to $\widetilde{r}$ is close to $\pi/3$ if the orientation is preserved, and it is close to $-\pi/3$ if the orientation is reversed.}
		\label{fig:orientation}
	\end{figure}

	Next we introduce the rigidity result for triangles that underlies our argument. Informally it says that for each
	triangle $T_{i,j,t}$, if its spring energy is small, then
	the deformed triangle $\widetilde{T}_{i,j,t}$ closely resembles $T_{i,j,t}$ in shape, either with preserved or reversed orientation. 
	We illustrate this in Figure \ref{fig:orientation}, where the reference triangle is equilateral with side length $1$. If the deformation $u$ has small elastic energy, then the images of its three edge 
	vectors ($\widetilde{b}, \widetilde{r}$, and 
	$\widetilde{r} - \widetilde{b}$) remain close to unit length. Moreover, the angle from $\widetilde{b}$ to $\widetilde{r}$ is approximately $\pi/3$ when the orientation is preserved, or approximately $-\pi/3$ when the orientation is reversed. Putting this assertion slightly differently: the mismatch between $\widetilde{r}$ and $R_{\frac{\pi}{3}} \widetilde{b}$ (or $R_{-\frac{\pi}{3}} \widetilde{b}$ in the reversed case) should be small when the orientation is preserved (or reversed). The key point is that \emph{the mismatch is bounded above by the spring energy}. We state this result more formally (for any isosceles 
	triangle) as follows.
	\begin{proposition}\label{prop:iso-triangle}
		Consider an isosceles triangle $\Delta ABC$ with $|AB| = |BC|=1$ and angle $\alpha$ between edges $AB$ and $CB$. For any deformation $u$ (defined at the nodes $A,B,C$), let
		\begin{equation}\label{eqn:spr-energy-iso}
			E^\text{spr}(u) := \Big(|u(A)-u(B)|-1\Big)^2 + \Big(|u(C)-u(B)|-1\Big)^2 + \Big(|u(A)-u(C)|-|A-C|\Big)^2 .
		\end{equation}
		Then there is a constant $c$, depending only on $\alpha$, such that
		\begin{equation}\label{eqn:lower-bd-spr-energy}
			E^\text{spr}(u) \geq \begin{cases}
				c |\widetilde{r} - R_\alpha \widetilde{b}|^2, & \det(\nabla u) \geq 0,\\
				c |\widetilde{r} - R_{-\alpha} \widetilde{b}|^2, & \det(\nabla u) < 0,\\
			\end{cases}
		\end{equation}
		(Here the vectors $\widetilde{b},\widetilde{r}$ are the images of $ b=\overrightarrow{CB}$ and $r=\overrightarrow{AB}$.)
	\end{proposition}
	When we apply this proposition we shall of course take 
	$\alpha = \pi/3$ for the Kagome metamaterial and $\alpha = \pi/2$ 
	for the Rotating Squares metamaterial. The proof of 
	Proposition \ref{prop:iso-triangle} is relatively simple: it uses the 
	Law of Cosines to get the cosine of the angle between $\widetilde{b}$ and $\widetilde{r}$ from the lengths of the sides, and the monotonicity of the cosine function to know that if the cosine doesn't change much then neither does the angle. The details, however, are somewhat tedious, so we have chosen to put the proof of 
	Proposition \ref{prop:iso-triangle} in Appendix \ref{app:proof-iso-triangle}.
	
	Now we use Proposition \ref{prop:iso-triangle} to prove Proposition \ref{prop:isotropic-bd}.
	\begin{proof}[Proof of Proposition \ref{prop:isotropic-bd}]
		The idea is that when the averaged energy is small and the penalty constant $\eta$ is small (the penalty magnitude $1/\eta$ is large), most triangles have their orientations preserved;
		Proposition \ref{prop:iso-triangle} assures us, for such triangles, 
		that the mismatch between $\widetilde{r}_{i,j}^t$ and $R_\alpha \widetilde{b}_{i,j}^t$ is small. We'll show that if $\eta$ is small then the contribution of the orientation-reversed triangles can be treated as an error term. This will lead to the conclusion that if the averaged energy is small then averaged deformed lattice vector $\widetilde{a}_2$ in \eqref{eqn:avg-br-def} is close to $R_\alpha \widetilde{a}_1$. The desired lower bound will then follow using arguments parallel to the ones used in Section \ref{subsec:proof-of-geometric-lemma}. 
		
		Starting the proof now, we apply Proposition \ref{prop:iso-triangle} to get 
		that for every triangle $T_{i,j,t}$, the associated spring energy satisfies
		\begin{equation}\label{eqn:spring-lower-bound-triangle}
			E_{i,j,t}^\text{spr}(u) \geq c|\widetilde{z}_{i,j}^t |^2,
		\end{equation}
		where the mismatch vectors $\widetilde{z}_{i,j}^t$ are defined as
		\begin{equation}\label{eqn:deviation-z}
			\widetilde{z}_{i,j}^t := \begin{cases}
				\widetilde{r}^t_{i,j} - R_{\alpha} \widetilde{b}^t_{i,j}, & (i,j)\in OP_t,\\
				\widetilde{r}^t_{i,j} - R_{-\alpha} \widetilde{b}^t_{i,j}, & (i,j)\in OR_t .
			\end{cases}
		\end{equation}

		We recall that when $u(x) = \lambda x + \psi(x)$ with $\psi$ $k$-periodic, the averaged deformed vectors $\widetilde{b}_{i,j}^t$ and $\widetilde{r}_{i,j}^t$ are equal to $\lambda a_1, \lambda a_2$ respectively (see \eqref{eqn:avg-br-def} and \eqref{eqn:lambda-avec}). To distinguish between the orientation-preserved and orientation-reversed triangles, we write this as
		\begin{equation*}
			\begin{aligned}
				\lambda a_1 &= \widetilde{a}_1 =\frac{1}{k^2} \sum_{t\in \mathcal{T}} \sum_{i,j=0}^{k-1} \widetilde{b}^t_{i,j} = \frac{1}{k^2} \Big[\sum_{t\in \mathcal{T}} \sum_{(i,j) \in OP_t} \widetilde{b}^t_{i,j} + \sum_{t\in \mathcal{T}} \sum_{(i,j) \in OR_t} \widetilde{b}^t_{i,j}\Big],\\
				\lambda a_2 &= \widetilde{a}_2 = \frac{1}{k^2} \sum_{t\in \mathcal{T}} \sum_{i,j=0}^{k-1} \widetilde{r}^t_{i,j} = \frac{1}{k^2} \Big[\sum_{t \in \mathcal{T}} \sum_{(i,j) \in OP_t} \widetilde{r}^t_{i,j} + \sum_{t \in \mathcal{T}} \sum_{(i,j) \in OR_t} \widetilde{r}^t_{i,j}\Big].
			\end{aligned}
		\end{equation*}
		We rewrite the latter using the mismatch terms $\widetilde{z}_{i,j}^t$ introduced in \eqref{eqn:deviation-z}:
		\begin{equation*}
			\lambda a_2 = \frac{1}{k^2} \Big[\sum_{t \in \mathcal{T}} \sum_{(i,j) \in OP_t} \big(R_{\alpha}\widetilde{b}^t_{i,j} + \widetilde{z}_{i,j}^t\big) + \sum_{t\in \mathcal{T}} \sum_{(i,j) \in OR_t} \big(R_{-\alpha}\widetilde{b}^t_{i,j} + \widetilde{z}_{i,j}^t\big)\Big].
		\end{equation*}
		Subtracting $R_{\alpha} \lambda a_1$ from both sides gives
		\begin{equation}\label{eqn:difference-flipping}
			\lambda a_2 - R_{\alpha} \lambda a_1= \frac{1}{k^2} \sum_{t \in \mathcal{T}} \sum_{i,j=0}^{k-1} \widetilde{z}_{i,j}^t + \frac{1}{k^2} \sum_{t \in \mathcal{T}} \sum_{(i,j) \in OR_t} \big(R_{-\alpha} - R_{\alpha}\big) \widetilde{b}_{i,j}^t. 
		\end{equation}
		By \eqref{eqn:direct-computation}, 
		$| \lambda a_2 - R_{\alpha} \lambda a_1|^2$ is a constant times the right hand side of the desired inequality \eqref{eqn:istropic-bd}. Therefore the proof will be done if we can show that the right hand side of \eqref{eqn:difference-flipping} is controlled by the averaged energy in the sense that
		\begin{equation}\label{eqn:average-iso-lower-bound}
			\Big|\frac{1}{k^2} \sum_{t \in \mathcal{T}} \sum_{i,j=0}^{k-1} \widetilde{z}_{i,j}^t + \frac{1}{k^2} \sum_{t \in \mathcal{T}} \sum_{(i,j) \in OR_t} \big(R_{-\alpha} - R_{\alpha}\big) \widetilde{b}_{i,j}^t\Big| \leq {C} \sqrt{\overline{E^\eta}}(\lambda,\psi,kU)
		\end{equation}
		with a constant $C$ that is independent of $\lambda$, $k$, $\psi$, and $\eta$.
		
		To show \eqref{eqn:average-iso-lower-bound}, we shall bound the averaged mismatch term by the averaged spring energy and the other term by the averaged penalty+spring energy. The bound on the mismatch is straightforward: we know from \eqref{eqn:spring-lower-bound-triangle} that $|\widetilde{z}_{i,j}^t|^2$ is controlled by the spring energy $E_{i,j,t}^\text{spr}(u)$. Therefore, by the Cauchy-Schwartz inequality, the average of $\widetilde{z}_{i,j}^t$ is controlled by the averaged spring energy:
		\begin{equation}\label{eqn:avg-spring-bound}
			\Big|\frac{1}{k^2} \sum_{t \in \mathcal{T}} \sum_{i,j=0}^{k-1} \widetilde{z}_{i,j}^t \Big| \leq \frac{\sqrt{|\mathcal{T}|} }{\sqrt{c}}\sqrt{\frac{1}{k^2} \sum_{t\in \mathcal{T}} \sum_{i,j=0}^{k-1} E_{i,j,t}^\text{spr}(u)}.
		\end{equation}
		
		We turn now to the second term on the left side of \eqref{eqn:average-iso-lower-bound}, which clearly satisfies 
		\begin{equation} \label{eqn:second-term-start}
			\Big|\frac{1}{k^2} \sum_{t \in \mathcal{T}} \sum_{(i,j) \in OR_t} \big(R_{-\alpha} - R_{\alpha}\big) \widetilde{b}_{i,j}^t\Big| \leq \frac{M_\alpha}{k^2} \sum_{t \in \mathcal{T}} \sum_{(i,j) \in OR_t}  \Big|\widetilde{b}_{i,j}^t\Big|
		\end{equation}
		where $M_\alpha$ is the 2-norm of the matrix $R_{-\alpha} - R_\alpha$. When the elastic energy is small, $\Big|\widetilde{b}_{i,j}^t\Big|$ is of order $1$ except perhaps in a few exceptional cases (where the spring energy $E_{i,j,t}^\text{spr}$ could be large). So this term is closely related to
		\begin{equation*}
			N_f = \mbox{the number of triangles $T_{i,j,t}$ whose 
				orientations are reversed.}
		\end{equation*}
		We observe that
		\begin{enumerate}
			\item[(1)] $N_f \leq |\mathcal{T}|k^2$ since the total number of triangles under consideration is 
			$|\mathcal{T}|k^2$, and 
			
			\item[(2)] $N_f = \frac{\eta}{c_0} \sum_{i,j=0}^{k-1} \sum_{t \in \mathcal{T}} E_{i,j,t}^{\eta, \text{pen}}(u)$ with $c_0 = \frac{\sqrt{3}}{4}$ for Kagome and $\frac{1}{2}$ for Rotating Squares. (The constant $c_0$ is the area factor in the definition of $E^{\eta,\text{pen}}_{i,j,t}$, see \eqref{eqn:energy-ijs-kagome} and \eqref{eqn:energy-ijs-rs}); recall that our penalty is $\frac{c_0}{\eta}$ for each triangle that changes its orientation.)
		\end{enumerate}
		
		We also note that $|\widetilde{b}_{i,j}^t| \leq 1 + \sqrt{E_{i,j,t}^\text{spr}(u)}$. Using this, we can  control the right side of \eqref{eqn:second-term-start} as follows:
		\begin{equation*}
			\begin{aligned}
				&\Big|\frac{1}{k^2} \sum_{t \in \mathcal{T}} \sum_{(i,j) \in OR_t} \big(R_{-\alpha} - R_{\alpha}\big) \widetilde{b}_{i,j}^t\Big| \leq \frac{M_\alpha N_f}{k^2} + \frac{M_\alpha}{k^2} \sum_{t \in \mathcal{T}} \sum_{(i,j) \in OR_t} \sqrt{E_{i,j,t}^\text{spr}(u)}\\        \leq \;& M_\alpha \frac{N_f}{k^2} + \frac{M_\alpha}{k^2} \sqrt{N_f} \sqrt{\sum_{t \in \mathcal{T}} \sum_{(i,j) \in OR_t} E_{i,j,t}^\text{spr}(u)} \leq M_\alpha \frac{N_f}{k^2} + M_\alpha \sqrt{|\mathcal{T}|} \sqrt{\frac{1}{k^2} \sum_{t \in \mathcal{T}} \sum_{i,j=0}^{k-1} E_{i,j,t}^\text{spr}(u)},
			\end{aligned}
		\end{equation*}
		where the last inequality holds due to observation (1). 
		The preceding result implies a bound in terms of the 
		averaged penalty and spring energies:
		\begin{equation*}
			\begin{aligned}
				&\Big|\frac{1}{k^2} \sum_{t \in \mathcal{T}} \sum_{(i,j) \in OR_t} \big(R_{-\alpha} - R_{\alpha}\big) \widetilde{b}_{i,j}^t\Big| \leq M_\alpha \sqrt{|\mathcal{T}|} \sqrt{\frac{N_f}{k^2}} + M_\alpha \sqrt{|\mathcal{T}|} \sqrt{\frac{1}{k^2} \sum_{t \in \mathcal{T}} \sum_{i,j=0}^{k-1} E_{i,j,t}^\text{spr}(u)}\\
				= \; & M_\alpha \sqrt{|\mathcal{T}|} \sqrt{\frac{\eta}{c_0}} \sqrt{\frac{1}{k^2} \sum_{t\in\mathcal{T}}\sum_{i,j=0}^{k-1} E_{i,j,t}^{\eta,\text{pen}}(u)} + M_\alpha \sqrt{|\mathcal{T}|}\sqrt{\frac{1}{k^2} \sum_{t \in \mathcal{T}} \sum_{i,j=0}^{k-1} E_{i,j,t}^\text{spr}(u)},
			\end{aligned}
		\end{equation*}
		where in the first line we used observation (1) to know that 
		$\frac{N_f}{k^2} \leq \sqrt{|\mathcal{T}|}\sqrt{\frac{N_f}{k^2}}$,
		and in the second line we used observation (2).
		Thus by choosing the threshold $\eta_0$ so that 
		\begin{equation} \label{eqn:eta0}
			\frac{\eta_0}{c_0} \leq 1
		\end{equation}
		we obtain the desired control of the second term on the left of \eqref{eqn:average-iso-lower-bound}:
		\begin{equation} \label{eqn:second-term-end}
			\begin{aligned}
				\Big|\frac{1}{k^2} \sum_{t \in \mathcal{T}} \sum_{(i,j) \in OR_t} \big(R_{-\alpha} - R_{\alpha}\big) \widetilde{b}_{i,j}^t\Big|
				&\leq M_\alpha\sqrt{|\mathcal{T}|}\Bigg(\sqrt{\frac{1}{k^2} \sum_{t \in \mathcal{T}} \sum_{i,j=0}^{k-1} E_{i,j,t}^\text{spr}(u)} + \sqrt{\frac{1}{k^2} \sum_{t\in\mathcal{T}}\sum_{i,j=0}^{k-1} E_{i,j,t}^{\eta,\text{pen}}(u)} \Bigg)\\
				&\leq \sqrt{2} M_\alpha \sqrt{|\mathcal{T}|} \sqrt{|U|} \sqrt{ \overline{E^\eta}(\lambda,\psi,kU) },
			\end{aligned}
		\end{equation}
		where in the second line we have used 
		\eqref{eqn:avg-energy-rewritten}. Combining this with  \eqref{eqn:avg-spring-bound} (and noting that in both estimates, the constants on the right are independent of $\lambda$, $\psi$, $k$, and $\eta$) we obtain the desired inequality \eqref{eqn:average-iso-lower-bound}. This completes our proof of the proposition. 
	\end{proof}
	
	\subsection{The compressive bound}\label{subsec:compress-bd}
	Our goal in this subsection is to prove Proposition \ref{prop:compressive-bd}. As a warm-up,
	we start by proving the compressive lower bound \eqref{eqn:compressive-bd} in the 
	simplest case: for the Rotating Squares metamaterial, when
	$\lambda = \text{diag}(\lambda_1, \lambda_2)$ is diagonal. While the argument is 
	relatively simple in this case, it reveals some (though not all) the ideas needed
	to prove Proposition \ref{prop:compressive-bd}.
	
	Here is a precise statement of our warm-up result.
	
	\begin{proposition}\label{prop:horizontal-bound}
		For the Rotating Squares metamaterial, consider a deformation of the form 
		$u(x)=\lambda x + \psi(x)$ with $\psi$ $k$-periodic, and assume that 
		$\lambda = \text{diag}(\lambda_1, \lambda_2)$ with $\lambda_1, \lambda_2 \geq 0$.
		Then we have 
		\begin{equation}\label{eqn:rs-diag-bd}
			\frac{1}{k^2|\mathcal{T}|} \sum_{i,j=0}^{k-1} \sum_{t\in \mathcal{T}} \Bigg[\Big(|\widetilde{b}_{i,j}^t|-1\Big)^2 + \Big(|\widetilde{r}_{i,j}^t|-1\Big)^2\Bigg] \geq (\lambda_1 - 1)_+^2 + (\lambda_2 - 1)_+^2,
		\end{equation}
		where $\widetilde{b}_{i,j}^t, \widetilde{r}_{i,j}^t$ are the images of the horizontal and vertical vectors in the triangle $T_{i,j,t}$  . 
	\end{proposition}
	
	\begin{proof}
		We shall take advantage of the fact that the Rotating Squares system has straight lines
		of horizontal springs and straight lines of vertical springs. If $\lambda_1 > 1$ (resp. $\lambda_2 > 1)$) then there is macroscopic horizontal (resp. vertical) stretching, 
		and we expect that this requires microscopic horizontal (resp. vertical) stretching. We shall turn this intuition into an energy bound by an application of Jensen's inequality.
		
		We first combine \eqref{eqn:avg-br-def} and \eqref{eqn:av-relation} with
		Lemma \ref{lemma:atilde-vs-lambda} to see that
		\begin{equation}\label{eqn:lambda-e1-average}
			\lambda e_1 = \frac{1}{|\mathcal{T}|}\widetilde{a}_1 =\frac{1}{k^2|\mathcal{T}|}\sum_{i,j=0}^{k-1} \sum_{t \in \mathcal{T}} \widetilde{b}_{i,j}^t.
		\end{equation}
		Since $(|\xi| - 1)_+^2$ is a convex function of $\xi \in \mathbb{R}^2$. Jensen's inequality
		gives
		\begin{equation}\label{eqn:avg-rs-e1}
			(\lambda_1 - 1)_+^2 = (|\lambda e_1| - 1)_+^2 = \Big( \Big|\frac{1}{k^2|\mathcal{T}|} \sum_{i,j=0}^{k-1} \sum_{t \in \mathcal{T}} \widetilde{b}^t_{i,j}\Big| - 1\Big)_+^2 \leq \frac{1}{k^2 |\mathcal{T}|} \sum_{i,j=0}^{k-1} \sum_{t \in \mathcal{T}} \Big(|\widetilde{b}^t_{i,j}| - 1\Big)_+^2.
		\end{equation}
		The same argument applies to the vertical springs $\widetilde{r}^t_{i,j}$: we have
		\begin{equation}\label{eqn:lambda-e2-average}
			\lambda e_2 = \frac{1}{|\mathcal{T}|}\widetilde{a}_2 =\frac{1}{k^2|\mathcal{T}|}\sum_{i,j=0}^{k-1} \sum_{t \in \mathcal{T}} \widetilde{r}_{i,j}^t .
		\end{equation}
		and Jensen's inequality gives 
		\begin{equation}\label{eqn:avg-rs-e2}
			(\lambda_2 - 1)_+^2 = (|\lambda e_2| - 1)_+^2 \leq \frac{1}{k^2 |\mathcal{T}|} \sum_{i,j=0}^{k-1} \sum_{t \in \mathcal{T}} \Big(|\widetilde{r}^t_{i,j}| - 1\Big)_+^2.
		\end{equation}
		Adding \eqref{eqn:avg-rs-e1} and \eqref{eqn:avg-rs-e2} gives the desired 
		bound \eqref{eqn:rs-diag-bd}.
	\end{proof}
	
	To prove our compressive bound for general $\lambda$ and for Kagome as well as Rotating
	Squares we must deal with an issue that our warm-up avoided: the principal directions of
	$\lambda$ will not, in general, be directions in which the lattice has a straight line 
	of springs. Our isotropic bound (Proposition \ref{prop:isotropic-bd}) will help us deal with this
	issue. 
	
	\begin{proof}[Proof of Proposition \ref{prop:compressive-bd}]
		
		We begin with the Kagome metamaterial. As usual, we let $\lambda_1, \lambda_2$ be the principal stretches associated with $\lambda$ (which can now be any $2 \times 2$ matrix),
		and we assume without loss of generality that $\lambda_1 \geq \lambda_2 \geq 0$. Our 
		argument involves two key assertions: we will show that when $u(x) = \lambda x + \psi(x)$ with $\psi$ $k$-periodic,
		\begin{equation}\label{eqn:non-diag-1}
			\frac{1}{k^2|\mathcal{T}|} \sum_{i,j=0}^{k-1} \sum_{t \in \mathcal{T}} E_{i,j,t}^\text{spr}(u) \geq \left(\sqrt{\frac{3}{4} \lambda_1^2 + \frac{1}{4} \lambda_2^2} - 1\right)_+^2
		\end{equation}
		(where $E_{i,j,t}^{\text{spr}}(u)$ is defined by \eqref{eqn:energy-ijs-kagome}); and 
		we will show that 
		\begin{equation}\label{eqn:non-diag-2}
			(\lambda_1 - \lambda_2)^2 + \left(\sqrt{\frac{3}{4} \lambda_1^2 + \frac{1}{4} \lambda_2^2} - 1\right)_+^2 \geq \frac{1}{4} \Big[(\lambda_1 - 1)_+^2 + (\lambda_2 - 1)_+^2\Big].
		\end{equation}
		Before proving these assertions, let us explain why they imply the desired result 
		\eqref{eqn:compressive-bd}. Combining Proposition \ref{prop:isotropic-bd} with the elementary 
		fact that $(\lambda_1 + \lambda_2)^2 \geq (\lambda_1-\lambda_2)^2$, we know (for
		$\eta < \eta_0$) that
		$$
		\frac{1}{C_1} \overline{E^\eta}(\lambda, \psi, kU) \geq (\lambda_1 - \lambda_2)^2 .
		$$
		Also, combining \eqref{eqn:non-diag-1} with 
		\eqref{eqn:avg-energy-rewritten}--\eqref{eqn:energy-ijs-kagome} we have
		$$
		\frac{|U|}{|\mathcal{T}|} \overline{E^\eta}(\lambda, \psi, kU) \geq
		\left(\sqrt{\frac{3}{4} \lambda_1^2 + \frac{1}{4} \lambda_2^2} - 1\right)_+^2 .
		$$
		By adding these inequalities then using \eqref{eqn:non-diag-2} we conclude that
		$$
		4 \left( \frac{1}{C_1} + \frac{|U|}{|\mathcal{T}|} \right)  
		\overline{E^\eta}(\lambda, \psi, kU) \geq 
		\Big[(\lambda_1 - 1)_+^2 + (\lambda_2 - 1)_+^2\Big] .
		$$
		Thus \eqref{eqn:compressive-bd} holds with 
		$C_2 = \frac{1}{4} \left( \frac{1}{C_1} + \frac{|U|}{|\mathcal{T}|} \right)^{-1}$.
		
		Turning now to the proof of \eqref{eqn:non-diag-1}, 
		we recall that $\lambda e_1, \lambda e_2, \lambda e_3$ are the averages of $\widetilde{b}_{i,j}^t$, $\widetilde{r}_{i,j}^t$, and $\widetilde{r}_{i,j}^t - \widetilde{b}_{i,j}^t$, where $e_1, e_2, e_3$ are the unit vectors in the horizontal, 60 and 120 degree directions. Therefore, by applying Jensen's inequality in the three lattice directions we obtain
		\begin{align} \label{eqn:kagome-lambda-bd}
			\frac{1}{k^2|\mathcal{T}|} \sum_{i,j=0}^{k-1} \sum_{t \in \mathcal{T}} \Big(|\widetilde{b}^t_{i,j}| - 1\Big)^2& \geq (|\lambda e_1|-1)_+^2, \nonumber\\
			\frac{1}{k^2|\mathcal{T}|} \sum_{i,j=0}^{k-1} \sum_{t \in \mathcal{T}}\Big(|\widetilde{r}^t_{i,j}| - 1\Big)^2 & \geq (|\lambda e_2|-1)_+^2,\\
			\frac{1}{k^2|\mathcal{T}|} \sum_{i,j=0}^{k-1} \sum_{t \in \mathcal{T}}\Big(|\widetilde{r}^t_{i,j} -\widetilde{b}^t_{i,j} | - 1\Big)^2& \geq (|\lambda e_3|-1)_+^2. \nonumber
		\end{align}
		Adding, we get a lower bound for the averaged spring energy:
		\begin{equation}\label{eqn:app-avg-comp-1}
			\frac{1}{k^2|\mathcal{T}|} \sum_{i,j=0}^{k-1} \sum_{t \in \mathcal{T}} E_{i,j,t}^\text{spr}(u) \geq \Big[(|\lambda e_1|-1)_+^2 + (|\lambda e_2|-1)_+^2 + (|\lambda e_3|-1)_+^2\Big]. 
		\end{equation}
		We claim now that 
		\begin{equation}\label{eqn:non-diagonal-1}
			(|\lambda e_1|-1)_+^2 + (|\lambda e_2|-1)_+^2 + (|\lambda e_3|-1)_+^2 \geq \left(\sqrt{\frac{3}{4} \lambda_1^2 + \frac{1}{4} \lambda_2^2} - 1\right)_+^2 \quad \text{ when } \lambda_1 \geq \lambda_2 \geq 0 .
		\end{equation}
		Note that \eqref{eqn:non-diag-1} follows directly from \eqref{eqn:app-avg-comp-1} and \eqref{eqn:non-diagonal-1}. 
		
		To prove \eqref{eqn:non-diagonal-1}, we use a singular value decomposition of $\lambda = USV^T$ with $U,V \in O(2)$, $S = \text{diag}(\lambda_1, \lambda_2)$ and $\lambda_1 \geq \lambda_2$. We claim that there exists $\widetilde{U}, \widetilde{V} \in SO(2)$ such that $\lambda = \widetilde{U}D\widetilde{V}^T$ and $D = \text{diag}(d_1, d_2)$ where $|d_1| = \lambda_1$ and $|d_2|  = \lambda_2$. In fact, when the orthonormal matrix $U$ has negative determinant, we choose $\widetilde{U} = U F$ with $F = \text{diag}(1,-1)$ and $D = FS$; similarly, when the orthonormal matrix $V$ has negative determinant, we choose $\widetilde{V} = V F$ and $D = SF$.
		
		Continuing the proof of \eqref{eqn:non-diagonal-1}, we use the decomposition $\lambda = \widetilde{U} D \widetilde{V}^T$ to see that
		\begin{equation*}
			|\lambda e_1|= |\widetilde{U} D \widetilde{V}^T e_1|= |D \widetilde{V}^T e_1| =\sqrt{d_2^2 + (d_1^2 - d_2^2) \cos^2 \theta} = \sqrt{\lambda_2^2 + (\lambda_1^2 - \lambda_2^2) \cos^2 \theta},
		\end{equation*}
		where $\widetilde{V}^T e_1 = (\cos\theta, \sin \theta)$. We also get similar formulas for $|\lambda e_2|$ and $|\lambda e_3|$:
		\begin{equation*}
			|\lambda e_2| = \sqrt{\lambda_2^2  + (\lambda_1^2 - \lambda_2^2) \cos^2 (\theta+\frac{\pi}{3})},\qquad
			|\lambda e_3| = \sqrt{\lambda_2^2 + (\lambda_1^2 - \lambda_2^2) \cos^2 (\theta+\frac{2\pi}{3})}.
		\end{equation*}
		Thus the left hand side of \eqref{eqn:non-diagonal-1} becomes 
		\begin{equation}\label{eqn:non-diagonal-desired-form}
			\begin{aligned}
				&\big(|\lambda e_1| - 1)_+^2 + (|\lambda e_2| - 1)_+^2 + (|\lambda e_3| - 1\big)_+^2 = \big(\sqrt{\lambda_2^2 + (\lambda_1^2 - \lambda_2^2) \cos^2 \theta} - 1\big)_+^2 \\
				& + \left(\sqrt{\lambda_2^2  + (\lambda_1^2 - \lambda_2^2) \cos^2 (\theta+\frac{\pi}{3})}-1\right)_+^2 + \left(\sqrt{\lambda_2^2 + (\lambda_1^2 - \lambda_2^2) \cos^2 (\theta+\frac{2\pi}{3})}-1\right)_+^2. 
			\end{aligned}
		\end{equation}
		Since $\lambda_1 \geq \lambda_2$, the left hand side of \eqref{eqn:non-diagonal-1} is lower bounded by the largest of the three terms on the right hand side of \eqref{eqn:non-diagonal-desired-form}:
		\begin{equation*}
			\big(|\lambda e_1| - 1)_+^2 + (|\lambda e_2| - 1)_+^2 + (|\lambda e_3| - 1\big)_+^2 \geq \left(\sqrt{\lambda_2^2 + (\lambda_1^2 - \lambda_2^2) t(\theta)} - 1\right)_+^2,
		\end{equation*}
		where $t(\theta) = \max\big\{\cos^2 \theta, \cos^2(\theta + \frac{\pi}{3}), \cos^2(\theta + \frac{2\pi}{3})\big\}$. By a standard calculation that separates the region $[0,2\pi]$ into pieces where only one of $\cos^2 \theta, \cos^2(\theta + \frac{\pi}{3}), \cos^2(\theta + \frac{2\pi}{3})$ dominates, one verifies that $t(\theta) \geq \frac{3}{4}$ for all $\theta \in [0,2\pi)$. Combining this with the monotonicity of the function $x \in \mathbb{R}_+ \rightarrow (\sqrt{x}-1)_+^2$ gives the desired inequality\eqref{eqn:non-diagonal-1}, completing the proof of our first key assertion 
		\eqref{eqn:non-diag-1}.
		
		Our other key assertion \eqref{eqn:non-diag-2} is just a fact about any pair of real numbers satisfying $\lambda_1 \geq \lambda_2 \geq 0$. We prove it by considering the three cases (a) $\lambda_1, \lambda_2 \leq 1$, 
		(b) $\lambda_1 \geq 1 \geq \lambda_2$, and (c) $\lambda_1,\lambda_2 \geq 1$.
		Case (a) is trivial, because the right hand side of \eqref{eqn:non-diag-2} vanishes. Case (b) is also easy, since in this case 
		\begin{equation*}
			(\lambda_1 - \lambda_2)^2 +\Bigg(\sqrt{\frac{3}{4}\lambda_1^2 + \frac{1}{4} \lambda_2^2} - 1\Bigg)_+^2 \geq (\lambda_1 - \lambda_2)^2 \geq (\lambda_1 - 1)^2 = (\lambda_1 - 1)_+^2 + (\lambda_2 - 1)_+^2.
		\end{equation*}
		As for case (c): when $\lambda_1 \geq \lambda_2 \geq 1$ we have 
		$\sqrt{\frac{3}{4}\lambda_1^2 + \frac{1}{4} \lambda_2^2} \geq \lambda_2 \geq 1$, so 
		\begin{equation*}
			\begin{aligned}
				(\lambda_1 - \lambda_2)^2 + \Bigg(\sqrt{\frac{3}{4}\lambda_1^2 + \frac{1}{4} \lambda_2^2} - 1\Bigg)_+^2 &= (\lambda_1 - \lambda_2)^2 + \Bigg(\sqrt{\frac{3}{4}\lambda_1^2 + \frac{1}{4} \lambda_2^2} - 1\Bigg)^2 \geq (\lambda_1 - \lambda_2)^2 + (\lambda_2 - 1)^2\\
				& \geq \frac{1}{2}(\lambda_1 - 1)^2 \geq \frac{1}{4} \bigg[(\lambda_1 - 1)_+^2 + (\lambda_2-1)_+^2\bigg].
			\end{aligned}
		\end{equation*}
		(The first inequality on the second line follows from $x^2 + y^2 \geq \frac{1}{2}(x + y)^2$, taking $x = \lambda_1 - \lambda_2$ and $y = \lambda_2 - 1$); the last inequality is obvious since $\lambda_1 \geq \lambda_2 \geq 1$.) This establishes our second key assertion \eqref{eqn:non-diag-2}, thereby completing the proof of \eqref{eqn:compressive-bd} for the Kagome metamaterial.
		
		For the Rotating Squares metamaterial, the argument is only slightly different. It relies on two key assertions which are analogous to 
		\eqref{eqn:non-diag-1}--\eqref{eqn:non-diag-2}: 
		\begin{equation}\label{eqn:avg-spr-bd-rs}
			\frac{1}{k^2|\mathcal{T}|} \sum_{i,j=0}^{k-1} \sum_{t \in \mathcal{T}} E_{i,j,t}^\text{spr}(u) \geq \left(\sqrt{\frac{1}{2} \lambda_1^2 + \frac{1}{2} \lambda_2^2} - 1\right)_+^2,
		\end{equation}
		where $E_{i,j,t}^\text{spr}(u)$ is defined by \eqref{eqn:energy-ijs-rs}; and
		\begin{equation}\label{eqn:rs-lambda-bd}
			(\lambda_1 - \lambda_2)^2 + \left(\sqrt{\frac{1}{2} \lambda_1^2 + \frac{1}{2} \lambda_2^2} - 1\right)_+^2 \geq \frac{1}{4} \Big[(\lambda_1 - 1)_+^2 + (\lambda_2 - 1)_+^2\Big] 
		\end{equation}
		when $\lambda_1 \geq \lambda_2 \geq 0$. The proof of the latter is entirely parallel to the justification we just gave for \eqref{eqn:non-diag-2}, so it
		can safely be left to the reader. The proof that these inequalities imply
		\eqref{eqn:compressive-bd} is parallel to the argument by which we deduced
		\eqref{eqn:compressive-bd} from \eqref{eqn:non-diag-1}--\eqref{eqn:non-diag-2}, so this too is left to the reader. We will, however, provide the proof
		of \eqref{eqn:avg-spr-bd-rs}. We start by applying Jensen's inequality to get
		\begin{equation}\label{eqn:app-comp-rs}
			\frac{1}{k^2 |\mathcal{T}|} \sum_{i,j=0}^{k-1} \sum_{t \in \mathcal{T}} (|\widetilde{b}_{i,j}^t|-1)^2 + (|\widetilde{r}_{i,j}^t| - 1)^2\geq (|\lambda e_1| - 1)_+^2 + (|\lambda e_2| - 1)_+^2,
		\end{equation}
		where $e_1 = (1,0)^T$ and $e_2 = (0,1)^T$. Using the special singular value decomposition introduced above ($\lambda = \widetilde{U} D \widetilde{V}^T$ with $\widetilde{U}, \widetilde{V} \in SO(2)$, $D = \text{diag}(d_1,d_2)$ and $|d_1| = \lambda_1, |d_2| = \lambda_2$), we calculate 
		$\lambda e_1$ and $\lambda e_2$ directly:
		\begin{equation*}
			|\lambda e_1| = |D V^T e_1| = \sqrt{\lambda_2^2 + (\lambda_1^2 - \lambda_2^2) \cos^2\theta},\qquad |\lambda e_2| = |D R_{\frac{\pi}{2}} V^T e_1| = \sqrt{\lambda_2^2 + (\lambda_1^2 - \lambda_2^2) \cos^2(\theta+\frac{\pi}{2})}
		\end{equation*}
		where $V^T e_1 = (\cos\theta, \sin \theta)^T$. By arguments parallel to those used for Kagome we get
		\begin{equation*}
			(|\lambda e_1| - 1)_+^2 + (|\lambda e_2| - 1)_+^2 \geq \left(\sqrt{\lambda_2^2 + (\lambda_1^2 - \lambda_2^2) t(\theta)} - 1\right)_+^2,
		\end{equation*}
		where $t(\theta) = \max\big\{\cos^2 \theta, \cos^2(\theta + \frac{\pi}{2})\big\}$. By a standard calculation, one finds that $t(\theta) \geq \frac{1}{2}$ for all $\theta \in [0,2\pi)$. Combining these elements gives 
		\begin{equation}\label{eqn:comp-rs-1}
			\frac{1}{k^2|\mathcal{T}|}\sum_{i,j=0}^{k-1} \sum_{t \in \mathcal{T}} E^{\text{spr}}_{i,j,t}(u) \geq \frac{1}{k^2 |\mathcal{T}|} \sum_{i,j=0}^{k-1} \sum_{t \in \mathcal{T}} \Big(|\widetilde{b}_{i,j}^t|-1\Big)^2 + \Big(|\widetilde{r}_{i,j}^t| - 1\Big)^2 \geq \left(\sqrt{\frac{1}{2}\lambda_1^2 + \frac{1}{2}\lambda_2^2} - 1\right)_+^2 .
		\end{equation}
		This completes our justification of \eqref{eqn:avg-spr-bd-rs}, thereby completing the proof of \eqref{eqn:compressive-bd} for the Rotating Squares metamaterial.
	\end{proof}
	
	\section{Other conformal metamaterials}\label{sec:other-conformal}
	
	The methods developed in this paper are applicable to many other examples besides the Kagome and Rotating Squares metamaterials. To show this, we now present a few examples. They share the properties that 
	\begin{enumerate}
		\item[(i)] the framework of \cite{li2025effective} is applicable, so there is a well-defined effective energy $\overline{W}^\eta$; and 
		\item[(ii)] the methods of this paper are also applicable, giving
		direct analogues of the results proved in sections \ref{sec:geo-argument} and \ref{sec:lower-bd} for the Kagome and Rotating Squares metamaterials.
	\end{enumerate}
	
\begin{figure}[!htb]
		\centering
		{\includegraphics[width=0.85\linewidth]{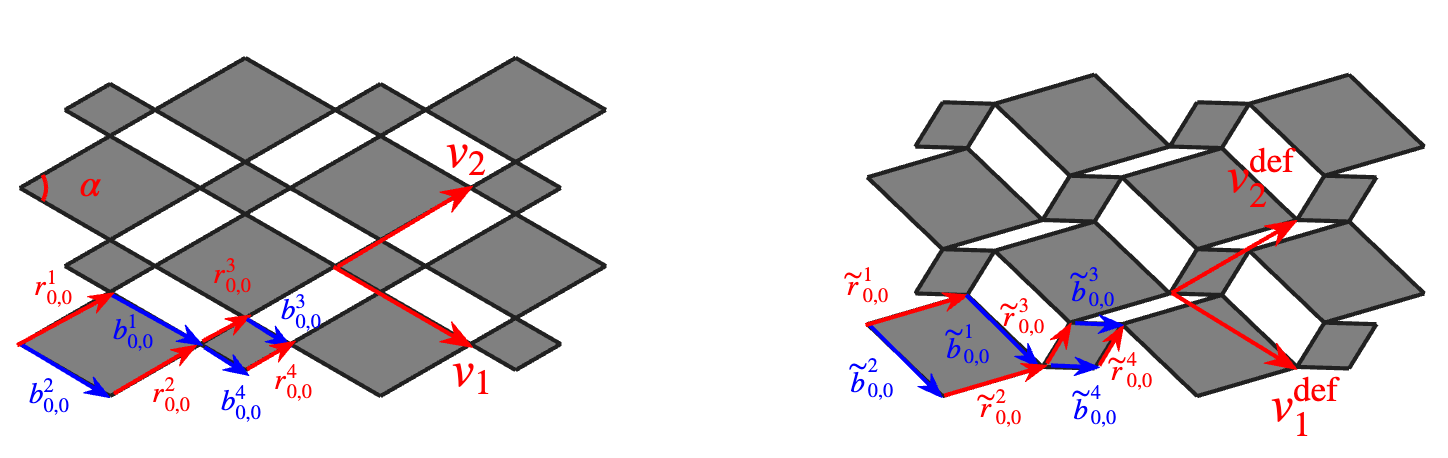}}
		\caption{A variant of the Rotating Squares metamaterial and its periodic mechanism: the left figure shows the reference state, and the right figure shows the deformed state.}
		\label{fig:deformed-rs-mechanism}
	\end{figure}

	Elaborating briefly upon (ii): for all our examples, the macroscopic effect of a $k$-periodic mechanism is an isotropic compression (thus: the analogue of Proposition \ref{lemma-geometric} holds)\footnote{In certain conformal metamaterials, such as the one shown on the right in Figure \ref{fig:deformed-rs-mechanism}, the parallelogram slits do not close simultaneously, implying that the compression ratio cannot reach zero. In these cases, when establishing an argument analogous to Proposition \ref{lemma-geometric}, the compression ratio will be restricted to a subset of $[0,1]$ instead of $[0,1]$.}; and for all our examples, the effective energy satisfies the lower bound \eqref{eqn:eff-lower-bd} (thus: the analogue of Theorem \ref{thm:lower-bound} holds). In particular, our examples are all conformal metamaterials. 
	
	We shall be rather informal: instead of providing detailed proofs for 
	each example, we'll just provide enough information to make it clear 
	how our methods can be used.
	\medskip 
	
	\noindent
	{\sc A variant of the Rotating Squares metamaterial.}
	Our first example is shown in Figure \ref{fig:deformed-rs-mechanism}. Its elastic regions are not squares but rather rhombuses (that is, parallelograms with four equal sides). There are in fact two families of rhombuses, which must  be \emph{scaled copies of one another}. It is important that 
	the edges of the rhombuses form two families of straight lines. Of course we
	want to view this example as a spring network: each edge of each rhombus is 
	a spring, and each rhombus also has a diagonal spring. 
	
	To show that the macroscopic effect of a periodic mechanism is an isotropic compression, it is convenient to use the choice of 
	$b_{i,j}^t, r_{i,j}^t$ (and the associated $\widetilde{b}_{i,j}, \widetilde{r}_{i,j}^t$) shown in Figure \ref{fig:deformed-rs-mechanism}. 
	For a $k$-periodic mechanism, there is no spring energy and each rhombus has its orientation preserved. Therefore we have
	\begin{equation*}
		r_{i,j}^t = R_\alpha b_{i,j}^t, \qquad \widetilde{r}_{i,j}^t = R_\alpha \widetilde{b}_{i,j}^t, \quad i,j=0,\dots,k-1, \quad t\in \mathcal{T},
	\end{equation*}
	where $\mathcal{T} = \{1,2,3,4\}$. Moreover, the averaged vectors defined by \eqref{eqn:avg-br-ref} and \eqref{eqn:avg-br-def} satisfy
	\begin{equation*}
		a_2 = R_\alpha a_1, \qquad \widetilde{a}_2 = R_\alpha \widetilde{a}_1.
	\end{equation*}
	Using the relationship $\widetilde{a}_1 = \lambda a_1$ and $\widetilde{a}_2 = \lambda a_2$ (whose proof follows that of Lemma \ref{lemma:atilde-vs-lambda}) we conclude (exactly as in Section \ref{sec:geo-argument}) that if $u(x) = \lambda x + \psi(x)$ is a periodic mechanism then
	$\lambda = cR$ for some $c \geq 0$ and $R \in SO(2)$.  
	To see that $c \leq 1$ we can again argue as in Section \ref{sec:geo-argument}: since a mechanism preserves the length of each spring, the straight lines of springs in the $b$ and $r$ directions are mapped, in general, to zig-zags (each piece keeping its length). The macroscopic effect of such a map cannot be extension, so $c \leq 1$.
	
	Turning now to our lower bound on the effective energy: the arguments we 
	used for the isotropic lower bound (Proposition \ref{prop:isotropic-bd}) extend straightforwardly to this example. However, our proof of the 
	compressive lower bound for Kagome and Rotating Squares (Proposition \ref{prop:compressive-bd}) used that those structures have lines of 
	springs \emph{that all have the same length} -- which isn't true for
	the structure under discussion now.
	
	The argument we must generalize is the one that led to \eqref{eqn:kagome-lambda-bd} and \eqref{eqn:app-comp-rs}. To explain the generalization, we focus on the vectors $b_{i,j}^t$ (the treatment of $r_{i,j}^t$ being entirely analogous). Our situation is that
	all the $b_{i,j}^t$ are multiples of a single unit vector (which we shall call $e_1$), but their lengths
	$$
	\ell_{i,j}^t =: |b_{i,j}^t|
	$$
	are not all the same. We shall show that if
	\begin{equation} \label{eqn:spring-law-M}
		E^{\text{spr}}_{i.j.t} \geq M \ell_{i,j}^t \left( \frac{|\widetilde{b}^t_{i,j}|}{\ell_{i,j}^t} -1 \right)^2
	\end{equation}
	for some constant $M$, then
	\begin{equation} \label{eqn:compressive-estimate-generalized}
		(|\lambda e_1| - 1)_+^2 \leq
		\frac{1}{M\ell_\text{avg}} \, \frac{1}{k^2 |\mathcal{T}|}
		\sum_{t \in \mathcal{T}} \sum_{i,j=0}^{k-1} E^\text{spr}_{i,j,t}
	\end{equation}
	where $\ell_\text{avg}$ is the average of $\ell_{i,j}^t$:
	\begin{equation} \label{eqn:lavg-defn}
		\ell_\text{avg} =: \frac{1}{k^2 |\mathcal{T}|}
		\sum_{t \in \mathcal{T}} \sum_{i,j=0}^{k-1} \ell_{i,j}^t .
	\end{equation}
	We note that
	\begin{enumerate}
		\item[(1)] Since the metamaterial is periodic, $\ell_{i,j}^t$ is actually independent of $i$ and $j$, and
		therefore $\ell_\text{avg} = \frac{1}{|\mathcal{T}|} \sum_{t \in \mathcal{T}} \ell_{0,0}^t$ is independent
		of $k$.
		
		\item[(2)] If the spring associated with $b_{i,j}^t$ has energy $(|\widetilde{b}_{i,j}^t| - \ell_{i,j}^t)^2$ then
		(using periodicity) \eqref{eqn:spring-law-M} is satisfied with $M = \min_{t \in \mathcal{T}} \{ l_{0,0}^t \}$.
	\end{enumerate}
	To prove \eqref{eqn:compressive-estimate-generalized} we shall use the fact that $\widetilde{a}_1 = \lambda a_1$, where
	$$
	a_1 = \frac{1}{k^2} \sum_{t \in \mathcal{T}} \sum_{i,j=0}^{k-1} b_{i,j}^t = |\mathcal{T}| \ell_{\text{avg}} e_1 \quad
	\text{and} \quad \widetilde{a}_1 = \frac{1}{k^2} \sum_{t \in \mathcal{T}} \sum_{i,j=0}^{k-1} \widetilde{b}_{i,j}^t .
	$$
	Preparing to apply Jensen's inequality, we observe that
	$$
	\sum_{t \in \mathcal{T}} \sum_{i,j=0}^{k-1} \frac{\ell_{i,j}^t}{\ell_\text{avg} k^2 |\mathcal{T}|} = 1
	$$
	and
	\begin{align*}
		\sum_{t \in \mathcal{T}} \sum_{i,j=0}^{k-1}
		\frac{\ell_{i,j}^t}{\ell_\text{avg} k^2 |\mathcal{T}|} \frac{\widetilde{b}_{i,j}^t}{\ell_{i,j}^t} &=
		\frac{1}{\ell_\text{avg} k^2 |\mathcal{T}|} \sum_{t \in \mathcal{T}} \sum_{i,j=0}^{k-1} \widetilde{b}_{i,j}^t \\
		&= \frac{1}{\ell_\text{avg} |\mathcal{T}|} \lambda a_1 = \lambda e_1 .
	\end{align*}
	Therefore applying Jensen's inequality to the convex function $(|x|-1)_+^2$ and using \eqref{eqn:spring-law-M} gives
	\begin{align*}
		(|\lambda e_1| - 1)_+^2 &\leq
		\sum_{t \in \mathcal{T}} \sum_{i,j=0}^{k-1} \frac{\ell_{i,j}^t}{\ell_\text{avg} k^2 |\mathcal{T}|}
		\left( \frac{|\widetilde{b}^t_{i,j}|}{\ell_{i,j}^t} -1 \right)^2 \\
		& \leq \frac{1}{M\ell_\text{avg} k^2 |\mathcal{T}|} \sum_{t \in \mathcal{T}} \sum_{i,j=0}^{k-1} E^\text{spr}_{i,j,t} ,
	\end{align*}
	which is the desired relation \eqref{eqn:compressive-estimate-generalized}.
	
	Using \eqref{eqn:compressive-estimate-generalized} in place of \eqref{eqn:app-comp-rs}, the arguments used for our compressive
	bound (Proposition \ref{prop:compressive-bd}) extend straightforwardly to the structure shown in Figure \ref{fig:deformed-rs-mechanism}.
	\medskip
	
\begin{figure}[!htb]
		\centering
		\includegraphics[width=0.85\linewidth]{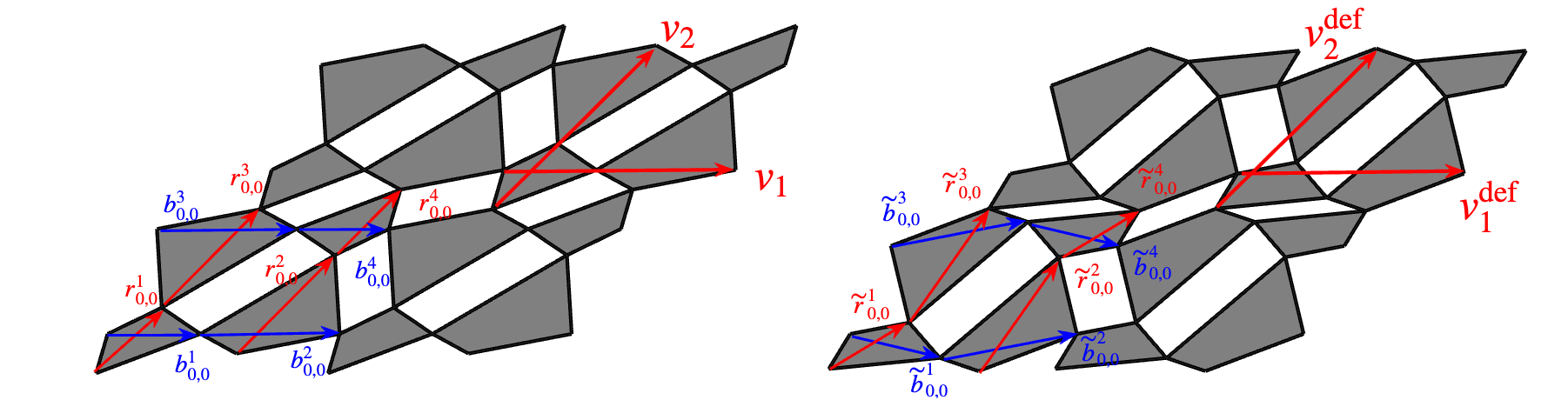}
		\caption{Another variant of the Rotating Squares metamaterial and its periodic mechanism: the left figure shows the reference state, and the right figure shows the deformed state.}
		\label{fig:deformed-quad-mechanism}
	\end{figure}

	\noindent 
	{\sc Another variant of the Rotating Squares metamaterial.} 
	Our second variant of the Rotating Squares metamaterial is shown in
	Figure \ref{fig:deformed-quad-mechanism}. Its elastic regions are specially designed quadrilaterals whose two diagonals have equal length and whose diagonals share a common angle $\alpha$ across all the quadrilaterals.
	To analyze this example it is convenient to choose the vectors $b_{i,j}^t$ and $r_{i,j}^t$ along the diagonals of the elastic regions, as shown in the figure. As usual, our metamaterial is a spring network; but this time we require that \emph{both diagonals of the elastic regions are springs} (so that each of the vectors
	$b_{i,j}^t$ and $r_{i,j}^t$ is associated with a spring). Of course, the 
	edges of the elastic regions must also be springs. 
	
	The geometric conditions stated above can be expressed as $r_{i,j}^t = R_\alpha b_{i,j}^t$. For a $k$-periodic mechanism, the associated 
	deformed vectors satisfy $\widetilde{r}_{i,j}^t = R_\alpha \widetilde{b}_{i,j}^t$. It follows by arguing as in Section \ref{sec:geo-argument}
	that a $k$-periodic mechanism of this metamaterial must be an 
	isotropic contraction.
	
	The results in Section \ref{sec:lower-bd} also extend to this example. There is, however, a new twist. One must show that each $i,j,t$, if $E^\text{spr}_{i,j,t}$ is
	small then $\widetilde{r}_{i,j}^t$ is close to $R_\alpha \widetilde{b}_{i,j}^t$ if the orientation of the associated elastic region
	is preserved. (More precisely, one needs an estimate analogous to 
	\eqref{eqn:spring-lower-bound-triangle}--\eqref{eqn:deviation-z}.) This
	requires a result analogous to Proposition \ref{prop:iso-triangle} for general 
	triangles (whereas for simplicity we stated and proved Proposition \ref{prop:iso-triangle} only for isosceles triangles).
	\medskip 
	
\begin{figure}[!htb]
		\centering
		\includegraphics[width=0.8\linewidth]{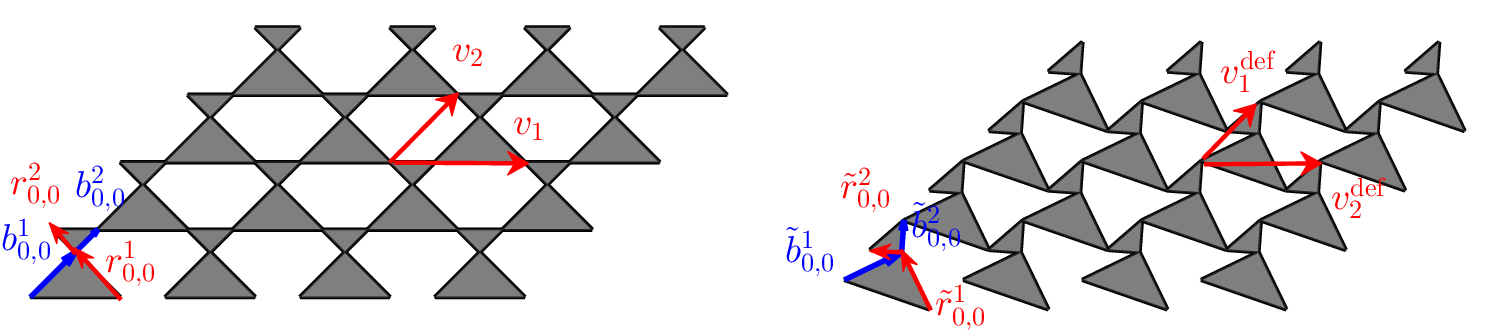}
		\caption{Our simple variant of the Kagome metamaterial and its one-periodic mechanism: the left figure shows the reference state, and the right figure shows the deformed state.}
		\label{fig:deformed-kagome-mechanism-revised}
	\end{figure}

	\noindent {\sc A simple variant of the Kagome metamaterial.} 
	Figure \ref{fig:deformed-kagome-mechanism-revised} shows a relatively simple 
	generalization of the Kagome metamaterial, in which the uniform equilateral triangles of the Kagome structure are replaced with isosceles triangles of identical shapes but different sizes. 
	
	It is convenient in this case to choose the vectors $b_{i,j}^t$ and $r_{i,j}^t$ to lie along the two equal sides of a 
	triangle, as shown in Figure \ref{fig:deformed-kagome-mechanism-revised}. Since all the triangles have the same
	shape, we have $r_{i,j}^t = R_\alpha b_{i,j}^t$ where $\alpha$ is 
	the relevant angle of the triangle. It follows, as usual, that when $u$ 
	is $k$-periodic mechanism $\widetilde{r}_{i,j}^t = R_\alpha \widetilde{b}_{i,j}^t$. We deduce using the methods of 
	Section \ref{sec:geo-argument} that a $k$-periodic mechanism must be an 
	isotropic contraction. The methods of Section \ref{sec:lower-bd} are also 
	applicable, so the effective energy satisfies our lower bound \eqref{eqn:eff-lower-bd}.
	\medskip 
	
\begin{figure}[!htb]
		\centering
		\includegraphics[width=0.8\linewidth]{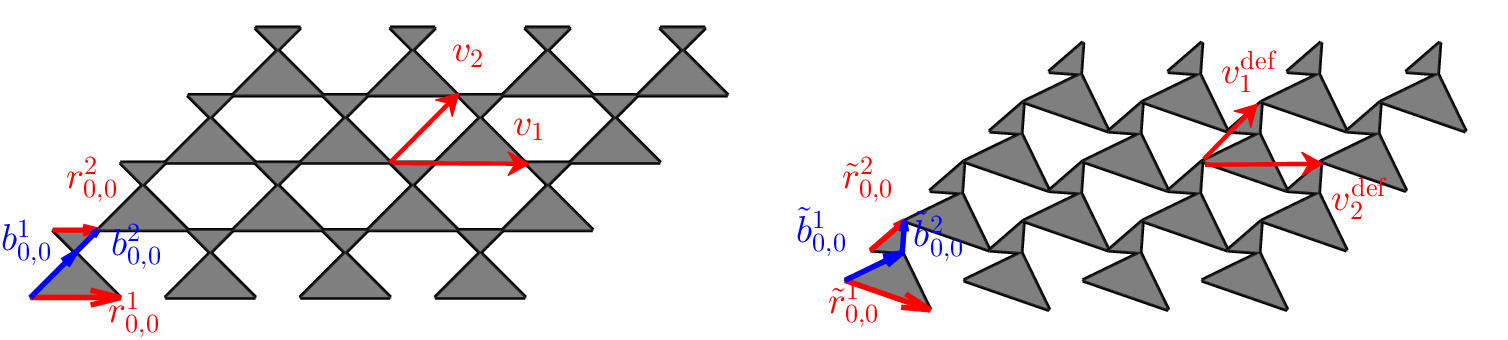}
		\caption{Our simple variant of the Kagome metamaterial with a choice
			of $b_{i,j}^t$ and $r_{i,j}^t$ satisfying \eqref{eqn:br-unequal-length} with $c \neq 1$.}
		\label{fig:deformed-kagome-mechanism}
	\end{figure}

\begin{figure}[!htb]
		\centering
		\includegraphics[width=0.85\linewidth]{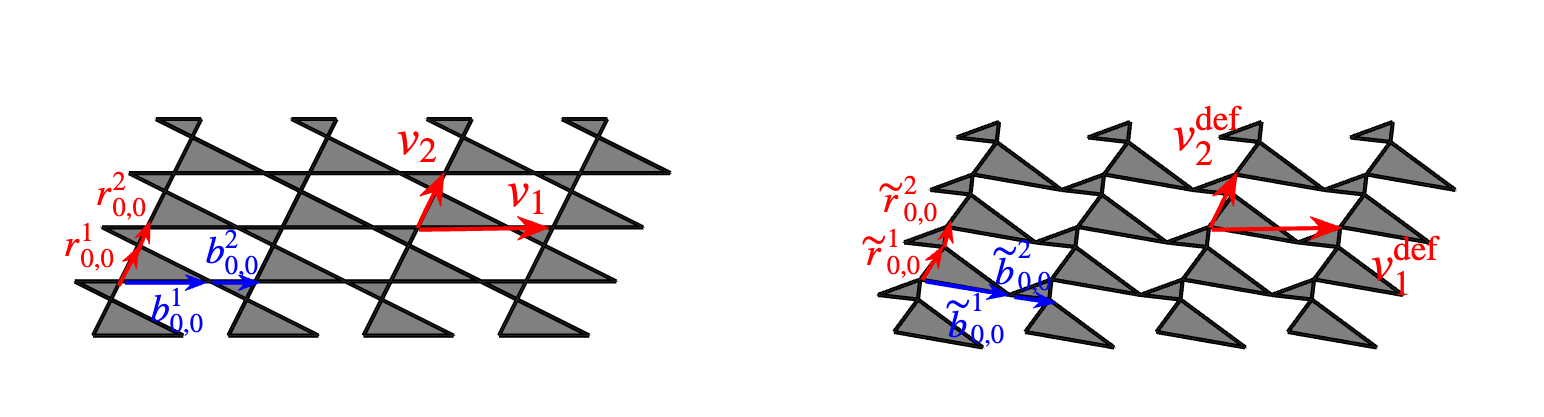}
		\caption{Our arguments apply even if the triangles are not isosceles, provided that they have the same shape. One such example is shown here, along with the image of its one-periodic mechanism.}
		\label{fig:deformed-kagome-general}
	\end{figure}

	\noindent {\sc A more general variant of the Kagome metamaterial.}
	All the variants discussed so far had the property that 
	$r_{i,j}^t = R_\alpha b_{i,j}^t$; in particular, the vectors $r_{i,j}^t$ and $b_{i,j}^t$ had the same length. Our methods work, however, under the 
	weaker condition that 
	\begin{equation} \label{eqn:br-unequal-length}
		r_{i,j}^t = c R_\alpha b_{i,j}^t \quad \mbox{for some constant $c$ that is indendent of $i$, $j$, and $t$.}
	\end{equation}
	Our simple variant of Kagome (the structure shown in 
	Figure \ref{fig:deformed-kagome-mechanism-revised}) provides
	an example if we use a different choice of $b_{i,j}^t, r_{i,j}^t$ as shown in Figure \ref{fig:deformed-kagome-mechanism}.
	However, the more flexible framework provided 
	by \eqref{eqn:br-unequal-length} permits a more general class
	of examples, since the triangles no longer need to be isosceles (they just 
	need to have the same shape). For example, our arguments apply to the 
	structure shown in Figure \ref{fig:deformed-kagome-general}.
	
	To see that the arguments in Section \ref{sec:geo-argument} are applicable when \eqref{eqn:br-unequal-length} holds, we observe that when 
	$u(x) = \lambda x + \psi(x) $ is a $k$-periodic
	mechanism the images of $b_{i,j}^t$ and $r_{i,j}^t$ must satisfy 
	\begin{equation}\label{eqn:iso-geo-2}
		\widetilde{r}_{i,j}^t = cR_\alpha \widetilde{b}_{i,j}^t, \qquad i,j = 0, \dots, k-1, \quad t \in \mathcal{T}.
	\end{equation}
	Therefore the averaged vectors $a_1, a_2$ and $\widetilde{a}_1, \widetilde{a}_2$ satisfy
	\begin{equation*}
		a_2 = c R_\alpha a_1, \qquad \widetilde{a}_2 = c R_\alpha \widetilde{a}_1.
	\end{equation*}
	We still have $\widetilde{a}_i = \lambda a_i$, $i=1,2$ (since the proof
	of Lemma \ref{lemma:atilde-vs-lambda} works when the vectors $r_{i,j}^t$ lie along a family of parallel lines through the lattice and the vectors $b_{i,j}^t$ lie along another family of parallel lines). We therefore
	get 
	\begin{equation*}
		\widetilde{a}_2 = cR_\alpha \widetilde{a}_1 = c R_\alpha \lambda a_1 = \widetilde{a}_2 \quad \mbox{and} \quad \lambda a_2 = \lambda c R_\alpha a_1 \qquad \Rightarrow \qquad \lambda R_\alpha a_1 = R_\alpha \lambda a_1.
	\end{equation*}
	Now arguing exactly as in Section \ref{sec:geo-argument} we conclude that $\lambda = c R$ for some $c \geq 0$ and $R \in SO(2)$. 
	The fact that $c \leq 1$ comes, as usual, from the fact that the structure (in its reference configuration) has straight lines of springs. 
	
	Our lower bound \eqref{eqn:eff-lower-bd} also extends to this setting. No 
	new ideas are required for its proof. 
	
	\begin{remark}[About variants of the Kagome metamaterial]
		Our variants are a subclass of a larger family of
		variants of the Kagome metamaterial, sometimes called ``deformed Kagome metamaterials'' (see e.g. \cite{rocklin2017transformable}). They are periodic structures whose 
		overall topology is similar to that of the Kagome structure, but whose unit cell has two triangles of any shape (see Figure 1 in \cite{rocklin2017transformable}).  
		
		For most deformed Kagome metamaterials, there is \textbf{only one} periodic mechanism. The examples in 
		Figure \ref{fig:deformed-kagome-mechanism-revised} -- Figure \ref{fig:deformed-kagome-general} are special: they have straight lines of springs. We expect that, like the standard Kagome metamaterial (made using with uniform equilateral triangles, as shown in Figure \ref{fig:kagome-mechanisms}), they possess infinitely many periodic mechanisms. Nevertheless, our arguments show that for these structures, every periodic mechanism is an isotropic compression.
	\end{remark}
	
	\vspace{5em}

	\noindent 
	\textit{Acknowledgements.} The authors are grateful to Katia Bertoldi and Bolei Deng for help exploring 
	the mechanical properties of the Kagome metamaterial numerically; those explorations were a major 
	inspiration for the theory presented here. Both authors gratefully acknowledge support from the Simons Foundation through its Collaboration on Extreme Wave Phenomena (Grant No. 733694), and from the National Science Foundation (Grant No. DMS-2009746).
	
	\appendix
	\section*{Appendices}
	\addcontentsline{toc}{section}{Appendix}
	\setcounter{section}{1}  
	\renewcommand{\theremark}{\Alph{section}.\arabic{remark}}
	
	\renewcommand{\thesubsection}{\Alph{subsection}}
	\renewcommand{\theequation}{\thesubsection.\arabic{equation}}

	\setcounter{equation}{0}
	\subsection{Proof of Proposition \ref{prop:iso-triangle}}\label{app:proof-iso-triangle}
\begin{figure}[tbh!]
		\centering
		\includegraphics[width=0.35\linewidth]{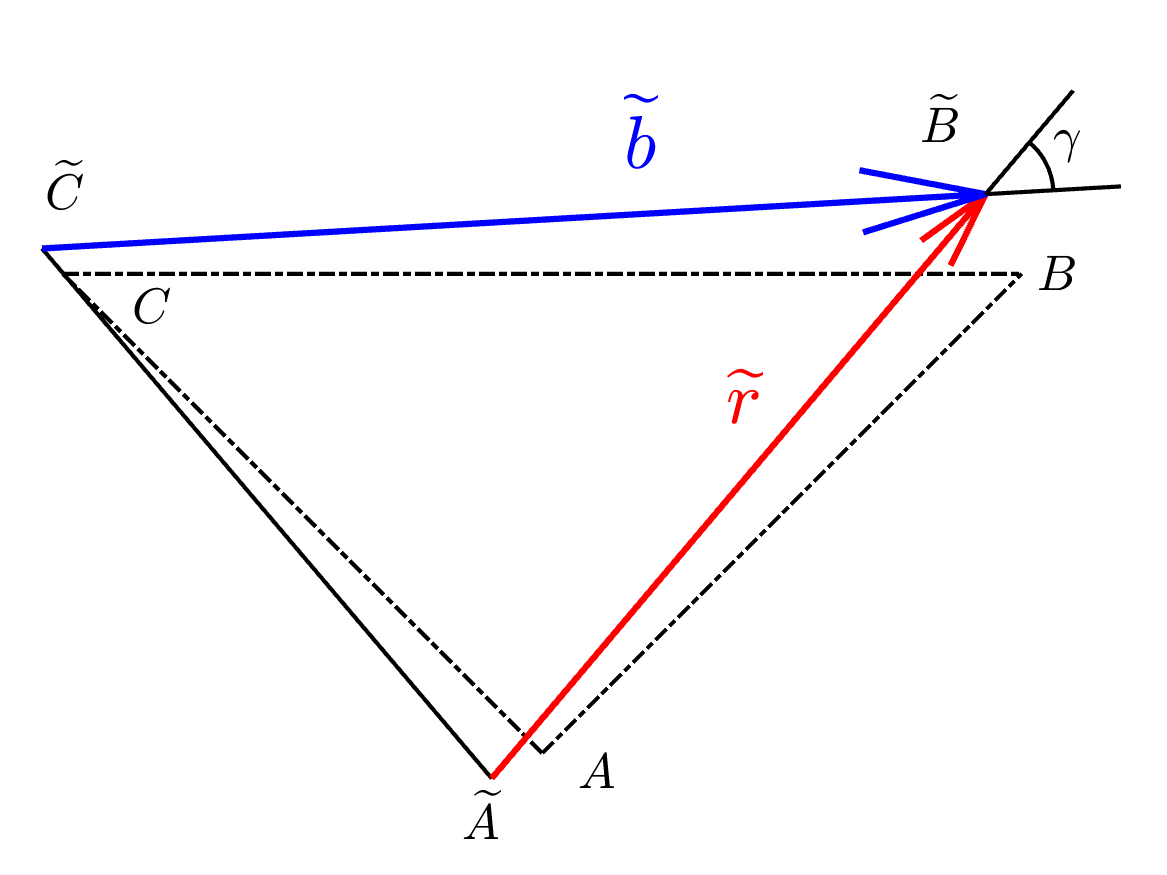}
		\caption{The triangle with nodes $A,B,C$ in the reference state and its deformed state with nodes $\widetilde{A}, \widetilde{B}, \widetilde{C}$. The angles in the reference state and deformed states are $\angle ABC = \alpha$ and $\angle \widetilde{A}\widetilde{B}\widetilde{C} = \gamma$.}
		\label{fig:triangle-1}
	\end{figure}

	This appendix gives the proof of Proposition \ref{prop:iso-triangle}. We recall our notation: the vectors $\widetilde{b}$ and $\widetilde{r}$ are the images of vectors $\overrightarrow{CB}$ and $\overrightarrow{AB}$ (see Figure \ref{fig:triangle-1} as an illustration), i.e.
	\begin{equation*}
		\widetilde{b} = \overrightarrow{\widetilde{C}\widetilde{B}}, \qquad \widetilde{r} = \overrightarrow{\widetilde{A}\widetilde{B}},
	\end{equation*}
	where $\widetilde{A} = u(A), \widetilde{B} = u(B), \widetilde{C} = u(C)$ are the image of nodes. The spring energy $E^\text{spr}(u)$ is
	\begin{equation*}
		E^\text{spr}(u) = \Big(|\widetilde{b}|-1\Big)^2 + \Big(|\widetilde{r}|-1\Big)^2 + \Big(|\widetilde{r} -\widetilde{b}|-\sqrt{2-2\cos\alpha}\Big)^2,
	\end{equation*}
	where $\sqrt{2-2\cos\alpha}$ is the length of edge $AC$, computed using the Law of Cosines. We shall show that
	\begin{equation}\label{eqn:app-iso-bd}
		E^\text{spr}(u) \geq \begin{cases}
			c |\widetilde{r} - R_\alpha \widetilde{b}|^2, & \det(\nabla u) \geq 0,\\
			c |\widetilde{r} - R_{-\alpha} \widetilde{b}|^2, & \det(\nabla u) < 0,\\
		\end{cases}
	\end{equation}
	for some constant $c$ independent of $u$.
	
	To prove \eqref{eqn:app-iso-bd}, we consider two cases: 
	\begin{enumerate}
		\item[(1)] when the energy $E^\text{spr}(u)$ is small (in the sense 
		that $E^\text{spr}(u) \leq \epsilon_0^2$ for some small $\epsilon_0$), 
		and
		\item[(2)] when $E^\text{spr}(u)$ is large (in the sense that 
		$E^\text{spr}(u) \geq \epsilon^2_0$) . 
	\end{enumerate}
	We start by observing that the desired result \eqref{eqn:app-iso-bd} is very easy in case (2). Indeed, since $E^\text{spr}(u)$ is the sum of the spring energies of the edges of $\Delta ABC$, the length of each deformed spring is controlled by  $\sqrt{E^\text{spr}(u)}$:
	\begin{equation}\label{eqn:app-spring-bd}
		\begin{aligned}
			& \big||\widetilde{b}| - 1\big| \leq \sqrt{E^\text{spr}(u)} \qquad \Rightarrow \qquad 1 -\sqrt{E^\text{spr}(u)} \leq |\widetilde{b}| \leq 1 + \sqrt{E^\text{spr}(u)}\\
			& \big||\widetilde{r} | - 1\big| \leq \sqrt{E^\text{spr}(u)} \qquad  \Rightarrow \qquad1 + \sqrt{E^\text{spr}(u)} \leq |\widetilde{r} | \leq 1 + \sqrt{E^\text{spr}(u)},
		\end{aligned}
	\end{equation}
	Therefore, when $E^\text{spr}(u) \geq \epsilon_0^2$, we have $\frac{2}{\epsilon_0} \sqrt{E^\text{spr}(u)} \geq 2$ and
	\begin{equation}\label{eqn:large-iso-bd}
		\Big|\widetilde{r} - R_{\alpha} \widetilde{b}\Big| \leq |\widetilde{r}| + |\widetilde{b}| \leq 2 + 2 \sqrt{E^\text{spr}(u)} \leq \Big(2 + \frac{2}{\epsilon_0}\Big) \sqrt{E^\text{spr}(u)}.
	\end{equation}
	It is worth mentioning that \eqref{eqn:large-iso-bd} holds even when the orientation of the triangle is reversed. Similarly, we have that when $E^\text{spr}(u) \geq \epsilon_0^2$,
	\begin{equation}\label{eqn:large-iso-bd-flip}
		\Big|\widetilde{r} - R_{-\alpha} \widetilde{b}\Big| \leq |\widetilde{r}| + |\widetilde{b}| \leq 2 + 2 \sqrt{E^\text{spr}(u)} \leq \Big(2 + \frac{2}{\epsilon_0}\Big) \sqrt{E(u)},
	\end{equation}
	even when the orientation of the triangle is preserved. Thus, the lower bound \eqref{eqn:app-iso-bd} holds when $E^\text{spr}(u) \geq \epsilon_0^2$.
	
	The rest of this appendix addresses case (1), when 
	$E^\text{spr}(u) \leq \epsilon_0^2$ . (Some smallness conditions on $\epsilon_0$ will emerge in the course of the argument.) It is intuitively clear that the deformed vector $\widetilde{r}$ should be close to 
	$R_{\alpha} \widetilde{b}$ if the orientation of the triangle is preserved, and  close to $R_{-\alpha} \widetilde{b}$ if the orientation is reversed; our task is to make this quantitative. Our strategy is simple: we first show that the angle $\gamma$ between the two vectors $\widetilde{b}$ and $\widetilde{r}$ (see Figure \ref{fig:triangle-1}) deviates from $\pm \alpha$ by a small amount that can be controlled by $\sqrt{E^\text{spr}(u)}$, specifically we shall show that
	\begin{equation}\label{eqn:angle-relation}
		\begin{aligned}
			& \Big|\gamma - \alpha\Big| \leq c \sqrt{E^\text{spr}(u)}, \qquad \gamma \geq 0, \\
			& \Big|\gamma + \alpha \Big| \leq c \sqrt{E^\text{spr}(u)}, \qquad \gamma < 0, 
		\end{aligned}
	\end{equation}
	for some constant $c$ that does not depend on $u$; then we give a relatively simple argument to show that \eqref{eqn:angle-relation} implies the desired result \eqref{eqn:app-iso-bd}.
	
	To show \eqref{eqn:angle-relation}, we use the Law of Cosines for the angle $\gamma$. For simplicity, we shall write $x,y,z$ for the lengths of the deformed edges: 
	\begin{equation*}
		|\widetilde{b}| = x, \qquad |\widetilde{r}|= y, \qquad|\widetilde{r} - \widetilde{b}| = z .
	\end{equation*}
	The Law of Cosines says
	\begin{equation}\label{eqn:law-cosine}
		\cos \gamma = \frac{x^2 + y^2 - z^2}{2xy}.
	\end{equation}
	We claim that when $\sqrt{E^\text{spr}(u)} \leq \frac{1}{2}$, the following bounds hold
	\begin{equation}\label{eqn:cos-rule-bd}
		|x^2 + y^2 - z^2 -2\cos\alpha| \leq 9 \sqrt{E^\text{spr}(u)}, \qquad
		|xy-1| \leq 3 \sqrt{E^\text{spr}(u)}.
	\end{equation}
	The proof of \eqref{eqn:cos-rule-bd} is a direct calculation, which we now discuss. We know that $x,y,z$ are controlled by the spring energy:
	\begin{equation*}
		|x-1| \leq \sqrt{E^\text{spr}(u)}, \qquad |y-1|\leq \sqrt{E^\text{spr}(u)}, \qquad |z-\sqrt{2-2\cos\alpha}| \leq \sqrt{E^\text{spr}(u)}.
	\end{equation*}
	To simplify our proof, we introduce names for the deviations:
	\begin{equation*}
		\epsilon_1 = x-1, \qquad \epsilon_2 = y-1,\qquad \epsilon_3 = z-\sqrt{2-2\cos\alpha}.
	\end{equation*}
	Using $\sqrt{E^\text{spr}(u)} \leq \frac{1}{2}$, we obtain that (i) $|\epsilon_i| \leq \sqrt{E^\text{spr}(u)} \leq \frac{1}{2}$ for $i=1,2,3$; (ii) $\epsilon_i^2 \leq \frac{1}{2} |\epsilon_i|$ for $i=1,2,3$; and (iii) $|\epsilon_i \epsilon_j| \leq \frac{1}{4} (|\epsilon_i| + |\epsilon_j|)$ for all $i,j = 1,2,3$.
	By writing $x^2 + y^2 - z^2$ in terms of $\epsilon_1, \epsilon_2, \epsilon_3$, we have
	\begin{equation*}
		\begin{aligned}
			x^2 + y^2 - z^2 &= (1+\epsilon_1)^2 + (1+\epsilon_2)^2 - (\sqrt{2-2\cos\alpha}+\epsilon_3)^2\\
			&= 2\cos\alpha + 2(\epsilon_1 + \epsilon_2) -2\epsilon_3\sqrt{2-2\cos\alpha} + \epsilon_1^2 + \epsilon_2^2 - \epsilon_3^2.
		\end{aligned}
	\end{equation*}
	Using the bounds mentioned in (i) and (ii), we obtain
	\begin{equation*}
		\begin{aligned}
			2\cos \alpha  -8\sqrt{E^\text{spr}(u)} - \frac{1}{2} \sqrt{E^\text{spr}(u)} \leq x^2 + y^2 - z^2  & \leq 2\cos\alpha + 8\sqrt{E^\text{spr}(u)} + \sqrt{E^\text{spr}(u)},
		\end{aligned}
	\end{equation*}
	where $ \pm 8\sqrt{E^\text{spr}(u)}$ serves as an upper and lower bound for the linear terms in $\epsilon_1, \epsilon_2, \epsilon_3$. Therefore, we obtain $|x^2 + y^2 - z^2 - 2\cos\alpha| \leq 9 \sqrt{E^\text{spr}(u)}$. Arguing similarly for the $xy$ term, we write it as $xy = (1+\epsilon_1) (1+\epsilon_2) = 1 + (\epsilon_1 + \epsilon_2) + \epsilon_1 \epsilon_2$. Using (i) and (iii), we obtain the bound
	\begin{equation*}
		1 - 2\sqrt{E^\text{spr}(u)} - \frac{1}{2} \sqrt{E^\text{spr}(u)}\leq xy \leq 1 + 2\sqrt{E^\text{spr}(u)} + \frac{1}{2} \sqrt{E^\text{spr}(u)},
	\end{equation*}
	which leads to the desired bound $|xy-1| \leq 3\sqrt{E^\text{spr}(u)}$.
	
	Returning to the formula  \eqref{eqn:law-cosine} for $\cos\gamma$, by choosing $(1-3\sqrt{E^\text{spr}(u)})(1+6\sqrt{E^\text{spr}(u)}) \geq 1$ (equivalently, choosing $\sqrt{E^\text{spr}(u)} \leq 1/6$), we can bound $\cos \gamma$ above and below as follows:
	\begin{equation*}
		\begin{aligned}
			2 \cos \gamma & \leq \frac{2 \cos \alpha +9\sqrt{E^\text{spr}(u)}}{1-3\sqrt{E^\text{spr}(u)}} \leq (2\cos \alpha +9\sqrt{E^\text{spr}(u)}) (1+6\sqrt{E^\text{spr}(u)}) \leq 2\cos \alpha + c_1 \sqrt{E^\text{spr}(u)},\\
			2 \cos \gamma & \geq \frac{2\cos \alpha-9\sqrt{E^\text{spr}(u)}}{1+3\sqrt{E^\text{spr}(u)}} \geq (1-3\sqrt{E^\text{spr}(u)}) (2\cos \alpha-9\sqrt{E^\text{spr}(u)}) \geq 2 \cos \alpha-c_1 \sqrt{E^\text{spr}(u)},
		\end{aligned}
	\end{equation*}
	for some constant $c_1$ independent of $\sqrt{E^\text{spr}(u)}$. (Here we have used that $\frac{1}{1-t} \leq 1+2t$ and $\frac{1}{1+t} \geq1 - t$ when $0 \leq t \leq \frac{1}{2}$.) Thus, we have shown
	\begin{equation}\label{eqn:cos-rule}
		\Big|\cos\gamma - \cos\alpha\Big| \leq \frac{c_1}{2} \sqrt{E^\text{spr}(u)}.
	\end{equation}
	Since the function $\arccos: \cos \gamma \rightarrow \gamma \in [0,\pi]$ is Lipschitz near $\cos\alpha$, by restricting $\cos\gamma$ to be in a neighborhood of $\cos \alpha$ (equivalently, by requiring that the term $\frac{c_1}{2} \sqrt{E^\text{spr}(u)}$ be small enough), we can find a Lipschitz constant $c_2$ such that
	\begin{equation}\label{eqn:acos-lipschitz}
		\begin{aligned}
			&|\gamma - \alpha| \leq c_2 \Big|\cos\gamma - \cos \alpha\Big|, & & \text{ when } \gamma \geq 0,\\
			&|\gamma + \alpha| \leq c_2 \Big|\cos\gamma - \cos 
			\alpha\Big|, & & \text{ when } \gamma < 0.\\
		\end{aligned}
	\end{equation}
	Therefore, by combining \eqref{eqn:cos-rule} with \eqref{eqn:acos-lipschitz}, we obtain the desired lower bounds \eqref{eqn:angle-relation}.
	
	\begin{remark}
		Though we did not choose the energy threshold $\epsilon_0$ explicitly in the proof, our argument requires only that the energy $E^\text{spr}(u)$ be small enough that (i) $\sqrt{E^\text{spr}(u)} \leq \frac{1}{6}$ and (ii) the Lipschitz constant $c_2$ can be found.
	\end{remark}
	
	Lastly, we use the angle bound \eqref{eqn:angle-relation} to prove the desired bound \eqref{eqn:app-iso-bd}. We begin by noticing that the sign of $\gamma$ aligns with the orientation of the deformed triangle. To see why, we observe that the deformation $D u$ on the triangle $\Delta ABC$ is
	\begin{equation*}
		D u = \begin{pmatrix}
			\overrightarrow{\widetilde{C} \widetilde{B}} & \overrightarrow{\widetilde{A} \widetilde{B}}
		\end{pmatrix} \begin{pmatrix}
			\overrightarrow{CB} & \overrightarrow{AB}
		\end{pmatrix}^{-1} .
	\end{equation*}
	Since $\det\begin{pmatrix}
		\overrightarrow{CB} & \overrightarrow{AB}
	\end{pmatrix} > 0$, the signs of $\det(D u)$ and $\det\begin{pmatrix}
		\overrightarrow{\widetilde{C} \widetilde{B}} & \overrightarrow{\widetilde{A} \widetilde{B}}
	\end{pmatrix}$ are the same. We also have the following relationship between $\sin \gamma$ and $\det\begin{pmatrix}
		\overrightarrow{\widetilde{C} \widetilde{B}} & \overrightarrow{\widetilde{A} \widetilde{B}}
	\end{pmatrix}$:
	\begin{equation*}
		\bigg|\overrightarrow{\widetilde{C} \widetilde{B}}\bigg| \bigg|\overrightarrow{\widetilde{A} \widetilde{B}}\bigg| \sin \gamma = \det \begin{pmatrix}
			\overrightarrow{\widetilde{C} \widetilde{B}} & \overrightarrow{\widetilde{A} \widetilde{B}}
		\end{pmatrix}.
	\end{equation*}
	Therefore, the sign of $\det(D u)$ is the same as the sign of $\gamma$. 
	
	Proceeding now toward the desired lower bound \eqref{eqn:app-iso-bd}, we observe that $R_\gamma \widetilde{b}$ is parallel to $\widetilde{r}$, i.e.
	\begin{equation}\label{eqn:rotation}
		R_{\gamma} \widetilde{b} = \frac{|\widetilde{b}|}{|\widetilde{r}|} \widetilde{r}.
	\end{equation}
	Using \eqref{eqn:app-spring-bd} and $\sqrt{E^\text{spr}(u)} \leq \frac{1}{2}$, we obtain bounds on the ratio $\frac{|\widetilde{b}|}{|\widetilde{r}|}$ :
	\begin{equation*}
		\begin{aligned}
			\frac{|\widetilde{b}|}{|\widetilde{r}|} &\leq \frac{1+\sqrt{E^\text{spr}(u)}}{1-\sqrt{E^\text{spr}(u)}} \leq (1+\sqrt{E^\text{spr}(u)}) (1+2\sqrt{E^\text{spr}(u)}) \leq 1 + 4\sqrt{E^\text{spr}(u)}, \\
			\frac{|\widetilde{b}|}{|\widetilde{r}|} &\geq \frac{1-\sqrt{E^\text{spr}(u)}}{1+\sqrt{E^\text{spr}(u)}} \geq (1-\sqrt{E^\text{spr}(u)})^2 \geq 1 - 2\sqrt{E^\text{spr}(u)}.
		\end{aligned}
	\end{equation*}
	Consequently, we obtain
	\begin{equation}\label{eqn:upper-bd-rot}
		\Big|\frac{|\widetilde{b}|}{|\widetilde{r}|} - 1\Big| \leq 4 \sqrt{E^\text{spr}(u)}.
	\end{equation}
	Thus, combining \eqref{eqn:rotation} with \eqref{eqn:upper-bd-rot}, we get
	\begin{equation}\label{eqn:upper-bd-rot-2}
		\Big|R_\gamma \widetilde{b} - \widetilde{r}\Big| = |\widetilde{r}| \Big|\frac{|\widetilde{b}|}{|\widetilde{r}|} - 1\Big| \leq 4 \sqrt{E^\text{spr}(u)} \Big(1+\sqrt{E^\text{spr}(u)}\Big)  \leq 4\sqrt{E^\text{spr}(u)} + 4E^\text{spr}(u) \leq 6\sqrt{E^\text{spr}(u)}.
	\end{equation}
	Finally, to provide a bound for $R_{\alpha} \widetilde{b} - \widetilde{r}$, we use the fact that the rotation matrix $R_\gamma$ is Lipschitz as a function of $\gamma$, i.e. there exists a constant $c_3 > 0$ such that
	\begin{equation}\label{eqn:rot-mat-lipschitz}
		\begin{aligned}
			& \Big|R_\alpha - R_{\gamma}\Big| \leq c_3 |\alpha - \gamma|, & & \text{when } \gamma \geq 0 \quad (\det(D u) \geq 0),\\
			& \Big|R_{-\alpha} - R_{\gamma}\Big| \leq c_3 |\alpha + \gamma|, & & \text{when } \gamma < 0 \quad (\det(D u) < 0).
		\end{aligned}
	\end{equation}
	Using \eqref{eqn:angle-relation} and \eqref{eqn:upper-bd-rot-2}, we obtain the desired lower bound \eqref{eqn:app-iso-bd} when $\det(Du)\geq 0$
	\begin{equation*}
		\Big|R_\alpha \widetilde{b} - \widetilde{r}\Big| \leq \Big|R_{\gamma} \widetilde{b} - \widetilde{r}\Big| + \Big|R_\alpha - R_{\gamma}\Big||\widetilde{b}|\leq 6\sqrt{E^\text{spr}(u)} + c c_3 \sqrt{E^\text{spr}(u)} \Big(1 + \sqrt{E^\text{spr}(u)}\Big) \leq C_1 \sqrt{E^\text{spr}(u)},
	\end{equation*}
	for some constant $C_1$. Similarly, we obtain the following bound when $\det(D u) < 0$
	\begin{equation*}
		\Big|R_{-\alpha} \widetilde{b} - \widetilde{r}\Big| \leq C_1 \sqrt{E(u)},
	\end{equation*}
	which completes the proof of Proposition \ref{prop:iso-triangle}.
	
	\setcounter{equation}{0}
	\subsection{A non-periodic mechanism} \label{subsec:non-periodic}
\begin{figure}[!htb]
		\centering
		\subfloat[]{\includegraphics[width=0.7\linewidth]{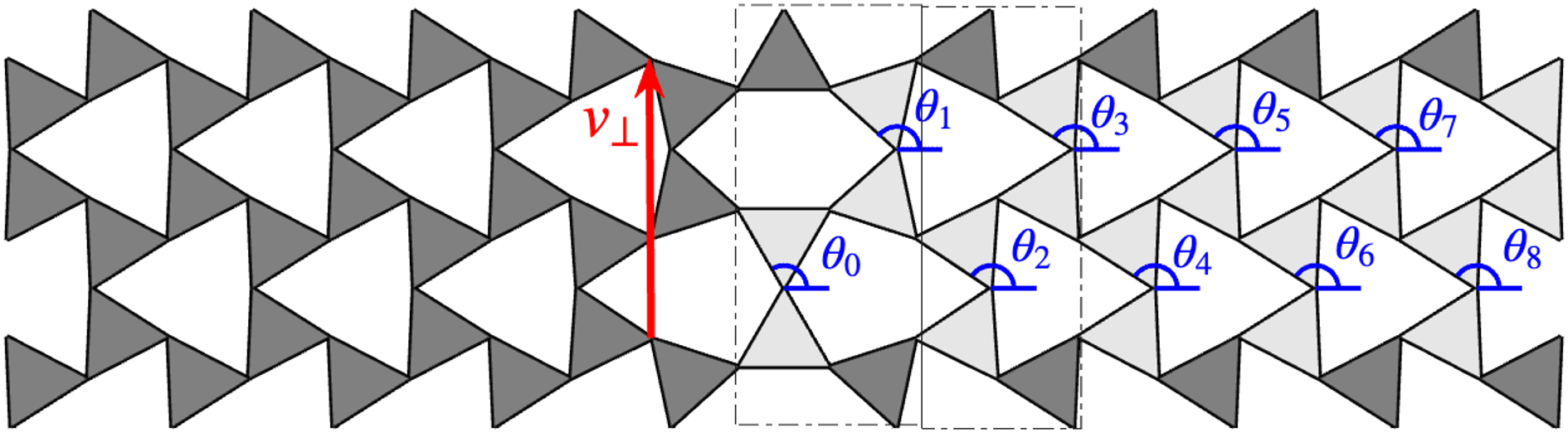}}\hfil
		\subfloat[]{\includegraphics[width=0.2\linewidth]{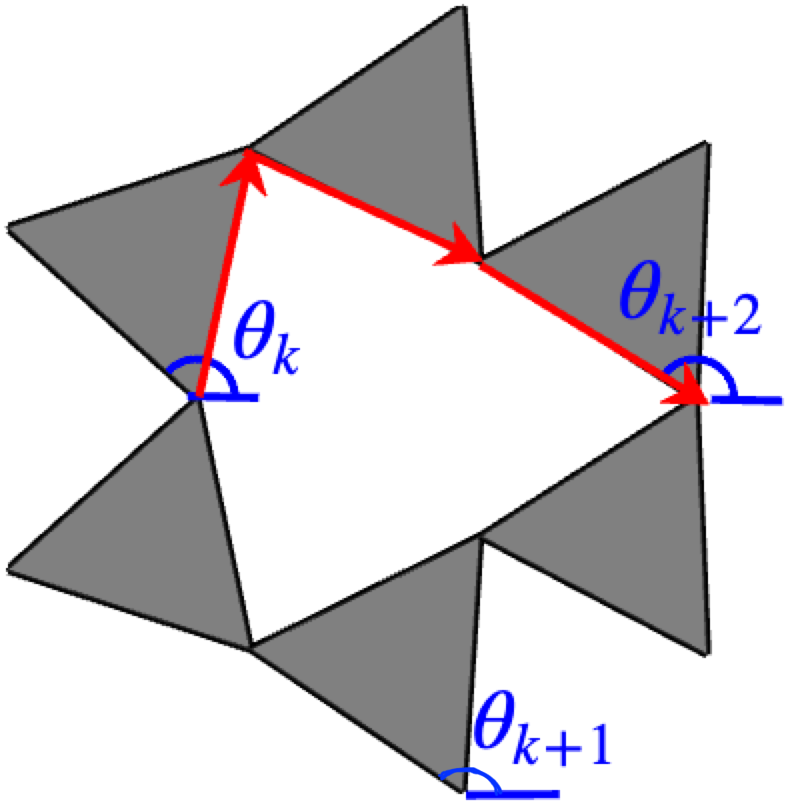}}
		\caption{The non-periodic mechanism: (a) the deformation is periodic in the y-direction with translation vector $v_\perp$ and is symmetric about the pair of triangles marked by the angle $\theta_0$; (b) the relationship between the angles $\theta_k, \theta_{k+1}, \theta_{k+2}$.}
		\label{fig:domain-wall-unit}
	\end{figure}

	We give the details of the non-periodic mechanism shown in Figure \ref{fig:kagome-mechanisms}c. This mechanism is periodic in the vertical direction and symmetric in the horizontal direction (see Figure \ref{fig:domain-wall-unit}a). It features a central domain wall, with the left and right regions being mirror images of each other.
	
	The non-periodic mechanism is built from vertical layers of twisted pairs of equilateral triangles. Due to its periodicity in the $y$-direction and symmetry in the $x$-direction, it suffices to describe only half of the unit cell -- the right side of the horizontal slice in Figure \ref{fig:domain-wall-unit}a, shown in lighter gray. In fact, this half horizontal unit cell can be divided into columns, each containing four triangles: for $i = 1, 2, \dots$, the $i$-th column consists of two pairs of twisted triangles, with twist angles $\theta_{2i-2}$ and $\theta_{2i-1}$. The mechanism is fully specified once all these angles are determined.
	
	The sequence of angles $\theta_0, \theta_1, \dots$ are not arbitrary. To form the domain wall in the middle, we need $\theta_0 = \frac{2\pi}{3}$. Once $\theta_1$ is specified, all remaining angles are uniquely determined since they satisfy the following relationship
	\begin{equation}\label{eqn:angle-relation-energy-free-wall}
		\sin \left(\theta_k - \frac{\pi}{3}\right) = \sin \left(\theta_{k+1} - \frac{\pi}{3}\right) - \sin \theta_{k+1} + \sin \theta_{k+2}, \qquad \forall \: k = 0,1,2,\dots.
	\end{equation}
	To understand why this relationship holds, note that the hexagon formed by the twisted triangles with angles $\theta_{k}, \theta_{k+1}, \theta_{k+2}$ is closed. Consequently, the marked vectors in Figure \ref{fig:domain-wall-unit}b must sum to a horizontal vector, and the requirement that the y-component vanishes yields \eqref{eqn:angle-relation-energy-free-wall}. We also note that if all angles are equal to an arbitrary value, i.e. $\theta_0 = \theta_1 = \theta_2 = \dots$ (without the requirement $\theta_0 = 2\pi/3$), the relationship \eqref{eqn:angle-relation-energy-free-wall} is automatically satisfied. This configuration corresponds to the one-periodic mechanism shown in Figure \ref{fig:kagome-mechanisms}a.
	
	Lastly, it is worth noting that once $\theta_1$ is specified, the sequence $\theta_0, \theta_1, \theta_2, \dots$ converges rapidly. Numerically, for most $\theta_1$ values between $2\pi/3$ and $\pi$, the sequence has nearly converged by $\theta_{10}$. The far-field limit is, of course, the one-periodic mechanism (but with different values of its parameter at the far left and the far right). The two far-field patterns are reflections of one another; in particular, they achieve the same macroscopic compression.
	
	\bibliography{ref}
	\bibliographystyle{amsalpha}

\end{document}